\documentclass[12pt]{amsart}
\usepackage{amsfonts,amssymb,amscd,amsmath}
\usepackage{latexsym,amsbsy,bbold,stmaryrd}
\usepackage{tocvsec2}
\usepackage{mathrsfs}

\usepackage[all]{xy}

\usepackage[
dvips=true, 
hyperindex=true,pagebackref=true,bookmarks=true,
colorlinks=true,linkcolor=blue,citecolor=red]
{hyperref}
\def\~{{\rm --}}

\settocdepth{subsubsection}

\renewcommand{\tilde}{\widetilde}
\renewcommand{\hat}{\widehat}

\newcommand{\BZ}{{\mathbb Z}}

\newcommand{\cT}{\hbox{$\mathcal T$}}

\newcommand{\Z}{{\mathbb Z}}

\newcommand{\N}{{\mathbb N}}
\newcommand{\C}{{\mathbb C}}

\def\HH{\mbox{${\mathcal H}$\kern-5.2pt${\mathcal H}$}}

\newtheorem{theorem}{Theorem}[section]

\newtheorem{proposition}[theorem]{Proposition}
\newtheorem{definition}[theorem]{Definition}
\newtheorem{lemma}[theorem]{Lemma}

\newtheorem{theorem }{Theorem}[section]
\newtheorem{maintheorem }[theorem]{Main Theorem}
\newtheorem{proposition }[theorem]{Proposition}
\newtheorem{definition }[theorem]{Definition}
\newtheorem{lemma }[theorem]{Lemma}
\newtheorem{corollary }[theorem]{Corollary}
\newtheorem{notation }[theorem]{Notation}
\newtheorem{remark }[theorem]{Remark}
\newtheorem{example }[theorem]{Example}

\newtheorem{ maintheorem }[theorem]{Main Theorem}
\newtheorem{ theorem}{Theorem}[section]
\newtheorem{ proposition}[theorem]{Proposition}
\newtheorem{ definition}[theorem]{Definition}
\newtheorem{ lemma}[theorem]{Lemma}
\newtheorem{ corollary}[theorem]{Corollary}
\newtheorem{ notation}[theorem]{Notation}
\newtheorem{ remark}[theorem]{Remark}
\newtheorem{ example}[theorem]{Example}

 \newcommand{\rmk}{{\bf Comment.\ }}

\def\for{\  \hbox{ for } \ }

\def\where{\  \hbox{ where } \ }
\def\and{\  \hbox{ and } \ }
\def\and{\  \hbox{ and } \ }

\def\equal{\stackrel{\,\mathbf{def}}{= \kern-3pt =}}

\def\la{\lambda}
\def\La{\Lambda}
\def\om{\omega}
\def\Om{\Omega}

\def\al{\alpha}
\def\be{\beta}
\def\ga{\gamma}
\def\ep{\epsilon}

\def\de{\delta}

\def\si{\sigma}

\def\Ga{\Gamma}

\def\vph{\varphi}

\def\vep{\varepsilon}

\def\tga{\tilde{\gamma}}

\def\B{\mathbf{B}}

\def\S{\mathbf{S}}
\def\F{\mathbf{F}}

\def\0{\mathbf{0}}

\def\f{\mathcal{F}}
\def\çF{\mathcal{F}}

\def\t{\mathcal{T}}

\def\l{\mathcal{L}}

\def\p{\mathcal{P}}

\def\a{\mathcal{A}}
\def\h{\mathcal{H}}

\def\y{\mathcal{Y}}
\def\e{\mathcal{E}}
\def\v{\mathcal{V}}

\def\x{\mathcal{X}}
\def\s{\mathcal{S}}

\def\w{\mathcal{W}}

\def\lan{\langle}

\def\ran{\rangle}

\newcommand{\sq}{\phantom{1}\hfill$\qed$}

\newcommand{\lr}{\langle}
\newcommand{\rr}{\rangle}

\newcommand{\sgn}{\mbox{sgn}}

\def\HH{\mathfrak{H}}

\def\HH{\hbox{${\mathcal H}$\kern-5.2pt${\mathcal H}$}}

\def\#{\sharp}

\begin{document}
\newcommand{\comment}[1]{}

\begin{abstract}
This paper is devoted to the theory of nil-DAHA for
the root system $A_1$ and its applications to symmetric
and nonsymmetric (spinor) global $q$-Whittaker functions.
These functions integrate the $q$-Toda eigenvalue problem
and its Dunkl-type nonsymmetric version.

The global symmetric function can be interpreted as
the generating function of the Demazure characters
for dominant weights, which describe the algebraic-geometric
properties of the corresponding affine Schubert varieties.
Its Harish-Chandra-type asymptotic expansion appeared directly 
related to the solution of the $q$-Toda eigenvalue problem obtained 
by Givental and Lee in the quantum $K$-theory of flag varieties.
It provides an exact mathematical relation between the corresponding 
physics $A$-type and $B$-type models.

The spinor global functions extend the symmetric ones to the case of
all Demazure characters (not only those for the dominant weights);
the corresponding Gromov-Witten theory is not known.
The main result of the paper is a complete algebraic theory of these
functions in terms of the induced modules of the core subalgebra of
nil-DAHA.  It is the first instance of the DAHA theory of
canonical-crystal bases, quite non-trivial even for $A_1$.
\end{abstract}

\title[Nil-DAHA and Whittaker functions]
{One-dimensional nil-DAHA and Whittaker functions}
\author[Ivan Cherednik]{Ivan Cherednik $^\dag$}
\author[Daniel Orr]{Daniel Orr}

\address[I. Cherednik]{Department of Mathematics, UNC
Chapel Hill, North Carolina 27599, USA\\
chered@email.unc.edu}
\address[D. Orr]{Department of Mathematics, UNC
Chapel Hill, North Carolina 27599, USA\\
danorr@email.unc.edu}

\def\bysame{{\bf --- }}
\def\~{{\bf --}}
\renewcommand{\tilde}{\widetilde}
\renewcommand{\hat}{\widehat}
\newcommand{\dagx}{\hbox{\tiny\mathversion{bold}$\dag$}}
\newcommand{\ddagx}{\hbox{\tiny\mathversion{bold}$\ddag$}}

\thanks{$^\dag$  
\ \ \ Partially supported by NSF grant
DMS--1101535.\\
This paper will appear in {\em Transformation Groups} in two parts.}
\maketitle
\renewcommand{\baselinestretch}{1.2}
{
\smallskip
\tableofcontents
\smallskip
}
\renewcommand{\baselinestretch}{1.2}
\vfill\eject

\renewcommand{\natural}{\wr}
\renewcommand{\baselinestretch}{1.2}
\setcounter{section}{-1}
\setcounter{equation}{0}
\section{\sc Introduction}

This paper is devoted to the theory of nil-DAHA $\overline{\HH}$ for
the root system $A_1$ and its applications to the symmetric and 
nonsymmetric $q$\~Whittaker functions, integrating the $q$\~Toda 
eigenvalue problem and its Dunkl counterpart introduced in \cite{C102}.

In the symmetric case, we establish an exact relation between the
generating functions of the $A$\~type model ($K$\~theoretic
Gromov-Witten invariants of flag varieties) and the $B$\~type model 
(sheaf cohomology of affine Schubert varieties) in the case of $SL(2)$.


\medskip
{\bf Major connections.}
In greater detail, the global {\em symmetric $q$\~Whittaker
functions}, defined in \cite{ChW} for any reduced root
systems, can be interpreted as generating functions of 
level one Demazure characters (equivalently, $q$\~Hermite 
polynomials) for dominant weights; these characters describe 
the spaces of global sections of the canonical level one line bundle 
restricted to the corresponding affine Schubert subvarieties of 
the affine flag variety. From this perspective, the symmetric
$q$\~Whittaker functions are examples of the generating 
functions of $B$\~type models in physics.  

On the other hand, the Harish-Chandra asymptotic decompositions of 
these functions from \cite{ChW} are directly related to the 
solutions of the $q$\~Toda eigenvalue problem found in \cite{GiL}.
These solutions are important instances of generating functions for
$A$\~models, as they describe $K$\~theoretic genus zero Gromov-Witten
invariants. The Appendix, Section \ref{sec:qKth}, contains more 
details on this connection.

Thus we arrive at an exact mathematical relation between these
two theories. It is described in this paper (mathematically)
in full detail for type $A_1$. This construction is a explicit example 
of the Whittaker limit of the $q,t$\~generalization of Harish-Chandra 
theory of asymptotic expansions started in \cite{ChW} and recently 
finalized in \cite{Sto2}.
\medskip

{\bf Nonsymmetric functions.}
One of the main results of \cite{C102} is the construction of
the {\em spinor Dunkl operators} for the $q$\~Toda eigenvalue problem
in type $A_1$. It was a surprising development, since the
$q$\~Toda operators are already nonsymmetric and the standard
approach based on the symmetrization of the Dunkl operators
was not expected to work.

A related development of \cite{C102} is the
definition of the {\em nonsymmetric $q$\~Whittaker function},
which integrates the eigenvalue problem for the spinor Dunkl operators.
A variant of the Ruijsenaars-Etingof limiting
procedure \cite{Rui, Et1} was used to obtain it from the
nonsymmetric global $q,t$\~spherical function defined in \cite{C5}.

In the $B$\~type model interpretation,
the nonsymmetric $q$\~Whittaker function extends its
symmetric counterpart to the case of all Demazure characters (not only
those for the dominant weights) and arbitrary Schubert subvarieties
of the affine flag variety. The ``nonsymmetric'' counterpart of the
$K$\~theoretic Gromov-Witten theory is not known.

The natural expectations are that the nonsymmetric
$q$\~Whittaker functions will impact the classical real and
$p$\~adic theories of Whittaker functions, as well
as the recent Whittaker theories based on the Kac-Moody
algebras and Quantum Groups.
See \cite{C102}, \cite{GiL}, \cite{GLO}, \cite{FJM},
and \cite{BF} for this and other known and expected applications of
$q$\~Whittaker functions.
In this connection, let us also mention that the nil-DAHA is
directly related to the constructions of \cite{KK} in the case of
the affine flag variety.

The Matsumoto $p$\~adic spherical functions,
the simplest example of the nonsymmetric theory, demonstrate
potential power of such development.
DAHA already resulted in the theory of real counterparts of
the Matsumoto functions, which are the nonsymmetric Macdonald
polynomials, and unified the real and $p$\~adic theories in
one $q,t$\~theory. Now it can be done for Whittaker functions!
\medskip

{\bf Core subalgebra.}
The main result of this paper is an entirely algebraic theory of
the spinor Dunkl operators and nonsymmetric Whittaker function for 
$A_1$ in terms of the induced modules of a new object, {\em the core
subalgebra of nil-DAHA}.

The core subalgebra is essentially the span
of the basic creation and annihilation operators for the
nonsymmetric $q$\~Hermite polynomials and the $Y$\~operators
(which diagonalize these polynomials).
Its existence {\em as a subalgebra of} $\overline{\HH}$
is an important special feature of the nil-case.
Conceptually, {\em the nonsymmetric Whittaker function is the
reproducing kernel of the transform sending the
nonsymmetric $q$\~Hermite polynomials
to the corresponding creation operators}. See Theorem
\ref{EPOLTILDE},
a simple (at least in $A_1$) but fundamental fact; this connection
seems to be the main algebraic message of this paper.

It is expected that the core subalgebra will play the
prime role in the theory of {\em canonical-crystal} bases of DAHA.
Its graded algebra is a nil-nil variant of the extended
Weyl algebra (a non-commutative torus extended by the
Weyl group $W$), which is quite interesting in its own right.

This paper contains a systematic analysis of modules for the
core subalgebra induced from one-dimensional representations
of {\em core affine Hecke subalgebras}. Some of these modules
totally collapse, i.e., remain one-dimensional.
The induced modules that are of the right size
(infinite-dimensional, of course) can be extended
to modules over the whole $\overline{\HH}$. This is
one of the key results of the paper, Theorem \ref{TILDEINTER}.
The core-induced modules have a natural filtration of submodules;
the corresponding irreducible constituents are one-dimensional
(for $A_1$) and can be calculated in full detail.

A key point is that the core subalgebra contains all 
intertwining operators necessary to decompose
such and similar induced modules. It is an important change from
the theory of intertwiners for AHA and DAHA. Actually, it is a new
and interesting version of the technique of intertwiners in the
nil-case.

The core subalgebra solves the problem of finding
an induced interpretation of the spinor-polynomial module.
This module was discovered in \cite{C102} {\em analytically}, via
reading the coefficients of the spinor $q$\~Whittaker function,
obtained there as a limit of the global spherical nonsymmetric
function. Now we can introduce it entirely {\em algebraically};
it is isomorphic to an induced module of the core
subalgebra upon its extension to $\overline{\HH}$. Correspondingly,
the spinor $q$\~Whittaker function can be now defined and
calculated within $\overline{\HH}$, i.e., without any reference
to the $q,t$\~theory. See Theorems \ref{TILDESPIN} and \ref{TILDEOM}.
\medskip

{\bf Analytic aspects.}
There nonsymmetric global spherical functions defined in
\cite{C5} come in two variants, one for $|q|<1$ and one for $|q|>1$.
We note that these functions are quite different, though there is
a connection. Both are needed in the analytic theory, for instance,
to obtain the theory for all $q$, including the unit circle.

The nonsymmetric Whittaker function defined in \cite{C102}
corresponds to $|q|<1$. The case $|q|>1$ appeared more involved;
in fact, the limiting procedure must be twisted. In this paper,
we define the appropriate limiting procedure and use it to produce
the nonsymmetric Whittaker function for $|q|>1$.

Interestingly, in contrast to the $q,t$\~case, 
both nonsymmetric $q$\~Whittaker functions can be extended
to all $|q|\neq 1$ and {\em almost all} $q$ at the unit circle.
However the relation between these two functions becomes more
sophisticated than in the $q,t$\~case, since they correspond
to different limits, namely $t\to 0$ and $t\to\infty$ (and they solve
different equations).

In the symmetric setting, the Whittaker functions are connected
by a simple conjugation sending $q$ to $q^{-1}$.
Theorem \ref{CONJUGG} below establishes
the analogue of this connection in the nonsymmetric case.
The two functions extend
each other through a variant of the Hardy-Littlewood
continuation theory from $|q|\neq 1$ to $|q|=1$
(almost everywhere). Here roots of unity must be avoided,
but not only such points. The most direct reference is 
\cite{Lub}. It is likely that the roots of unity will be the 
only obstacle for such continuation in the final theory, but
a different approach will be necessary.
\smallskip

{\bf Some other results.}
The paper contains a systematic algebraic theory of the
nil-DAHA and its subalgebras. It includes several
PBW-theorems (some are of unusual type), the theory of
inner products and study of the polynomial and various
induced representations. An important fact is that the
core subalgebra is bi-graded and also has a natural filtration
with the corresponding graded algebra that is a nil-Weyl
algebra. The bi-grading and the filtration are the key
simplifications over the $q,t$\~case. Let us also mention
Theorem \ref{SHAPOL}, which establishes the coincidence
of the bilinear quadratic form in the polynomial representation
defined for the product of the Gaussian and mu-function with
the simplest Shapovalov form. This relation
is very deep in the $q,t$\~theory. It generalizes the technique
of picking residues due to Arthur-Heckman-Opdam;
see \cite{O5},\cite{ChL}.

\comment{
Part II of this paper (in process) will continue the theory of
nonsymmetric global functions, including the analytic aspects
and the Harish-Chandra theory.
It will also contain a complete theory at roots of unity
in the case of $A_1$.
}

\comment{
Due to novelty of technique of spinors in the nonsymmetric theory,
we perform the key verifications in detail
and provide quite a few direct calculations.
}

\comment{
Paper \cite{C102} contains quite a few themes and
references that are not discussed and provided here, including
the origins of the technique of spinors,
related DAHA topics, the $K$\~theoretic connections,
the Demazure level one Kac-Moody characters,
the $p$\~adic and differential theories.
}

\medskip
{\bf Acknowledgements.}
The first author is thankful to D.~Kazhdan for useful
discussions, to E.~Opdam, J.~Stokman and N.~Wallach for 
various talks on the Harish-Chandra and Whittaker theory. 
We also thank S.~Kumar and R.~Rimanyi for 
helpful conversations.
The first author thanks IHES and Caltech for the hospitality.
\medskip

\setcounter{equation}{0}
\section{\sc Polynomial representation}
\subsection{\bf Main definitions}
We consider only the case of $A_1$ in the paper.
Let $\al=\al_1$, $s=s_1$, and $\om=\om_1$, the fundamental weight;
then $\alpha=\al_1=2\omega$ and
$\rho=\om$. The extended affine Weyl group
$\widehat{W}=<s,\om>$
is generated by $s$ and the involution  $\pi=\om s$.

The double affine Hecke algebra
$\HH$ is generated by $Y=Y_{\om_1}=\pi T, T=T_1, X=X_{\om_1}$
subject to the quadratic relation $(T-t^{1/2})(T+t^{-1/2})=0$
and the cross-relations:
\begin{align}\label{dahaone}
&TXT=X^{-1},\ T^{-1}YT^{-1}=Y^{-1},\ Y^{-1}X^{-1}YXT^2q^{1/2}=1.
\end{align}
Using $\pi=YT^{-1}$, the second relation becomes $\pi^2=1$.
The field of definition will be $\C_{q,t}\equal\C(q^{1/4},t^{1/2})$
although $\Z[q^{\pm 1/4},t^{\pm 1/2}]$ is sufficient
for many constructions;
$q^{\pm 1/4}$ will be mainly needed in the automorphisms
$\tau_{\pm}$ below. We will frequently treat $q,t$ as numbers;
then the field of definition will be $\C$.

It is important that $\HH$ at $t=1$ becomes the
Weyl algebra defined as the span
$\lan X,Y\ran /(Y^{-1}X^{-1}YXq^{1/2}=1)$ extended
by the inversion $s=T(t=1)$ sending
$X\mapsto X^{-1}$ and $Y\mapsto Y^{-1}$.

The affine Hecke subalgebra in terms of $Y$
can be written as $\mathcal{H}_Y=\lan Y, T\ran$.

The {\em polynomial representation} is defined as
$\mathscr{X}=\C_{q,t}[X^{\pm 1}]$ with $X$ acting by
the multiplication. The formulas for the
other generators are
\begin{eqnarray*}
T=t^{1/2}s+\frac{t^{1/2}-t^{-1/2}}{X^{2}-1}
\circ (s-1),\ \, Y=\pi T
\end{eqnarray*}
in terms of the (multiplicative) reflection $s(X^n)=X^{-n}$
and $\pi(X^n)=q^{n/2}X^{-n}$ for $n\in\Z$.

We will sometimes set $X=q^x$. Then
\begin{align}\label{formalgauss}
&s(x)=-x,\, \om(f(x))=f(x-1/2),\
\pi=\om s,\, \pi(x)=1/2-x.
\end{align}

\subsubsection{\sf Automorphisms}
The following map can be extended to an {\em anti-involution} on $\HH$:
\begin{align}\label{vphdef}
\varphi: X \leftrightarrow Y^{-1},\ T\mapsto T,\ q,t\mapsto q,t.
\end{align}
The first two
relations in (\ref{dahaone}) are obviously fixed by $\varphi$;
as for the third, check that
$\varphi(Y^{-1}X^{-1}YX)=Y^{-1}X^{-1}YX$. Switching $X$ and $Y$
can be also achieved using the {\em involution}
\begin{align}\label{vepanti}
\vep: X \leftrightarrow Y,\ T\mapsto T^{-1},
\ q^{1/4}\mapsto q^{-1/4},\
t^{1/2}\mapsto t^{-1/2}.
\end{align}

The conjugation by the Gaussian $q^{x^2}$
preserves $\HH$. The Gaussian obviously
belongs to a completion of $\mathscr{X}$.
It satisfies:
$$
\om(q^{x^2})= q^{1/4}X^{-1}q^{x^2},\ \,
\om(q^{-x^2})=q^{-1/4}X q^{-x^2}.
$$
The conjugation $A\mapsto \tau_+(A)=q^{x^2}\,A\,q^{-x^2}$
for $A\in \HH$ satisfies:
\begin{align}\label{tau+def}
\tau_+(X)=X,\ \tau_+(T)=T,\ \tau_+(Y)=q^{-1/4}XY,\
\tau_+(\pi)=q^{-1/4}X\pi.
\end{align}
To see this, use that
$$
Y=\om\circ(t^{1/2}+\frac{t^{1/2}-t^{-1/2}}{X^{-2}-1}
\circ (1-s)).
$$

Applying $\varphi$ we obtain an automorphism
\begin{align}\label{tau-def}
\tau_-=\varphi\tau_+\varphi,\
\tau_-(Y)=Y,\, \tau_-(T)=T,\, \tau_-(X)=q^{1/4}YX.
\end{align}

The generalized Fourier transform, corresponds to the
following automorphism of $\HH$ (it is not an
involution)\,:
\begin{align}
\si(X)= Y^{-1},\ \si(T)=T,\ &\si(Y)=q^{-1/2}Y^{-1}XY=XT^2,\
\si(\pi)=XT,
\notag\\
\si\ &=\  \tau_+\tau_-^{-1}\tau_+\ =\ \tau_-^{-1}\tau_+\tau_-^{-1}.
\label{tautautau}
\end{align}
Check that $\,\si\tau_+=\tau_-^{-1}\si,\ \,
\si\tau_+^{-1}=\tau_-\si.$ Also,
$$
\vph\tau_+=\tau_-\vph,\ \vph\si=\si^{-1}\vph,\
\vep\tau_+=\tau_-\vep,\ \vep\si=\si^{-1}\vep.
$$

Due to the group nature of the definition of $\HH$,
we have the inversion anti-involution $\HH\ni
H\mapsto H^\ast$\,:
$$
X^\ast=X^{-1},\, Y^\ast=Y^{-1},\, T^\ast=T^{-1},\,
(q^{1/4})^\ast=
q^{-1/4},\, (t^{1/2})^\ast=t^{-1/2}.
$$
It commutes with all automorphisms and anti-automorphisms
of $\HH$. Note that $\ast=\vph\vep=\vep\vph$.

The following {\em anti-involutions} preserving $X,T,q,t$
will be important below, $\diamond\equal\vph\si$ and
$\psi\equal\diamond\tau_+^{-1}$\,:
\begin{align}
&\pi^\diamond =\pi,\ Y^\diamond =TYT^{-1}=\pi Y\pi,\label{diamondef}
\ \, \psi(\pi)=\tilde{\pi}\equal\tau_+(\pi)\\
=&q^{-1/4}X\pi =q^{1/4}\pi X^{-1},\,
\psi(Y^{-1})= \tilde{X}\equal\tilde{\pi}T^{-1}=
q^{1/4}YX.\label{psidef}
\end{align}

\subsubsection{\sf Inner products}
\label{sec:innprod}
The polynomial representation can be supplied with
inner products in various ways. The main ones are in
terms of the function:
\begin{equation}\label{mutildemuone}
\mu(X;q,t)\equal\prod_{j=0}^\infty
\frac{(1-q^jX^2)(1-q^{j+1}X^{-2})}
{(1-tq^jX^2)(1-tq^{j+1}X^{-2})}.
\end{equation}
For the constant term functional
$$
\mathscr{X}\ni f\mapsto  \lan f \ran =\hbox{ct}(f)\in \C_{q,t},
$$
we define three symmetric inner products in $\mathscr{X}$:
\begin{align}\label{innerp}
&(\!(f,g)\!)=\lan fg^*\mu\ran,\
\lan f,g\ran=\lan fg\mu \ran,\
\lan f,g\ran'=\lan fg \tga' \mu \ran.
\end{align}
Here $\tga'\equal\sum_{m=-\infty}^\infty q^{mx+m^2/4}$
is an expansion of $q^{-x^2}$ in the following
sense. The product $q^{x^2}\tga'$ is a $\Z/2$\~periodic
function in terms of $x$ provided that $|q|<1$.
Recall that $X=q^x$.

The anti-involutions of $\HH\ni H$ corresponding to
these forms are
those from (\ref{diamondef},\ref{psidef}), namely:
\begin{align}\label{innerinv}
&(\!(f,H(g))\!)=(\!(H^*(f),g)\!),\
\lan f,H(g)\ran=\lan H^{\diamond}(f),g\ran,\\
&\lan f,H(g)\ran'=\lan H^\psi(f), g\ran' \for
f,g\in \mathscr{X},\, H^\psi=\psi(H).\notag
\end{align}
\smallskip

\subsection{\bf The E-polynomials}{\label{sect:macpoly}}
Let us assume that
$k$ is generic; we set $t=q^k$.
The definition of nonsymmetric polynomials is as follows:
\begin{align}\label{nonsymp}
YE_{n}=q^{-n_{\#}}E_{n}\for n\in \Z,&&&&\\
n_{\#}=\left\{\begin{array}{ccc}\frac{n+k}{2}
&  & n>0, \\\frac{n-k}{2} &  & n\le 0,\end{array}\right\},
\text{\, note that }\,  0_{\#}=-\frac{k}{2}.
\end{align}
The normalization is
$E_{n}=X^{n}+\text{ ``lower terms'' },$
where by ``lower terms'', we mean
polynomials in terms of $X^{\pm m}$ as $|m|<n$
and, additionally,  $X^{|n|}$ for negative $n$.
It gives a filtration in $\mathscr{X}$;
check that $Y$ preserves the filtration,
which justifies the definition from (\ref{nonsymp}).

The $E_{n} \,(n\in \Z)$ are called the
{\em nonsymmetric Macdonald polynomials} or simply
$E$\~polynomials.
Obviously,  $E_{0}=1,\, E_{1}=X$.

\subsubsection{\sf The intertwiners}
The first intertwiner comes from the AHA theory:
$$
\Phi\equal T+\frac{t^{1/2}-t^{-1/2}}{Y^{-2}-1}\,:\
\Phi Y=Y^{-1}\Phi.
$$
The second is $\Pi\equal q^{1/4}\tau_+(\pi)$;
obviously, $\Pi^2=q^{1/2}$.
Explicitly,
\begin{align}\label{Piformula}
\Pi=X\pi=q^{1/2}\pi X^{-1}\,:\
\Pi Y=q^{-1/2}Y^{-1}\Pi.
\end{align}
Use that $\phi(\Pi)=\Pi$ to deduce the latter
relation from $\Pi X\Pi^{-1}=q^{1/2}X^{-1}$.
The $\Pi$\~type intertwiner is due to Knop and Sahi for $A_n$
(the case of arbitrary reduced systems was considered in \cite{C1}).
Since $\Phi,\Pi$ ``intertwine" $Y^{\pm 1}$, they can be
used for generating the $E$\~polynomials. Namely,
\begin{align}{\label{signE}}
&E_{n+1}=q^{n/2}\Pi (E_{-n}) \for n\ge 0,\\
&E_{-n}=t^{1/2}(T+\frac{t^{1/2}-t^{-1/2}}{q^{2n_{\#}}-1})E_{n}
\label{interphi}
\end{align}
and, beginning with $E_0=1$, one can readily construct the
whole family of $E$\~polynomials.
For instance,
\begin{align*}
T(X)&=t^{1/2}X^{-1}+\frac{(t^{1/2}-t^{-1/2})(X^{-1}-X)}{X^{2}-1}\\
&=t^{1/2}X^{-1}-(t^{1/2}-t^{-1/2})X^{-1} = t^{-1/2}X^{-1}, \\
E_{-1}&=t^{1/2}(T+\frac{t^{1/2}-t^{-1/2}}{qt-1})E_{1}
=X^{-1}+\frac{1-t}{1-t q}X.
\end{align*}
Using $\Pi$,
\begin{align*}
&E_{2} = q^{1/2}\Pi E_{-1}=  X^{2}+q\frac{1-t}{1-tq}.
\end{align*}
Applying $\Phi$ and then $\Pi$ again,
\begin{align*}
&E_{-2} = X^{-2}+\frac{1-t}{1-tq^{2}}X^{2}
+\frac{(1-t)(1-q^{2})}{(1-tq^{2})(1-q)},\\
&E_{3}=X^{3}+q^2\frac{1-t}{1-tq^{2}}X^{-1}+q\frac{(1-t)(1-q^{2})}
{(1-tq)(1-q)}X.
\end{align*}
It is not difficult to find the general formula. See, e.g.,
(6.2.7) from \cite{Ma4} for integral $k$. However, recalculating
these formulas from integral $k$ to generic $k$ is not
too simple; we will provide the exact formulas for
the $E$\~polynomials
below (in the form we need them).

The following properties of $n_{\#}$ reflect (\ref{signE}):
\begin{align}\label{sharpsym}
&(1-n)_{\#}=1/2-n_{\#} \hbox{ for all } n\in \Z,\
(-n)_{\#}=-n_\# \hbox{ when } n\neq 0.
\end{align}

\subsubsection{\sf The E-Pieri rules}
For any $n\in \Z$, we have the {\em evaluation formula}
\begin{eqnarray}\label{evale}
E_{n}(t^{-1/2})=t^{-|n|/2}\prod_{0<j<|\tilde{n}|}
\frac{1-q^{j}t^{2}}{1-q^{j}t},
\end{eqnarray}
where
\begin{align}\label{tilden}
|\tilde{n}|=|n|+1 \hbox{\, if\, } n\leq 0 \and
|\tilde{n}|=|n| \hbox{\, if\, } n>0.
\end{align}

It is used to introduce the {\em nonsymmetric spherical polynomials}
$$\mathcal{E}_{n}=\frac{E_{n}}{E_{n}(t^{-1/2})}.$$
This normalization is important in many constructions
due to the
{\em duality formula}: $\e_m(q^{n_\#})=\e_n(q^{m_\#})$.
The Pieri rules are the simplest for the
$E$\~spherical polynomials:
\begin{eqnarray}\label{pierie}
&X\mathcal{E}_{n}
=\frac{t^{-1/2\pm1}q^{-n}-t^{1/2}}
{t^{\pm1}q^{-n}-1}\mathcal{E}_{n+1}
+\frac{t^{1/2}-t^{-1/2}}{t^{\pm1}q^{-n}-1}\mathcal{E}_{1-n},\\
&X^{-1}\mathcal{E}_{n}
=\frac{t^{1/2\pm1}q^{-n+1}-t^{-1/2}}
{t^{\pm 1}q^{-n+1}-1}\mathcal{E}_{n-1}
-\frac{t^{1/2}-t^{-1/2}}{t^{\pm1}q^{-n+1}-1}\mathcal{E}_{1-n}.
\label{pierie1}
\end{eqnarray}
Here the sign is $\pm=+$\, if\, $n\leq 0$\  and $\pm=-$ if $n>0$.
These formulas give an alternative approach to
constructing the $E$\~polynomials and establishing
their connections with other theories, for instance, with
the $p$\~adic one.
\smallskip

Let us provide the norm formulas for the spherical and standard
polynomials
\begin{align}
\lr \e_m \e_n^*\mu_\circ\rr&=
\de_{mn}\prod_{0<j<|\tilde{n}|}\frac{1-q^j}{t^{-1}-q^jt},
\label{fnormep}\\
\lr E_m E_n^*\mu_\circ\rr&=\de_{mn}
\prod_{0<j<|\tilde{n}|}\frac{(1-q^j)(1-q^j t^2)}
{(1-q^j t)(1-q^{j}t)},\label{fnorme}
\end{align}
where we use $|\tilde{n}|$ from (\ref{tilden}),
and $\mu_\circ \equal \mu/\lr\mu\rr$.
\medskip

\subsubsection{\sf Rogers' polynomials}{\label{sect:Rogers}}
Let us introduce the {\em Rogers polynomials}:
\begin{align}\label{pviae}
P_{n}&=(1+t^{1/2}T)\bigl(E_{n}\bigr)
=(1+s)\bigl(\frac{t-X^2}{1-X^2}E_n\bigr)\\
&=E_{-n}+\frac{t-tq^n}{1-tq^n}E_n \for n\ge 0.\notag
\end{align}
The leading term is $X^n+X^{-n}$:\ $P_n=X^n+X^{-n}+ $``lower terms".
They are eigenfunctions
of the following well-known operator
\begin{align}\label{Lopertaor}
\l=\frac{t^{1/2}X-t^{-1/2}X^{-1}}{X-X^{-1}}\Ga+
\frac{t^{1/2}X^{-1}-t^{-1/2}X}{X^{-1}-X}\Ga^{-1},
\end{align}
where we set $\Ga(f(x))=f(x+1/2),\, \Ga(X)=q^{1/2}X$,
The operator $\l$ is the restriction of $Y+Y^{-1}$ to symmetric
polynomials; this is the key point of the DAHA approach
to the theory of the Macdonald polynomials.

The exact eigenvalues are as follows:
\begin{equation}{\label{eqn:LP}}
\l(P_n)\ = \ (q^{n/2}t^{1/2}+q^{-n/2}t^{-1/2})\,P_n,\ n\ge 0.
\end{equation}
It is obvious from the latter that the $P$\~polynomials are
$\ast$\~invariant. Using directly (\ref{conjep}) and that
$\eta(T)=T^{-1}$, for $\eta$ defined by (\ref{etaqt}),
\begin{align*}
P_{n}^\ast=((1+t^{1/2}T)\bigl(E_{n}\bigr))^\ast=
(1+t^{-1/2}T^{-1})\bigl(E_{n}^\ast\bigr)=\\
(1+t^{-1/2}T^{-1})\bigl(t^{1/2}T(E_{n})\bigr)=
(1+t^{1/2}T)\bigl(E_{n}\bigr)=P_n.
\end{align*}

The evaluation formula reads:
$$
P_n(t^{\pm1/2})=t^{-n/2}\prod_{0\le j\le n-1}
\frac{1-q^{j}t^{2}}{1-q^{j}t}.
$$
The spherical $P$\~polynomials $\p_n\equal P_n/P_n(t^{1/2})$
satisfy the duality $\p_n(t^{1/2}q^{m/2})=\p_m(t^{1/2}q^{n/2})$.
The norm formula reads:
\begin{align}\label{normppolsbara1t}
&\lan P_{m}(X)P_{n}(X)
\mu_\circ\ran\
=\ \de_{mn}\prod_{j=0}^{n-1}\frac{(1-q^{j+1})
(1-t^2q^{j})}{(1-tq^{j+1})
(1-tq^{j})}\,,
\end{align}
as $m,n\geq 0.$

\subsubsection{\sf Explicit formulas}{\label{sect:formulas}}
Let us begin with the well-known formulas for the
Rogers polynomials ($n\ge 0$):
\begin{align}\label{exactp}
&P_{n}=X^{n}+X^{-n}+\sum_{j=1}^{[n/2]}M_{n-2j}\prod_{i=0}^{j-1}
\frac{(1-q^{n-i})}{(1-q^{1+i})}\,\frac{(1-tq^{i})\  }
{(1-tq^{n-i-1})},
\end{align}
where $M_{n}=X^{n}+X^{-n}\, (n>0)$ and $M_{0}=1$.

The formulas for the $E$\~polynomials are as follows ($n>0$):
\begin{align}\label{exacte-}
E_{-n}&=X^{-n}+X^n\frac{1-t}{1-tq^n}+
\sum_{j=1}^{[n/2]}X^{2j-n}\,\prod_{i=0}^{j-1}
\frac{(1-q^{n-i})}{(1-q^{1+i})}\,\frac{(1-tq^{i})}{(1-tq^{n-i})}
\notag\\
&+\sum_{j=1}^{[(n-1)/2]}X^{n-2j}\,
\frac{(1-tq^{j})}{(1-tq^{n-j})}\prod_{i=0}^{j-1}
\frac{(1-q^{n-i})}{(1-q^{1+i})}\,\frac{(1-tq^{i})}{(1-tq^{n-i})},
\end{align}
\begin{align}\label{exacte+}
E_{n}=X^{n}+&
\sum_{j=1}^{[n/2]}X^{2j-n}\,q^{n-j}\,
\frac{(1-q^{j})}{(1-q^{n-j})}\prod_{i=0}^{j-1}
\frac{(1-q^{n-i-1})}{(1-q^{1+i})}\,
\frac{(1-tq^{i}\ ) }{(1-tq^{n-i-1})}
\notag\\
&+\sum_{j=1}^{[(n-1)/2]}X^{n-2j}\,q^j\,
\prod_{i=0}^{j-1}\frac{(1-q^{n-i-1})}{(1-q^{1+i})}\,
\frac{(1-tq^{i}\ ) }{(1-tq^{n-i-1})}.
\end{align}
\medskip

\setcounter{equation}{0}
\section{\sc Global functions}
\subsection{\bf Spherical functions}
We will use $n_\#=(n+\sgn(n)k)/2$ for integers $n\neq 0$
and $0_\#=-k/2$. From now on $\sgn(0)\equal -1$, i.e., we will
always treat $0$ as a negative number; cf. (\ref{nonsymp}).
As above, let $\tga'\equal
\sum_{n=-\infty}^\infty q^{n^2/4}X^n, X=q^x$; recall that
the product $q^{x^2}\tga'(x)$ is a
$\Z/2$\~periodic function of $x$. We denote the
constant term functional, the coefficient of $X^0$, by
$\lan\,\cdot\,\ran$.
We set
\begin{align}\label{mucontone}
&\mu_\circ\equal\mu/\lr\mu\rr= 1+\frac{t-1}{1-qt}(X^2+qX^{-2})
+\ldots\ ,\\
&\where \lr \mu \rr\ =\ \prod_{j=1}^\infty \frac{(1-tq^j)^2}
{(1-t^2q^j)(1-q^j)}.\notag
\end{align}

Recall that
$$
X^\ast=X^{-1},\ (q^{\,1/4})^\ast=q^{-1/4},\
(t^{\,1/2})^\ast=t^{-1/2}.
$$
The series $\mu_\circ$ is $\ast$\~invariant.

\subsubsection{\sf Gauss-type integrals}
For arbitrary $m,n\in \BZ,$
\begin{align}\label{f81}
&\lr \e_n\e_m\tga'\mu_\circ\rr
=q^{\frac{m^2+n^2+2k(|m|+|n|)}{4}}
\e_m(q^{n_\#})\lr\tga'\mu_\circ\rr,\\
\label{f82}
&\lr \e_n\e_m^\ast\tga'
\mu_\circ\rr =q^{\frac{m^2+n^2+2k(|m|+|n|)}{4}}
\e_m^\ast(q^{n_\#})\lr \tga'\mu_\circ\rr.
\end{align}
In these formulas,
\begin{align}\label{mehmacone}
&\lr \tga'\mu_\circ\rr=
\prod_{j=1}^\infty \frac{1-q^j}{1-tq^{j}},\
\lr \tga'\mu\rr=\prod_{j=1}^\infty \frac{1-tq^j}{1-t^2q^{j}}.
\end{align}
See \cite{C101}, Theorem 2.7.1.

The $\ast$\~conjugations of the $\e$\~polynomials are as follows:
\begin{align}\label{conjepol}
&\e_m^\ast=t^{-\frac{1}{2}}T(\e_m)=
q^{-(m+(1+\hbox{\tiny sgn}(m))k)/2}\,\pi(\e_m)=
t^{-\frac{1}{2}}X^{-1}\e_{1-m},\\
& \where \bigl(m+(1+\hbox{sgn}(m))k\bigr)/2=m_{\#}+k/2.\notag
\end{align}
where $m\in\Z,\,\sgn(0)=-1$.
See, e.g., \cite{C101}, Proposition 2.5.13.

The conjugation is somewhat simpler in terms of the
$E$\~polynomials:
\begin{align}\label{conjep}
&E_m^\ast=
t^{\,\hbox{\tiny sgn}(m)/2}T(E_m)=q^{-\frac{m}{2}}\pi(E_m)=
X^{-1}E_{1-m} \for m\in \Z.
\end{align}

Using (\ref{evale}), we obtain:
\begin{align}\label{f81ee}
&\lr E_n E_m\tga'\mu_\circ\rr
=\prod_{j=1}^{|\tilde{n}|-1}
\frac{1-q^{j}t^2}{
1-q^{j}t}\
\prod_{j=1}^{\infty}
\frac{(1-q^j)}{
(1-q^{j}t)}\
q^{\frac{m^2+n^2+2k|m|}{4}}E_m(q^{n_\#}),\\
\label{f82ee}
&\lr E_n E_m^\ast\tga'\mu_\circ\rr
=
\prod_{j=1}^{|\tilde{n}|-1}
\frac{1-q^{j}t^2}{
1-q^{j}t}\
\prod_{j=1}^{\infty}
\frac{(1-q^j)}{
(1-q^{j}t)}\
q^{\frac{m^2+n^2+2k|m|}{4}}
E_m^\ast(q^{n_\#}).
\end{align}
Switching to $\mu$,
the  formulas from (\ref{f81},\ref{f82}) read:
\begin{align}\label{f81e}
&\lr E_n E_m\tga'\mu\rr
=\prod_{j=|\tilde{n}|}^\infty
\frac{1-q^{j}t}{1-q^{j}t^2}\
q^{\frac{m^2+n^2+2k|m|}{4}}E_m(q^{n_\#}),\\
\label{f82e}
&\lr E_n E_m^\ast\tga'\mu\rr
=\prod_{j=|\tilde{n}|}^\infty
\frac{1-q^{j}t}{1-q^{j}t^2}\
q^{\frac{m^2+n^2+2k|m|}{4}}
E_m^\ast(q^{n_\#}).
\end{align}
Use the second formula from (\ref{mehmacone}).

\subsubsection{\sf Fourier transforms}
The formulas for the Gauss integrals can be represented
algebraically, up to a (global)
coefficient of proportionality, as follows. In the theorem below,
$q^{\pm x^2}$ will be the {\em formal}
Gaussian satisfying the defining relations from
(\ref{formalgauss}):
\begin{align*}
&s(q^{\pm x^2})\,=\,q^{\pm x^2},\, \om(q^{\pm x^2})\,=\,
q^{\pm 1/4}X^{\mp 1}q^{\pm x^2}, \\
&\where s(x)=-x,\ \om(f(x))=f(x-1/2),\\
&Y(q^{\pm x^2})\,=\,\om(q^{\pm x^2})\,=
\,q^{\pm 1/4}X^{\mp 1} q^{\pm x^2},\\
&\where \pi=\om s,\ \pi(f(x))=f(1/2-x).
\end{align*}
For instance,
one can take $\tga'$ here as $q^{-x^2}$.


\begin{theorem}\label{FOURE}
(i) The $\C(q^{1/4},t^{1/2})$\~linear map defined by
$$
\F_\si: \x q^{-x^2}\,\ni\, E_m q^{-x^2}\ \mapsto \
q^{\frac{m^2+2k|m|}{4}}E_m q^{+x^2}\,\in\, \x q^{+x^2},
$$
where $m\in \Z$, induces the automorphism $\si$ on the
algebra $\HH$ naturally acting in $\x q^{\mp x^2}$.
Equivalently, $\si=\vph\,\diamond$.

(ii) Let
$$
\F_\vep: \x q^{+x^2}\,\ni\, E_m q^{+x^2}\ \mapsto\
q^{\frac{m^2+2k|m|}{4}}E_m^* q^{+x^2}\,\in\, \x q^{+x^2},
$$
where $q^{1/4}\mapsto q^{-1/4}$ and $t^{1/2}\mapsto t^{-1/2}$,
i.e., the constants are $\ast$\~conjugated under the action of
$\F_\vep$. This map
induces the involution $\vep$ on $\HH$;
equivalently, $\vep=\vph\ast$.

(iii) Formally conjugating (\ref{f82ee}), let
\begin{align*}
\F'_\vep:\, \x q^{-x^2}\,\ni\ &E_m q^{-x^2}\,\mapsto\,
t^{(\hbox{\tiny\rm sgn}(m)+1)/2}
q^{-\frac{m^2+2k|m|}{4}}E_m q^{-x^2}\,\in\, \x q^{-x^2},
\end{align*}
where $q^{1/4}\mapsto q^{-1/4}$ and $t^{1/4}\mapsto t^{-1/4}$.
This map induces the same involution $\vep$ on $\HH$.\sq
\end{theorem}

Here we closely follow \cite{C5}, Theorem 5.1. We will
clarify the way to calculate the corresponding DAHA
isomorphisms in the next section.

Note that the (formal) conjugation of ($i$) eventually
leads to the transform equivalent to that from ($iii$).
Use (\ref{conjepol}) and the first relation from
(\ref{sharpsym}).
\medskip

\subsubsection{\sf Reproducing kernels}
The general fact is that
if the map $f\mapsto \widehat{f}(m)=\lr \e_m, f\rr_\al$
can be naturally
extended to a morphism $\F$ of $\HH$\~modules, then it
corresponds to the automorphism $\be=\vph\al$, where  $\vph$
is from (\ref{vphdef}), $\al$
is the anti-involution of $\HH$ associated to
$\lr\,\cdot\,,\,\cdot\,\rr_\al$. The latter is assumed to be
nondegenerate symmetric inner product or an anti-symmetric one
(for the second component). Here $f$ may belong to various
function spaces, including $\mathscr{X}q^{lx^2} (l\in \Z)$, their
completions and their delta-counterparts.

The fundamental concept of Fourier analysis (and representation
theory) is the {\em reproducing kernel}. Given a transform $\B$,
it is defined as follows:
\begin{align}\label{reproker}
G(X,\La)&=\sum_n  f'_n(X)\,\B(f_n)(\La),
\where \lr f'_m, f_n\rr_\al=\de_{mn}.
\end{align}
The basis  $\{f_n\}$ in the
initial function space can be arbitrary; $\{f'_n\}$ is
its dual.
Provided the existence, we assume that $\B$ induces the isomorphism
$\be$ for the operators. Then
\begin{align}\label{reprokerpro}
\lr G(X,\La), f(X)\rr_\al=\B(f)(\La),\
H(G)=\widehat{\be\al(H)}(G)
\end{align}
for $H\in \HH$, where we consider $f_n$ and $f'_n$ as functions
of $X$ and suppose that $\B$ sends them to functions of
$\La$, i.e., $\B:f(X)\mapsto \widehat{f}(\La)$.
By $H, \widehat{H}$, we mean this operator acting in the domain
and range of $\B$, i.e., on functions of $X$ and $\La$ respectively.
For instance, $Y(G(X,\La))=\La^{-1}G(X,\La)$.

The formal independence of $G$ of the choice
of the basis $\{f_n\}$
becomes a non-trivial issue in the functional analysis.
The better the basis (and the space),
the better analytic theory of $G$.
In the $q$\~theory, the best choice is the
basis of the $E$\~polynomials multiplied by the Gaussian.

\subsubsection{\sf Nonsymmetric global functions}
Theorem \ref{FOURE} results in the following two
formulas (both are from \cite{C5}). Using ($i$) with
the pairing $\lr fg\mu_\circ\rr$ for
$$
f_n=E_n q^{-x^2}
\and f_n'=E_n^\ast q^{x^2}/\lr E_nE_n^\ast\mu_\circ\rr
$$
and ($ii$) with pairing  $\lr fg^*\mu_\circ\rr$ for
$$
f_n=E_n q^{x^2} \and f_n'=E_n q^{x^2}/\lr E_nE_n^\ast\mu_\circ\rr,
$$
we arrive at the following coinciding formulas:
\begin{align}\label{gxla}
\frac{\tga'(X)\tga'(\La)}{\tga'(t^{1/2})}
G(X;\La)&=\sum_{n=-\infty}^\infty\,q^{\frac{|n|^2}{4}}\,
t^{\frac{|n|}{2}}
\,\frac{E_n^*(X)E_n(\La)}{\lr E_n E_n^*\mu_\circ\rr} \\
&=\sum_{n=-\infty}^\infty\,q^{\frac{|n|^2}{4}}\,t^{\frac{|n|}{2}}
\,\frac{E_n(X)E_n^*(\La)}{\lr E_n E_n^*\mu_\circ\rr}\,, \notag
\end{align}
where we use $\tga'$ for $q^{-x^2}$.
The inner products
$\lr E_n E_n^\ast\mu_\circ\rr$ in these formulas and in those below
are provided in (\ref{fnorme}).

It is assumed here that $|q|<1$; $X,\La$ are
arbitrary apart from the zeros of $\tga'(X)\tga'(\La)$.
The formula converges for any $q$ not equal to a root of
unity and is meromorphic. When $|q|=1$, the analyticity
is understood with respect to the directions not tangent
to the unit circle ($q$ must not be a root of unity).

Checking the coincidence of these expressions for $G$
is not difficult because they satisfy the same relations
from (\ref{reprokerpro}). The basic ones are:
\begin{align}\label{relyxt}
&Y(G)=\La^{-1}\,G,\
X^{-1}\,G=Y_\La(G),\
T(G)=T_\La(G),
\end{align}
where the operators indexed by $\La$ act in terms of $\La$.
One can check this coincidence directly.
Substitute $E_m^\ast=t^{\,\hbox{\tiny sgn}(m)/2}T(E_m)$
and use that $T(G)=T_\La(G)$; equivalently, the relation
$E_m^*=X^{-1}E_{1-m}$ can be used.
\smallskip

Similarly, the formulas from ($iii$) coupled with
the inner product $\lr fg^*\mu_\circ\rr$ and the
basic elements
$$
f_n=E_n q^{-x^2}, \
n\in \Z, \ f_n'=E_n q^{-x^2}/\lr E_nE_n^\ast\mu_\circ\rr
$$
result in the second function:
\begin{align}\label{gxlaprime}
&\frac{\tga(X)\tga(\La)}{\tga(t^{1/2})}
G^\checkmark(X;\La)\\
&=\sum_{n=-\infty}^\infty\,
q^{-\frac{|n|^2}{4}}t^{-\frac{|n|}{2}}\,
t^{(\hbox{\tiny\rm sgn}(n)+1)/2}\,
\frac{E_n(X)E_n(\La)}{\lr E_n E_n^*\mu_\circ\rr}\,,\notag
\end{align}
where $|q|>1$ and
$\tga\equal\sum_{n=-\infty}^\infty q^{-n^2/4}X^n$.

It satisfies the relations from (\ref{reprokerpro})
and (\ref{relyxt}):
$$
Y(G^{\checkmark})=\La^{-1}\,G^{\checkmark},\,
X^{-1}\,G^{\checkmark}=\widehat{Y}(G^{\checkmark}),\,
T(G^{\checkmark})=\widehat{T}(G^{\checkmark}),
$$
i.e., the same relations as for $G$.

The functions $G,G^{\checkmark}$ are called the {\em global
nonsymmetric spherical functions}.

The symmetric (even) global functions read as follows:
\begin{align}\label{fxla}
\frac{\tga'(X)\tga'(\La)}{\tga'(t^{1/2})}
F(X;\La)&=\sum_{n=0}^\infty\,q^{\frac{n^2}{4}}\,
t^{\frac{n}{2}}
\,\frac{P_n(X)P_n(\La)}{\lr P_n P_n\mu_\circ\rr}\ , \for |q|<1,\\
\frac{\tga(X)\tga(\La)}{\tga(t^{1/2})}
F^\checkmark(X;\La)&=\sum_{n=0}^\infty\,q^{-\frac{n^2}{4}}\,
t^{-\frac{n}{2}}
\,\frac{P_n(X)P_n(\La)}{\lr P_n P_n\mu_\circ\rr}\ ,\, \,\ |q|>1.
\end{align}
Both satisfy the same
difference equation $\l (F)=(\La+\La^{-1})F$
for $\l$ from (\ref{Lopertaor}).
The second function is obtained from the first upon the formal
star-conjugation
$$
q\mapsto q^{-1},\, t\mapsto t^{-1},\, X\mapsto X^{-1},\,
\La\mapsto \La^{-1},
$$
preserving $\l$ ($\Ga$ remains fixed).

The action of $\ast$\~conjugation on the $G$\~functions is
more involved.

\begin{theorem}\label{CONJUGG}
Let $\tilde{G}=(G^\checkmark)^\ast$ be the result of the formal
conjugation of $q,t$ and $X,\La$. Namely, we
set $t=q^k, X=q^x, \La=q^\la$ and replace $q$ by $q^{-1}$
in the formula for $G^\checkmark$
without changing $k,x,\la$;
for instance, $\tga^*=\tga'$. The convergence of the series for
$\tilde{G}$ is for $|q|<1$ and
\begin{align}\label{conjgfnct}
& (t^{1/2}T^{-1})(\tilde{G}(X,\La))\ =\ G(X,\La)\ =\
t^{1/2}\La\pi(\tilde{G}(X,\La)).
\end{align}
\end{theorem}
{\it Proof.} The conjugation of functions
results in the application of the automorphism $\eta$
at level of operators (not $\ast$ as one may expect).
See \cite{C101}:
\begin{align}\label{etaqt}
\eta : \
&T \mapsto T^{-1},\ \pi \mapsto \pi,\ X \mapsto X^{-1},\
q^{1/4} \mapsto q^{-1/4},\, t^{1/2}\mapsto t^{-1/2},\\
&H(f^\ast)\, =\,
(\eta(H)(f))^\ast \for f \in \mathscr{X}=\C_{q,t}[X^{\pm 1}], \,
H\in \HH.\notag
\end{align}

It must be done in the $X$\~space and in the $\La$\~space.
For instance, the eigenvalue problem
$Y(G^\checkmark)=\La^{-1}G^\checkmark$ becomes
$\eta(Y)(\tilde{G})=\La\tilde{G}$. Let us calculate the
conjugations of the relations from (\ref{relyxt}).

Since $\eta(Y)=TY^{-1}T^{-1}$, we obtain:
\begin{align}\label{relyxteta}
&T(\tilde{G})=T_\La(\tilde{G}),\ \,
\eta(Y)(\tilde{G})=\La \tilde{G},\ \,
X\tilde{G}=\eta(Y)_\La(G),
\\
TY^{-1}&T^{-1}(\tilde{G})=\La\tilde{G}=
\La TT_\La^{-1}\tilde{G}\ \Rightarrow\
Y^{-1}(T^{-1}(\tilde{G}))=
\La \,T_\La^{-1}(\tilde{G}),\notag\\
X\tilde{G}&=(TY^{-1}T^{-1})_\La(\tilde{G})
\ \ \Rightarrow\ \
X\tilde{G}=
\bigl(T_\La Y_\La^{-1}T^{-1}_{\La}T_\La T^{-1}\bigr)(\tilde{G})
\notag\\
=T^{-1}\bigl(T_\La &Y_\La^{-1}\bigr)(\tilde{G})\ \Rightarrow\
X^{-1}(T^{-1}(\tilde{G}))=Y_\La (T_\La^{-1}(\tilde{G}))=
Y_\La (T^{-1}(\tilde{G})).
\notag
\end{align}

Therefore $T^{-1}(\tilde{G})$ satisfies {\em all} relations
from (\ref{relyxt}) and must coincide with $G$ up to
normalization. The normalization factor can be
readily determined. As a matter of fact, this general justification
is not necessary at all (for $A_1$) because of the following explicit
identification of the corresponding series:
\begin{align*}
\bigl(&\sum_{n=-\infty}^\infty\,
q^{-\frac{|n|^2}{4}}t^{-\frac{|n|}{2}}\,
t^{+(\hbox{\tiny\rm sgn}(n)+1)/2}\,
\frac{E_n(X)E_n(\La)}{\lr E_n E_n^*\mu_\circ\rr}\bigr)^\ast\\
=&\sum_{n=-\infty}^\infty\,
q^{+\frac{|n|^2}{4}}t^{+\frac{|n|}{2}}\,
t^{-(\hbox{\tiny\rm sgn}(n)+1)/2}\,
\frac{E^*_n(X)E^*_n(\La)}{\lr E_n E_n^*\mu_\circ\rr}\\
=&\sum_{n=-\infty}^\infty\,
q^{+\frac{|n|^2}{4}}t^{+\frac{|n|}{2}}\,
t^{-(\hbox{\tiny\rm sgn}(n)+1)/2}\,
\frac{t^{\hbox{\tiny\rm sgn}(n)/2} T(E_n(X)) E^*_n(\La)}
{\lr E_n E_n^*\mu_\circ\rr}.
\end{align*}

The second equality in (\ref{conjgfnct}) follows from
the first:
$$
t^{1/2}\pi(\tilde{G})\ =\
t^{1/2}(YT^{-1})(\tilde{G})\ =\
Y(t^{1/2}T^{-1}(\tilde{G}))\ =\ \La^{-1}\,G.
$$
\sq
\smallskip

We see that the functions $G$ and $G^\checkmark$ are
conceptually connected through the formal conjugation,
however, a certain twist is needed.
It is a general fact (which was not observed in \cite{C5}).
\smallskip

\comment{
\subsection{\bf P-adic limit}\label{sec:rank1}
The limiting case $q\to 0$ correspond to the
classical $p$\~adic theory. The following is the
Fourier-dual of the previous considerations under
this limit.

Let us introduce the {\em symmetrizer}
$$\mathscr{P}_{+}=\frac{1+t^{1/2}T}{1+t}.$$
For any $m\in \mathbb{Z}$, let $\delta_{m}=\delta_{m\om}$
and $\varepsilon_{m}=\delta_{m\om}^{\sharp}
=t^{-m/2}T_{m\om}\mathscr{P}_{+}$.

Then we have for $m\geq 0$,
\begin{align}
&T\varepsilon_{m}\ =\ t^{1/2}\varepsilon_{-m},\label{Tvepm}\\
&T\varepsilon_{-m}\ =\ t^{-1/2}
\varepsilon_{-m}+(t^{1/2}-t^{-1/2})\varepsilon_{m}.
\end{align}
Similarly, for $m\geq 0$,
\begin{align*}
&T^{-1}\varepsilon_{-m}\ =\ t^{-1/2}\varepsilon_{m},\\
&T^{-1}\varepsilon_{m}\ =\
(T-(t^{1/2}-t^{-1/2}))\varepsilon_{m}=t^{1/2}
\varepsilon_{-m}-(t^{1/2}-t^{-1/2})\varepsilon_{m}.
\end{align*}

\begin{lemma}\label{LEMPI}
For any $m\in \Z$, $\pi\varepsilon_{m}=\varepsilon_{1-m}$.
\end{lemma}
{\em Proof.}
Since $\pi^{2}=1$, it suffices to calculate
$\pi \varepsilon_{-m}$ for $m\le 0$.
Using that $Y\vep_m=t^{1/2}\vep_{m+1}$ (it results from the
definition of $\vep$ for such $m$),
$$\pi \varepsilon_{-m}=YT^{-1}\varepsilon_{-m}=
t^{-1/2}Y\varepsilon_{m}=
\varepsilon_{1-m}.$$
\sq

Let us apply the lemma to write down the action of $Y^{\pm1}$ on
$\varepsilon_{m},\vep_{-m}$ for $m\ge 0$:

\begin{eqnarray}
Y\varepsilon_{m}
&=&t^{1/2}\varepsilon_{m+1}\label{yepm},\\
Y\varepsilon_{-m}
&=&t^{-1/2}\varepsilon_{-m+1}+(t^{1/2}-t^{-1/2})\varepsilon_{m+1}
\label{yep-m},\\
Y^{-1}\varepsilon_{m+1}
&=&t^{-1/2}\varepsilon_{m}\label{y-epm},\\
Y^{-1}\varepsilon_{-m}
&=&t^{1/2}\varepsilon_{-m-1}-(t^{1/2}-t^{-1/2})\varepsilon_{m+1}.
\label{y-ep-m}
\end{eqnarray}
The last two formulas do not overlap as $m\ge 0$; the
pairs of formulas for the action of
$T,T^{-1}$ and $Y$
intersect (and  coincide) at $m=0$. The formulas for
the action of $Y,Y^{-1}$
are called {\em nonsymmetric Pieri rules}; they are {\em obviously}
sufficient to calculate the $\vep$\~functions (it holds
in any ranks).
These formulas are dual to (\ref{pierie},\ref{pierie1}).

Generally, the {\em technique of intertwiners} is more efficient
for calculating the $\vep$\~polynomials and their generalizations
than direct using the Pieri formulas (see, e.g., \cite{C101}).
In this example, formula (\ref{Tvepm}) is sufficient. Indeed,
for $m\ge 0$,
\begin{eqnarray}\label{vepformulas}
\vep_m & = & t^{-\frac{m}{2}}Y^m \hbox{\ \, implies \ that}\\
\vep_{-m}
&=& t^{-\frac{1}{2}}T\vep_m=t^{-\frac{m+1}{2}}T(Y^m)\notag\\
&=&t^{-\frac{m+1}{2}}(t^{\frac{1}{2}}Y^{-m}+
(t^{\frac{1}{2}}-t^{-\frac{1}{2}})\frac{Y^{-m}-Y^{m}}
{Y^{-2}-1}).\notag
\end{eqnarray}
}

\subsection{\bf Whittaker functions}
\label{sect:whitlim}
In this section, $t=q^k$ and
$|q|<1$ unless stated otherwise. We will use
the elementary difference operator $\Ga(X)=q^{1/2}X$
and also $\Ga_k(X)\equal q^{k/2}X$,

\subsubsection{\sf Whittaker limit}
Etingof states in
\cite{Et1} following Ruijsenaars \cite{Rui} that
$$
\lim_{k\to -\infty}q^{-kx}\Ga_{k}\,\mathcal{L}\,\Ga_{-k}q^{kx}
$$
becomes the so-called
$q$\~Toda (difference) operator. To be exact, they
considered the case of $A_n$.

Following \cite{ChW},
the basic limiting procedure
will be when $k$ approaches $\infty$
for $|q|<1$ ($|t|\to 0$), unless stated otherwise. Let
\begin{align}\label{remain}
\hbox{\ae}(\mathcal{L})&\equal
q^{kx}\Ga_{k}^{-1}\,\mathcal{L}\,\Ga_{k}q^{-kx},\
\t\equal R\!\!E(\l)=\lim_{k\to \infty}
\hbox{\ae}(\l),
\end{align}
where the second limit is
the {\em Ruijsenaars-Etingof procedure}.
At level of functions $F(X)$:
$$R\!\!E(F)=\lim_{k\to \infty}q^{kx}\,F(q^{-k/2}X)
=\lim_{k\to \infty}q^{kx}\Ga_{k}^{-1}(F).
$$

Later, the case $|q|>1$ ($|t|\to \infty$) will be
needed too. Then
\begin{align}\label{reprime}
\hbox{\ae}^\checkmark(\mathcal{L})&\equal
q^{-kx}\Ga_{k}^{-1}\,\mathcal{L}\,\Ga_{k}q^{kx},\
\t^\checkmark\equal
R\!\!E^\checkmark(\l)=\lim_{k\to \infty}
\hbox{\ae}^\checkmark(\l), \\
&R\!\!E^\checkmark(F)=\lim_{k\to \infty}q^{-kx}\,F(q^{-k/2}X)
=\lim_{k\to \infty}q^{-kx}\Ga_{k}^{-1}(F).\notag
\end{align}

Generally, the $R\!\!E$ procedures
require very specific functions $F$ to
be well defined. Formally, if
$\l(\Phi)=(\Lambda+\Lambda^{-1})\Phi$, then
\begin{align*}
&R\!\!E(\l)(\w )
=(\Lambda+\Lambda^{-1})\w \for \w=R\!\!E(\Phi)\\
&R\!\!E^\checkmark(\l)(\w^\checkmark)
=(\Lambda+\Lambda^{-1})\w^\checkmark \for
\w^\checkmark=R\!\!E^\checkmark(\Phi),\\
&\hbox{\, provided the existence of\, } \w ,\w^\checkmark.
\end{align*}
At level of operators,
\begin{align}\label{whitoper}
\hbox{\ae}(\mathcal{L})
&=\,\frac{X-X^{-1}}{t^{-1/2}X-t^{1/2}X^{-1}}t^{-1/2}\Ga
+\frac{tX^{-1}-t^{-1}X}{t^{1/2}X^{-1}-t^{-1/2}X}t^{1/2}\Ga^{-1}
\notag\\
&=\frac{X-X^{-1}}{X-tX^{-1}}\Ga
+\frac{t^{2}X^{-1}-X}{tX^{-1}-X}\Ga^{-1}.
\end{align}
Therefore when $|t|\to 0$,
\begin{align}\label{qToda}
\t=R\!\!E(\l)=\frac{X-X^{-1}}{X}\Ga+\Ga^{-1}=
(1-X^{-2})\Ga+\Ga^{-1}.
\end{align}

Similarly,
\begin{align}\label{whitoperprime}
\hbox{\ae}^\checkmark(\mathcal{L})
&=\,\frac{X-X^{-1}}{t^{-1/2}X-t^{1/2}X^{-1}}t^{1/2}\Ga
+\frac{tX^{-1}-t^{-1}X}{t^{1/2}X^{-1}-t^{-1/2}X}t^{-1/2}\Ga^{-1}
\notag\\
&=\frac{X-X^{-1}}{X-tX^{-1}}\Ga
+\frac{t^{2}X^{-1}-X}{tX^{-1}-X}\Ga^{-1}.
\end{align}
In this case $|q|>1$, so $|t|\to \infty$ and
\begin{align}\label{qTodaprime}
\t^\checkmark=R\!\!E^\checkmark(\l)
=\frac{X-X^{-1}}{-X^{-1}}\Ga+\Ga^{-1}=
(1-X^{2})\Ga+\Ga^{-1}.
\end{align}

One of the main results of \cite{ChW} is the formula
for the $R\!\!E$\~limits of the
{\em global symmetric
$q,t$\~spherical function}  (for arbitrary reduced root
systems; see the definitions there). In the $A_1$\~case,
the limit of $F(X;\La)$ from (\ref{fxla}) is as
follows:
\begin{align}\label{Whitsym}
&\w(X,\La)=\\
(\tga'(X)\tga'(\La))^{-1}
&\sum_{m=0}^{\infty}q^{m^{2}/4}\,\overline{P}_{m}(\La;q)\,
X^m\,\prod_{s=1}^{m}\frac{1}{1-q^{s}}\,,\notag
\end{align}
where  $|q|<1,\, \prod_{s=1}^0=1$,
$\La=q^\la$ as for $X$, $\overline{P}_m$ are
the symmetric $q$\~Hermite polynomials, to be discussed next.
It satisfies, among its other properties, the relation
$\t(\w)= (\La+\La^{-1})\w$.

Recall that $\tga'(X)\tga'(\La)$ divided by
$q^{-x^{2}-\la^{2}}$ is a $\Z/2$\~periodic function of $X$ and
of $\La$. Therefore, as far as the Toda eigenvalue problem
and other symmetries are
concerned,  $(\tga'(X)\tga'(\La))^{-1}$ can be replaced
by $q^{x^{2}+\la^{2}}.$ In this paper, we prefer to make
all functions in terms of Laurent variables, $X$ and $\La$;
this explains our choice to use $\tga'$ rather than $q^{-x^2}$.

We note that the series from (\ref{Whitsym}) was
introduced and discussed in \cite{Sus} (for $A_1$ only),
however, without
the Gaussians and without any reference to the Whittaker
theory. It appeared there as a quadratic exponential function
and as a quadratic generating function for one-dimensional
$q$\~Hermite polynomials.
\smallskip

When $|q|>1$, the $R\!\!E^\checkmark$\~procedure results in
\begin{align}\label{Whitsymprime}
&\w^{\,\checkmark}(X,\La)=\\
(\tga(X)\tga(\La))^{-1}
&\sum_{m=0}^{\infty}q^{-\frac{m^{2}}{4}}\,
\overline{P}_{m}(\La;q^{-1})\,
X^{-m}\,\prod_{s=1}^{m}\frac{1}{1-q^{-s}}\,,\notag
\end{align}
where $\tga$ is obtained from $\tga'$ by the
formal substitution $q\mapsto q^{-1}$.

The function $\w^{\,\checkmark}(X,\La)$
is obviously $\w(X,\La)$ upon the action of the automorphism
of $\C[q^{\pm 1/4}][X^{\pm 1}]$ defined by
$X^\checkmark= X^{-1},$ and $(q^{1/4})^\checkmark=q^{-1/4}$.
Since $x^\checkmark=x$ due to
$X=q^x$, the image $\Ga^\checkmark$ of $\Ga$,
sending $x\mapsto x+1/2$,
coincides with $\Ga$.
Thus $\t^\checkmark(\w^{\,\checkmark})=
(\La+\La^{-1})\w^{\,\checkmark}$,
as it is supposed to be.
\smallskip

Note that the formulas above are valid for any
$|q|\neq 1$ upon a straight algebraic transformation from
$q$ to $q^{-1}$. Moreover, the analytic continuation
to $|q|=1$ is possible for almost all such $q$ (roots
of unity must be excluded, but not only them); see below.
\smallskip

\subsubsection{\sf Harish-Chandra expansion}
A systematic analytic theory of the global functions will be
a subject of further papers; however, the following introduction
to the Harish-Chandra $q,t$\~theory seems quite relevant here.
The corresponding nonsymmetric theory will be presented in a
continuation of this paper. It is actually the best way to
proceed for arbitrary root systems; for $A_n$, the constructions
in the symmetric setting are relatively straightforward.

For $|q|<1$ and $|X|>|q|^{\,1/2\,}$ the asymptotic expansions
and the corresponding Harish-Chandra decomposition  read
in the $q$\~Whittaker case as follows:

\begin{align}\label{Whitexpan}
&\tga'(X)\tga'(\La)\,\w(X,\La)\\
=
\lan \overline{\mu}\ran&\,\overline{\si}(\La^{-1})\,\tga'(X\La)\,
\sum_{j=0}^\infty
q^{j}X^{-2j}\prod_{s=1}^j\,\frac{1}{(1-q^s)(1-q^s\La^{2})}
\notag\\
+\lan \overline{\mu}\ran\,&\overline{\si}(\La)\tga'(X\La^{-1})\,
\sum_{j=0}^\infty
q^{j}X^{-2j}\prod_{s=1}^j\frac{1}{(1-q^s)(1-q^s\La^{-2})},\notag\\
\hbox{where}&\ \
\overline{\si}(\La)\ =\ \prod_{j=0}^\infty (1-q^j \La^2)^{-1},\
\lan \overline{\mu}\ran\ =\ \prod_{j=1}^\infty (1-q^{j})^{-1}
\notag
\end{align}
are the $q$\~Whittaker version of the
Harish-Chandra $c$\~function and the constant term of
$\overline{\mu}$ from (\ref{mutildemuoneh}) below.
This identity follows from
the Whittaker part of the paper \cite{ChW}. The corresponding
expansions exist for arbitrary (reduced) root systems (in the twisted case),
but the explicit formulas for the coefficients can be involved.
\smallskip

Using that the terms in the right-hand solve the
$q$\~Toda eigenvalue problem, they can be meromorphically
continued to all values of $X$. The difference equation
in terms of $\La$ is also known (see \cite{ChW}), so the
right-hand side can be extended meromorphically to
all $X,\La$.

One can expect to exclude $q,\La$ such that
$q^s\La^{\pm 2}\neq 1$ for any $s\in \N$, but the poles from the
both terms in the right-hand side will cancel
each other at these points. Thus it is not necessary.
Generally speaking (for arbitrary root systems), the
cancelation of the poles in the Harish-Chandra decomposition
formulas is difficult to check directly. The existence of
the global function makes it immediate.

The inequalities $|X|>|q^{1/2}|$ guarantee the absolute
convergence of the right-hand side (the expansions) for $|q|<1$.
However, they are not needed for the
left-hand side (the global function), which is analytic
for any $X,\La$
when $|q|<1$ or when $|q|>1$.
\smallskip

Note that the convergence of the
summations in the right-hand side of (\ref{Whitexpan}) becomes
significantly better for $|q|>1$; $X,\La$ can be arbitrary
for such $q$ (as in the left-hand side). The same essentially
holds for the multipliers
$\lan \overline{\mu}\ran\,\overline{\si}(\La^{\pm 1})$. Indeed,
they can be redefined in terms of the $q$\~exponential series; then
similar arguments can be applied.
\medskip

The left-hand (global) side of
(\ref{Whitexpan}) can be extended even to
{\em almost all} points at the unit circle:
\begin{align}\label{ucircw}
&|q|=1 \hbox{\, provided }
|X\La^{\pm 1}|>1/R(q),\\
&\where R(q)\equal\liminf_{\,m\to \infty}\,
|1-q^m|^{1/m}.\notag
\end{align}

Here we must of course avoid the roots of unity, where $R(q)=0$,
but the essence of this approach is that it is not sufficient;
a greater set of points must be avoided.
See, e.g., \cite{Lub}, especially, (1.17)-(1.19) there.
The theory of such analytic continuations to the unit circle
is classical, due to Hardy-Littlewood and others.

The analytic continuation to the unit circle (at almost all
points) is, generally, only with respect to the non-tangent 
directions. However, the following sequences $\{q_n\}$ 
tending to $q$ are allowed. If $|q_n|=1$ then we must assume
that $|1-q_n^m|^{1/m}>\ep$ for an arbitrarily small $\ep>0$,
almost all $m>0$ and for sufficiently large $n$.

\medskip
\rmk
The right-hand side of (\ref{Whitexpan}) is
a weighted sum of the asymptotic series as $|X|\to \infty$ 
for the $q$\~Toda eigenvalue problem and that for
$\La\mapsto \La^{-1}$. This series solves
the Toda eigenvalue problem and its coefficients
can be uniquely determined from it.
Recall that the exact convergence condition for it
is $|X|>|q|^{1/2}$ for $|q|<1$, which 
becomes much better (any $X$) upon adding its
$\La^{-1}$\~counterpart. 

This series is an instance of the formula from \cite{GiL}
(the $GL(n)$\~case), representing a certain generating 
function in the quantum $K$\~theory of flag varieties. 
See the Appendix, Section \ref{sec:qKth}. We note that only
the asymptotic series alone appears in \cite{GiL},
not the weighted summation from (\ref{Whitexpan}).

Thus the Harish-Chandra formula (\ref{Whitexpan}) for 
the global Whittaker function connects that Givental-Lee theory 
with the algebraic geometry of affine
Schubert varieties (encoded in the level one affine Demazure character
formulas). The limit $q\to 1$ of (\ref{Whitexpan}) coincides with
the corresponding $c$\~weighted summation formula for the
classical (real) Whittaker function (see \cite{GW},\cite{Wa});
generally, the summation is over the Weyl group.
\sq
\smallskip

The spherical $q,t$\~generalization of (\ref{Whitexpan}) is actually
simpler to establish. The following
theorem is not too difficult to justify in the rank one case.
Generally, the nonsymmetric theory helps significantly.

\begin{theorem}\label{HCHEXP} For the
function $F(X;\La)$ from (\ref{fxla}), let us assume
that $|q|<1$  and $|X|<|t|^{1/2}|q|^{-1/2}$. Then
\begin{align}\label{Sphexpan}
&\frac{\tga'(X)\tga'(\La)}{\tga'(t^{1/2})}
\,F(X,\La)\\
=\lan \mu\ran\si(\La)&\tga'(X\La t^{-1/2})\,
\sum_{j=0}^\infty
(\frac{q}{t})^{\!{}^{j}}X^{2j}\prod_{s=1}^j\frac{(1-tq^{s-1})
(1-q^{s-1}t\La^{-2})}{(1-q^s)(1-q^s\La^{-2})}
\notag\\
+\lan \mu\ran\si(\La^{-1})&\tga'(X\La^{-1}t^{-1/2})\,
\sum_{j=0}^\infty
(\frac{q}{t})^{\!{}^{j}}X^{2j}\prod_{s=1}^j\frac{(1-tq^{s-1})
(1-q^{s-1}t\La^{2})}{(1-q^s)(1-q^s\La^{2})},
\notag
\end{align}
where $\si(\La)=\prod_{j=0}^\infty
\frac{1-tq^j \La^2}{1-q^j \La^2}$
is the $q,t$\~generalization of the Harish-Chandra $c$\~function;
$\lan \mu\ran$ is from (\ref{mucontone}). \sq
\end{theorem}

This identity follows from the $q,t$\~part of paper
\cite{ChW} and a straightforward calculation of
the expansion coefficients. The terms in the right-hand
satisfy the Macdonald  eigenvalue problem, so they
can be extended meromorphically to all $X,\La$.

Formula (\ref{Whitexpan})
is actually the result of the $R\!\!E$ limiting procedure applied
to (\ref{Sphexpan}); it is an instructional calculation.
We replace $X$ by $X^{-1}$ (using that $F$ is symmetric)
and then use that
$$
R\!\!E\Bigl(\frac{\tga'(X)\tga'(\La)}{\tga'(t^{1/2})
\tga'(X^{-1}\La t^{-1/2})}\Bigr)=
\frac{\tga'(X)\tga'(\La)}{\tga'(X^{-1}\La)}.
$$
It suffices to take here $q^{-x^2}$ instead of
$\tga'(X)$ (and for the other arguments).
Then this formula becomes quite obvious. The limits of
the other terms in (\ref{Sphexpan}) can be obtained by the
straight substitution $t\mapsto 0$.

The left-hand side here is an analytic function for all $X,\La$
when $|q|<1$, so its convergence is significantly
better than for the {\em classical basic hypergeometric
series} on the right-hand side.
The parameter $t$ is assumed sufficiently general to avoid the zeros
of $\lan P_n,P_n\ran$. The inequality $|X|<|t|^{1/2}|q|^{-1/2}$
from the theorem ensures the convergence of the right-hand side
for $|q|<1$.

Note that when $|q|>1$, the summations in the
right-hand side still converge if we impose the condition
$|X|<|q|^{1/2}|t|^{-1/2}$. Thus
the right hand side can be used to connect the domains
$|q|<1$ and $|q|>1$ (under the above inequality).
Here the function $\si(\La^{\pm 1})$ must be
redefined using the $q$\~exponential functions.
This continuation actually goes through the unit circle.

Concerning the left-hand side, it can be continued analytically to
\begin{align}\label{ucircwt}
|q|=1 \hbox{\, provided }
|X^{\pm 1}\La^{\pm 1}|<|t|^{1/2}R(q).
\end{align}
Compare with (\ref{ucircw}).
The continuation is for sequences $\{q_n\}$ approaching $q$
such that $|1-q_n^m|^{1/m}>\ep$ for an arbitrarily small $\ep>0$,
almost all $m>0$, and sufficiently large $n$.
\smallskip

Recall that $F$ in the left hand-side
remains a solution of the same
difference equation $\l (F)=(\La+\La^{-1})F$ upon the formal
star-conjugation
$$
q\mapsto q^{-1},\, t\mapsto t^{-1},\, X\mapsto X^{-1},\,
\La\mapsto \La^{-1}.
$$
This readily provides a formula defined for $|q|>1$. It is
straightforward to calculate its expansion decomposition.

The connection between the regions  $|q|<1$ and
$|q|>1$  in the $q,t$\~case
is significantly  different from that
in the Whittaker case. In the Whittaker case,
the conjugation changes the corresponding
difference equation (it does not in the $q,t$\~case);
however, the defining series of each of these two functions
(the left-hand size of (\ref{Whitexpan})\,) can be extended
from the region $|q|<1$ to the region
$|q|>1$ simply by recalculating the corresponding series
in terms of $q^{-1}$.
\smallskip

\rmk
Formula (\ref{Sphexpan}) (without explicit expressions for
coefficients) was announced by Jasper
Stokman for $GL(n)$ in his lecture (June, 2009). A complete
theory of the $q,t$\~version of the Harish-Chandra expansion
and decomposition theory for any reduced root systems 
(the twisted case) is in his recent \cite{Sto2}.
Similarly to the classical differential theory,
the difference decomposition formula results from the following
two ingredients:

(a) the calculation of the corresponding asymptotic limit,
which is the $c$\~function (\cite{ChW} in the $q,t$\~case),

(b) the coefficient-wise existence and uniqueness of the
asymptotic series and its meromorphic continuation.
\smallskip

The formal uniqueness is essentially sufficient to establish
the existence of this series in a neighborhood of infinity.
The exponential growth of the coefficients results from the
recurrence relations for them (the uniqueness).
Then the corresponding differential or difference
equations can be used for the meromorphic continuation.

Stokman's approach  is based on the Cherednik-Matsuo
map to/from the QAKZ-equation, where the calculation of the
asymptotic coefficients is very similar to that in the
differential AKZ-case (paper \cite{ChA} and previous first author's
papers on the $r$\~matrix Knizhnik-Zamolodchikov equations).

An alternative approach to the
coefficient-wise existence and uniqueness of the
asymptotic expansion is in deducing this
fact from the differential theory, namely,
from that due to Heckman and Opdam; see \cite{HO},
which is the limiting case $q\to 1$ of the $q,t$\~theory.
The standard deformation argument is used here.
This line was suggested by Opdam and the first author.
\smallskip

Paper \cite{MS} (the $GL(n)$\~case) and the papers
\cite{Me, Sto2} (the twisted case, arbitrary reduced root systems)
are devoted to the existence of the asymptotic series, their
convergence and symmetries. We remark that paper
\cite{ChW} was written in the twisted case too.

The advantage of the approach due
to van Meer and Stokman (vs. that based on the deformation to
the differential theory) is that the convergence
of the asymptotic series in a neighborhood of infinity readily
results in its global meromorphic continuation. The QAKZ
is a {\em difference} system of equations, which readily
provides the desired meromorphic continuation (in contrast to
the differential theory).

The Macdonald eigenvalue problem (used directly) is expected
to provide the same, but the corresponding tools are not properly
developed at the moment.

A disadvantage of this approach is that the Cherednik-Matsuo
map is the symmetrization of the vector-valued solutions
of QAKZ (of dimension $|W|$), which makes it practically
impossible to use for explicit calculation of the asymptotic
coefficients (even in the simplest cases). Finding the
fundamental solution of AQKZ is generally significantly
more difficult problem than direct finding the symmetric
solution of the Macdonald eigenvalue problem, algebraically
and analytically (in spite of the fact that Cherednik-Matsuo
maps are isomorphisms).

The Whittaker reduction of the $q,t$\~version of the
Harish-Chandra asymptotic decomposition formula and its
connection with the Givental-Lee theory was announced by
the first author. It also establishes the connection with
\cite{GW} and \cite{Wa} and with the classical $p$\~adic
Whittaker theory.

The first author is grateful to Jasper Stokman for the
information about his ongoing research and Eric 
Opdam for useful discussions.
\smallskip

\subsubsection{\sf Nonsymmetric q-Hermite polynomials}
\label{sect:QHermite}
For an $E$\~polynomial $E_{n}$,
let us define  its two limits:
$$
\overline{E}_{n}=\lim_{t\to 0}E_{n}
\and \overline{E}^\dag_{n}=\lim_{t\to\infty}E_{n}.
$$
Both limits exist (for instance, use the explicit formulas or the
intertwining operators) and are closely connected to each other.

\begin{proposition}
For $n\geq 0$,
\begin{align}\label{tildee}
\overline{E}^\dag_{-n}\,
=\,\bigl(q^{\frac{n}{2}}\overline{E}_{-n}(Xq^{\frac{1}{2}})\bigr)
\Big|_{q\to q^{-1}},\
\overline{E}^\dag_{n}\,
=\,\bigl(q^{-\frac{n}{2}}\overline{E}_{n}(Xq^{\frac{1}{2}})\bigr)
\Big|_{q\to q^{-1}}.
\end{align}
\end{proposition}
\sq

The polynomials $\overline{E}_n$ are called {\em nonsymmetric
(continuous)
$q$\~Hermite polynomials} (see \cite{ChW} and references
therein; they are considered there for arbitrary reduced root
systems). Their symmetrizations are the classical $q$\~Hermite
polynomials.
\smallskip

\rmk
Upon the substitution $X\mapsto X^{-1}$,
the polynomials $\overline{E}_n$ are directly
related to the Demazure characters
of level one Kac-Moody integrable modules, thus
are closely connected with the geometry of the affine
Schubert varieties (through the Kumar-Mathieu formula).

This connection is from \cite{San} in the $GL(n)$\~case;
for arbitrary root systems, it is from \cite{Ion}
(the twisted case).
\sq
\smallskip

More systematically, let us define
\begin{align}\label{overlinet}
\overline{T}\,\equal\, \lim_{t\to 0}t^{1/2}T=
\frac{1}{1-X^{2}}\circ (s-1),\
\overline{T}(\overline{T}+1)=0.
\end{align}
Using intertwiners, $\overline{E_0}=1$,
\begin{align}\label{PiTplus1}
&\overline{E}_{1+n}\,=\,q^{n/2}\Pi \overline{E}_{-n},\\
&\overline{E}_{-n}\,=\,(\overline{T}+1)\overline{E}_{n}\notag
\end{align}
for $n\ge 0$; the raising operator $\Pi=X\pi$ was defined
in (\ref{Piformula}).

From the divisibility condition
$\overline{T}+1=(s+1)\cdot\{\,\}$, we obtain that
$\overline{E}_{-n}$ is symmetric ($s$\~invariant) and
$\overline{P}_{n}=\overline{E}_{-n}$ for $n\ge 0.$

Explicitly,
\begin{align*}
\overline{E}_{-n-1}\,=\,((\overline{T}+1)\Pi q^{n/2})
\overline{E}_{-n},\\
(\overline{T}+1)\Pi
\,=\,\frac{X^{2}\Ga^{-1}-X^{-2}\Ga}{X-X^{-1}}.
\end{align*}

The bar-Pieri rules read as follows:
\begin{align}\label{pienilp+}
&X^{-1}\overline{E}_{-n}=\overline{E}_{-n-1}-\overline{E}_{n+1}
\ (n\ge 0),\\
&X^{-1}\overline{E}_n=(1-q^{n-1})\overline{E}_{n-1}+q^{n-1}
\overline{E}_{1-n}\ (n\ge 1),\notag\\
\label{pienilp-}
&X\overline{E}_{-n}=(1-q^n)\overline{E}_{1-n}+\overline{E}_{n+1}
\ (n\ge 0),\\
&X\overline{E}_n=\overline{E}_{n+1}-q^{n}
\overline{E}_{1-n}\ (n\ge 1).\notag
\end{align}

Let $\overline{Y}=\pi \overline{T}=\lim_{t\to 0}t^{1/2}Y$.
Recall that
\begin{eqnarray*}
YE_{n}=\left\{\begin{array}{ccc}t^{-1/2}q^{-n/2}E_n, &  & n>0, \\
t^{1/2}q^{-n/2}E_{n}, &  & n\leq 0.\end{array}\right.
\end{eqnarray*}

In the limit,
\begin{eqnarray}\label{Ynilp+}
\overline{Y}\,\overline{E}_{n}
=\left\{\begin{array}{ccc}q^{-|n|/2}\overline{E}_{n}, &  & n>0, \\
0, &  & n\le 0.\end{array}\right.
\end{eqnarray}

Since $\overline{Y}$ is not invertible, we need to introduce
$$
\overline{Y}\,'=\lim_{t\to 0}t^{1/2}Y^{-1}=
\lim_{t\to 0}t^{1/2}T^{-1}\pi=\overline{T}\,'\pi
$$
for $\overline{T}\,'=\overline{T}+1$. Then
$\overline{Y}\,\overline{Y}\,'=0=\overline{Y}\,'\,\overline{Y}$ and
\begin{eqnarray}\label{Ynilp-}
\overline{Y}\,'\,\overline{E}_{n}
=\left\{\begin{array}{ccc}q^{-|n|/2}\overline{E}_{n},
&  & n\le 0, \\
0, &  & n>0.\end{array}\right.
\end{eqnarray}
Finally,
\begin{align*}
\overline{\l}=\lim_{t\to 0}t^{1/2}\l=\overline{Y}\,'
+\overline{Y}
\,=\,\frac{1}{1-X^{2}}\Ga+\frac{1}{1-X^{-2}}\Ga^{-1}
\end{align*}
and $\overline{\l}\,\overline{P}_{n}=
q^{-n/2}\overline{P}_{n},\, n\ge 0$ (see (\ref{eqn:LP})); recall that
$\overline{P}_n=\overline{E}_{-n}$.

\subsubsection{\sf The series for the mu-function}
The Hermite-type degeneration of the function
$\mu$ from (\ref{mutildemuone}) is as follows:
\begin{align}\label{mutildemuoneh}
&\overline{\mu}(X;q)=\prod_{j=0}^\infty
(1-q^jX^2)(1-q^{j+1}X^{-2})=\\
&\prod_{j=1}^\infty (1-q^j)^{-1}\,
\sum_{n=-\infty}^\infty (-1)^nq^{\frac{n^2-n}{2}}X^{2n},\notag\\
&\lan \overline{\mu}\ran=\prod_{j=1}^\infty (1-q^j)^{-1},\
\om(\overline{\mu})=\overline{\mu}(Xq^{-1/2};q)=
(-X^2q^{-1})\overline{\mu}.\notag
\end{align}
Also, $\overline{\mu}(X^{-1})=-X^{-2}\overline{\mu}(X).$
Compare with $\om(\tga')=(q^{-1/4}X)\tga'$, where
\begin{align}\label{gauexp}
\tga'=&\sum_{n=-\infty}^\infty q^{n^2/4}X^n\\
=&\prod_{j=1}^\infty (1-q^{j/2})(1+q^{\frac{2j-1}{4}}X)
(1+q^{\frac{2j-1}{4}}X^{-1}).\notag
\end{align}

The function $\tga'\overline{\mu}$
serves the inner product with the Gaussian.
One has from (\ref{mehmacone}):
\begin{align}\label{mehmaconew}
&\lr \tga'\overline{\mu}_\circ\rr=
\prod_{j=1}^\infty (1-q^j),\
\lr \tga'\overline{\mu}\rr=1.
\end{align}
The complete expansion is as follows:
\begin{align}\label{mugaone}
&\tga'\overline{\mu}\ = \
\sum_{n=0}^\infty q^{n(n+2)/12}(X^{n+2}-X^{-n}),\where
n\neq 2 \hbox{ mod } 3.
\end{align}
More explicitly, the coefficient of $(X^{n+2}-X^{-n})$
here equals
$$
q^{\frac{(3m+2)m}{4}}\hbox{ for } n=3m,\ q^{\frac{(3m+1)(m+1)}{4}}
\hbox{ for } n=3m+1,\
0 \hbox{ otherwise}.
$$
Ignoring $n\neq 2 \hbox{ mod }3$ and substituting $n\mapsto 2n$,
$q\mapsto q^{3/2}$, the summation in (\ref{mugaone}) becomes that
from (\ref{mutildemuoneh}).

Adding the $\overline{E}$\~polynomials, the limit of
(\ref{f81ee}) as
$k\to \infty$ ($t\to 0$) is
\begin{align}\label{f81eew}
&\lr \overline{E}_n \overline{E}_m\tga'\overline{\mu}\rr=
\lim_{k\to \infty} q^{\frac{m^2+n^2+2k|m|}{4}}E_m(q^{n_\#})=\\
&q^{\frac{(m-n)^2}{4}}
\for \sgn(n)+\sgn(m)<2 \and 0 \hbox{\, otherwise\,}.\notag
\end{align}
Indeed, the term
$q^{\frac{k|m|}{2}}E_m(q^{n_\#})$ is
nonzero in the limit if and only if
$E_m$ contains the monomial $X^{-\hbox{\tiny sgn}(n)|m|}$.
This holds unless $n>0$ and $m>0$. Then the
coefficient of this monomial will be always $1$
in the bar-limit; we use that $\overline{E}_{-n}=\overline{P}_n$
for $n>0$.

The particular cases $0\le n,m \le 1$ are immediate from
(\ref{mugaone}).
\medskip

\setcounter{equation}{0}
\section{\sc Nil-DAHA}
\subsection{\bf Key definitions}
A systematic theory of $q$\~Hermite polynomials
and global Whittaker functions begins
with the following definition of the {\em
nil-DAHA\,} which can be
readily adjusted to arbitrary (reduced) root systems.

\begin{definition}\label{TWONIL}
(i) The {\sf nil-DAHA}\, $\overline{\HH}$ is
generated by $T,\pi, X^{\pm 1}$ over the ring $\C[q^{\pm 1/4}]$
with the defining relations: $T(T+1)=0$,
\begin{align}\label{nildahax}
\pi^2=1,\ \pi X\pi=q^{1/2}X^{-1},
\ TX-X^{-1}T=X^{-1}.
\end{align}
Equivalently, $TX=X^{-1}T'$
for  $\,T'\equal T+1\,$.
Setting $\,Y\equal\pi T\,$ and $\,Y'\equal T'\pi\,$ one has:
$\,T'Y=Y'T, TY'=0=YT'.$

(ii) Similarly, one can define
$\overline{\HH}^{\,\vph}=\C[q^{\pm 1/4}]
\lan T,\breve{\pi}, Y^{\pm 1}\ran$
subject to $TT'=0$ for $T'=T+1$\, and
\begin{align}\label{nildahay}
&\breve{\pi}^2=1,\ \breve{\pi} Y\breve{\pi}=
q^{-1/2}Y^{-1},\ TY^{-1}=YT'.
\end{align}
Setting $X\equal\breve{\pi} T',\,
X'\equal T\breve{\pi},$ one has:
$TX=X'T',\ T'X'=0=XT.$

(iii) The algebra $\,\overline{\HH}^{\,\vph}\,$ is the image of
the algebra
$\,\overline{\HH}\,$ under the
anti-isomorphism
\begin{align}\label{vphnil}
\varphi: T\mapsto T,\, \pi\mapsto\breve{\pi},\, X\mapsto Y^{-1}.
\end{align}
Correspondingly, $\varphi: Y\mapsto X', Y'\mapsto X$.
There is also an isomorphism $\si:$
$\,\overline{\HH}\,\to\,\overline{\HH}^{\,\vph}\,$ sending
\begin{align*}
&\si: T\mapsto T,\ \,X\mapsto Y^{-1},\ \,\pi\mapsto \breve{\pi},\\
&\si:\, Y\mapsto \breve{\pi} T,\ \,
Y'\mapsto T'\breve{\pi}.
\end{align*}

(iv) The automorphism $\tau_+$ fixing $T,X$ and sending
$Y\mapsto q^{-1/4}XY$ acts in $\,\overline{\HH}\,$.
Correspondingly,
$\tau_-\equal\varphi\tau_+\varphi^{-1}$ acts in
$\overline{\HH}^{\,\vph}$
preserving $T,Y$ and sending \,
$X\mapsto q^{1/4}YX$. One has
the relations
\begin{align}\label{sitausi}
&\si\tau_+\ =\ \tau_-^{-1}\si, \ \,
\si\tau_+^{-1}\ =\ \tau_-\si,
\end{align}
matching the identity from (\ref{tautautau}) in the
generic case. \sq
\end{definition}
Both algebras, $\,\overline{\HH}$ and $\overline{\HH}^{\,\vph}\,$,
satisfy the PBW Theorem because they are limits of $\HH$ (to be
discussed  below in detail),
so the $q,t$\~algebra is their {\em flat} deformation.
It includes the case when $q$ is a root of unity.
\comment{ However,
roots of unity must be avoided
in the construction of the $q$\~Hermite polynomials.
}

\subsubsection{\sf Polynomial representation}
Any $H\in \overline{\HH}$ can be represented as
$H=\sum_i c_iX^i h_i$ for constant $c_i (i\in \Z)$ and
$h_i\in \overline{\h}_Y\equal<T,\pi>$; this readily follows from
the defining relations. Moreover, this representations
is unique since $\overline{\HH}$ is a limit of $\HH$,
where the PBW Theorem holds. The existence of the
above representations automatically guarantees their
uniqueness. Furthermore, any element
$h\in \overline{\h}_Y$ can be uniquely
expressed as a linear combination of
$$
Y^m\pi, (Y')^{m+1}\pi, Y^m, (Y')^{m+1} \where m\ge 0.
$$

The bar-polynomial representation $\overline{\mathscr{X}}$ is
$\hbox{Ind}_{\overline{\h}_Y}^{\overline{\HH}}\C^0_+$
$=\C[q^{\pm 1/4}][X^{\pm 1}]$, where
$\C^{0,1}_{\,\pm}$ are one-dimensional
representations of $\overline{\h}_Y$ defined as follows:
\begin{align}\label{odimr}
&T'(1)=\ep,\ T(1)=\ep-1,\
\pi(1)=\pm 1,\\
& Y(1)=\pm(\ep-1),\ Y'(1)=\pm\ep \for  \ep=0,1.\notag
\end{align}

The bar-formulas discussed above give
an explicit description of the {\em polynomial}
(or {\em bar-polynomial\,})
representation of $\overline{\HH}$ in $\overline{\mathscr{X}};$
recall that $T,\pi, X^{\pm 1},Y,Y'$ are mapped to the operators
$\overline{T},\pi, X^{\pm 1}, \overline{Y},\overline{Y}\,'$.
\smallskip

\subsubsection{\sf The order function}
An important feature of the nil-case is that the nil-DAHA
are filtered algebras with Weyl-type algebras as their
associated graded algebras.
It significantly simplifies the PBW Theorem
and related issues.
Let the {\em order} be $m$ for the following elements
\begin{align}\label{ordmonom}
X^l Y^m\pi, X^l (Y')^m\pi (m>0), X^l Y^m, X^l (Y')^{m}(m>0),
\end{align}
where $l\in \Z, m\in \Z_+$.
The order ord$(H)$ of $H\in \overline{\HH}$ will be
defined as the maximum of orders of the terms in
the linear decomposition of $H$ with respect to
(\ref{ordmonom}).

\comment{Actually,
here it is not necessary to assume that such decomposition
is unique (which is true). One can take the minimum order
over all such linear decompositions of $H$ (if there were
more than one). The fact that the monomials (words)
above cannot belong to lower subspaces of the
ord-filtration will follow automatically.}

\begin{proposition}\label{ordbar}
The order above satisfies
$$\hbox{ord}(H_1 H_2)\le \hbox{ord}(H_1)+
\hbox{ord}(H_2).$$
The associated graded algebra
$\hbox{\rm gr}\overline{\HH}$ is isomorphic to
$$
\overline{\w}\equal <X,X^{-1},Y,Y',\pi>\,
\hbox{\ over\ }\, \C[q^{\pm 1/4}]
$$
subject to nil-versions of the Weyl-type relations
\begin{align}\label{nil-weyl}
YY'=0=Y'Y,\ XY=q^{1/2}YX \hbox{\ extended\ by\ }\\
 \pi^2=1,\ \pi X\pi= q^{1/2} X^{-1},\ \pi Y\pi= Y'.\notag
\end{align}
They result in $T^2=0$ for $T=\pi Y=Y'\pi$
and in $XY'=q^{-1/2}Y'X$. The corresponding grading is
given by $\hbox{ord}(Y)=1=\hbox{ord}(Y')$,
$\hbox{ord}(X)=0=\hbox{ord}(\pi).$
\end{proposition}
{\it Proof.} The relations $Y'T=T'Y, Y'=T'\pi$ become
$YT=TY$ and $Y'=T\pi$  in $\hbox{\rm gr}\overline{\HH}$
since $T'=T$ in the latter.
Similarly, $TX=X^{-1}T'$ becomes $XY=q^{1/2}YX$ upon
the substitution $T=\pi Y=T'$.\sq
\smallskip

\subsubsection{\sf Invariant symmetric forms}
Recall (\ref{mutildemuoneh}):
\begin{equation}\label{mutildemuonehh}
\overline{\mu}(X;q)=\prod_{j=0}^\infty
(1-q^jX^2)(1-q^{j+1}X^{-2}).
\end{equation}

The inner products
\begin{align}\label{innerpw}
&\lan f,g\ran_\circ=\lan fg\mu_\circ \ran,\
\lan f,g\ran'=\lan fg \tga' \mu \ran.
\end{align}
from (\ref{innerp}) can obviously be used in the
bar-polynomial representation $\overline{\mathscr{X}}$.
The remaining form, which involves the conjugation, will be
addressed later.
The kernels of these bilinear forms are the bar-limits,
$\overline{\mu}_\circ$ and $\tga\overline{\mu}$ from
(\ref{mutildemuoneh}) and (\ref{mugaone}).

Note that $ \lan 1,1\ran_\circ=1=\lan 1,1\ran'$ and
both forms are actually well-defined for any $q$ (including
roots of unity).
One can see this directly from (\ref{mutildemuoneh})
for $\lan f,g\ran_\circ$. The second inner product
was calculated in (\ref{f81eew}):
\begin{align}\label{f81eeww}
&\lr \overline{E}_n, \overline{E}_m\rr'=\\
&q^{\frac{(m-n)^2}{4}}
\for \sgn(n)+\sgn(m)<2 \and 0 \hbox{\, otherwise\,}.\notag
\end{align}

The corresponding anti-involutions of
the algebra $\overline{\HH}\ni H$ are
those from (\ref{diamondef},\ref{psidef}):
\begin{align}\label{innerinvw}
&\lan f,H(g)\ran_\circ=
\lan H^{\diamond}(f),g\ran_\circ,\
\lan f,H(g)\ran'=\lan H^\psi(f), g\ran',
\end{align}
where $f,g\in \overline{\mathscr{X}},\, H^\psi=\psi(H)$.
\medskip

\subsection{\bf Using the conjugation}
The nil-DAHA is a limit of $\HH$ as $t \to 0$.
The anti-involution $\ast$ plays an important role in
the $q,t$\~theory. Since $t^\ast=t^{-1}$, its nil-counterpart
requires considering the limit as $t \to \infty$ as well.

\subsubsection{\sf The limit at infinity}
We use $\overline{\HH}^\dag$ to denote the resulting
algebra. Explicitly, we have
$\overline{\HH}^\dag = <T^\dag, \pi, X^{\pm 1}>$
subject to the
defining relations:
\begin{align}\label{HHbardag}
&(T^\dag-1)T^\dag = 0,\ \pi^2 = 1,\ \pi X \pi =
q^{1/2}X^{-1},\\
&T^\dag X=X^{-1}(T^\dag-1)=X^{-1}(T^\dag)'
\for (T^\dag)'\equal T^\dag-1.
\notag
\end{align}
We set $T^\dag=\lim_{t\to \infty} t^{-1/2}T$ in the limit from
$\HH$. The polynomial representation
$\mathscr{X}$ survives in this limit; we denote the resulting
$\overline{\HH}^\dag$\~module by $\overline{\mathscr{X}}^\dag$. The
operators $\pi,\ X^{\pm 1}$ remain the same, but $T^\dag$ acts by
\begin{equation}
\overline{T}^\dag = s + \frac{1}{X^2-1}(s-1).
\end{equation}

The anti-involution $\ast$ naturally becomes the
anti-isomorphism sending
\begin{align}\label{astnew}
&\overline{\HH} \ni H \mapsto H^\ast \in
\overline{\HH}^\dag,\notag\\
&T \mapsto (T^\dag)'=T^\dag - 1,\ \pi \mapsto \pi,
\ X^{\pm 1} \mapsto X^{\mp 1},\ q \mapsto q^{-1}.
\end{align}
Its inverse will be denoted by the same symbol;
we note that $(T^\dag)^\ast=T'=T+1$
under this map.
\smallskip

\subsubsection{\sf Connection maps}
These two algebras are connected by the following
$\C[q^{\pm 1/4}]$\~linear isomorphisms:
\begin{align} \label{nudag}
&\nu_{\pm} : \overline{\HH} \rightarrow
\overline{\HH}^\dag :
T \mapsto -T^\dag,\ \pi \mapsto \pm\pi,\ X \mapsto X,\\
&\beta : \overline{\HH} \rightarrow
\overline{\HH}^\dag : T \mapsto
X^{-2}(1-T^\dag),\ \pi \mapsto \pi,\ X \mapsto X.
\label{betadag}
\end{align}

The map
$\beta$ is a adjusted to the polynomial representations in the
following sense: for any $f \in \C[q^{\pm 1/4}][X^{\pm 1}]$
and $H \in
\overline{\HH}$ we have
\begin{align}\label{betaction}
&\overline{\mathscr{X}}\ni \overline{H}(f) \, =\,
\overline{\beta(H)}(f) \in \overline{\mathscr{X}}^\dag,
\hbox{\ \ for instance,}\\
\overline{T}=
\frac{1}{1-X^{2}}&(s-1) =
X^{-2}(1-\overline{T}^\dag)=
X^{-2}(1-s + \frac{1}{1-X^2}(s-1)).\notag
\end{align}
See (\ref{overlinet}); we continue to use bar for the
operators acting in the standard polynomial representation.

We note that $\nu_{+}^{-1}\beta$ and $\beta\nu_{+}^{-1}$ are the
automorphisms of $\overline{\HH}$ and $\overline{\HH}^\dag$ given by
conjugation by $X^{-1}$ followed by $\tau_+^2$.
For instance, in $\overline{\HH}$ we have the identity
\begin{equation}
X^{-1}TX = X^{-2}T + X^{-2},
\end{equation}
and by applying $\nu_+$ to the right-hand side we
obtain the image of $T$ under $\beta$
($T, X$ are fixed under $\tau_+$, which we need here for $\pi$ only).

Finally, we define an isomorphism
\begin{align}\label{etaiso}
&\eta : \overline{\HH} \rightarrow \overline{\HH}^\dag: \\
&T \mapsto (T^\dag)'=T^\dag - 1,\ \pi
\mapsto \pi,\ X \mapsto X^{-1},\
q^{1/4} \mapsto q^{-1/4},\notag
\end{align}
and its inverse, sometimes denoted by the same symbol
$\eta$. This definition
is a direct nil-variant of that from (\ref{etaqt}).

\begin{proposition}\label{ASTPOLYN}
There exists a unique $\C$\~linear automorphism denoted
by $\ast$ from
$\overline{\mathscr{X}}$ to $\overline{\mathscr{X}}^\dag$
and the one in the opposite direction
sending $1\mapsto 1, q\mapsto q^{-1}$ and compatible
with  $\eta$:
\begin{align}\label{etaction}
&\overline{\mathscr{X}}\ni H(f^\ast)=
(\eta(H)(f))^\ast \in \overline{\mathscr{X}}^\dag
\for H\in \overline{\HH},\\
&\overline{\mathscr{X}}^\dag\ni H(f^\ast)=
(\eta(H)(f))^\ast \in \overline{\mathscr{X}}
\for H\in \overline{\HH}^\dag.\notag
\end{align}
Upon the standard identification of the polynomial
representations with  $\C[q^{\pm 1/4}][X^{\pm 1}]$,
it becomes the conjugation, namely,
$(X^n)^\ast=X^{-n}$ for all $n$ and $q^{1/4} \mapsto q^{-1/4}$.
Moreover,
\begin{align}\label{astepol}
(\overline{E}_n)^\ast=X^{-1}\overline{E}_{1-n}^\dag \and
(\overline{E}_n^\dag)^\ast=X^{-1}\overline{E}_{1-n} \for n\in \Z.
\end{align}
\sq
\end{proposition}
\smallskip

\subsubsection{\sf Pairing with conjugation}
The remaining inner product from formula (\ref{innerp})
in the $\circ$\~normalization is:
\begin{equation} \label{barinner}
(\!(f,g)\!)_\circ = \lan fg^\ast\mu_\circ\ran.
\end{equation}
In the limit, it becomes a pairing between
$f\in \overline{\mathscr{X}}$ and
$g\in \overline{\mathscr{X}}^\dag\,$ (or in the opposite order):
\begin{equation}
(\!(f,g)\!)_\circ = \lan fg^\ast\overline{\mu}_\circ\ran, \ \,
(\!(1,1)\!)_\circ = 1, \ \,
(\!(g,f)\!)_\circ = \lan gf^\ast\overline{\mu}^\dag_\circ\ran
=(\!(f,g)\!)^\ast_\circ.
\end{equation}
Here $g^\ast\in \overline{\mathscr{X}}$, so
$\mu_\circ$ in the first formula becomes $\overline{\mu}_\circ$,
which is the limit of $\mu_\circ$ as  $t\to 0$. In the opposite
order, $\overline{\mu}^\dag_\circ=\lim_{t\to \infty}\mu_\circ$
must be used; see (\ref{normppolsbara1t}).

Check that
\begin{equation}
\lim_{t\to 0} (g^\ast) =
\left(\lim_{t\to\infty} g\right)^\ast,
\end{equation}
whenever the limits exist. Recall that $\overline{\mathscr{X}}$ and
$\overline{\mathscr{X}}^\dag$ are both equal to
$\C[q^{\pm 1/4}][X^{\pm 1}]$, but
the former is defined as a module over $\overline{\HH}$,
the latter over $\overline{\HH}^\dag$.

Taking the limit of the first formula in
(\ref{innerinv}), we obtain
\begin{equation}\label{innerinvwc}
(\!(H(f),g)\!)_\circ \ =\  (\!(f,H^\ast(g))\!)_\circ,
\end{equation}
where $H \in
\overline{\HH}$ and hence $H^\ast \in \overline{\HH}^\dag$.

For instance, let us consider $T$. Then in $\HH$:
\begin{equation} \lan T(f)g^\ast\mu_\circ\ran
= (\!(T(f),g)\!)_\circ = (\!(f,T^\ast(g))\!)_\circ = \lan f
(T^{-1}(g))^\ast\mu_\circ\ran.
\end{equation}
This relation must be multiplied by $t^{1/2}$ followed
by $t \to 0$. On the left-hand side $t^{1/2}T \to
\overline{T}$, while on the far right-hand side we need to consider
the limit of $t^{1/2}(T^{-1}(g))^\ast = (t^{-1/2}T^{-1}(g))^\ast$.
Moving inside $\ast$, $t^{-1/2}T^{-1}$ tends to $T^\dag -1$
as $t \to \infty$. Thus
\begin{equation}
(\!(\overline{T}(f),g)\!)_\circ = (\!(f, (T^\dag-1)(g))\!)_\circ.
\end{equation}

We note that (\ref{fnorme}) readily results in
\begin{align}\label{innerepol}
&(\!(\overline{E}_m,\overline{E}_n^\dag)\!)_\circ=\de_{mn}
\prod_{0<j<|\tilde{n}|}(1-q^j),\\
&(\!(\overline{E}_m^\dag,\overline{E}_n)\!)_\circ=\de_{mn}
\prod_{0<j<|\tilde{n}|}(1-q^{-j}).\notag
\end{align}

\medskip
\subsection{\bf Tilde-subalgebras}
A surprising fact is that the construction of
non-symmetric Whittaker functions
naturally leads to a module over
$\overline{\HH}^{\,\vph}\,$, which differs significantly from
the bar-polynomial representation. We will call it
the {\em hat-polynomial\,} or {\em spinor}
representation. It naturally appears within
the technique of spinors to be discussed later.

The spinor
representation is not a standard induced $\overline{\HH}$\~module,
but can be interpreted as a sub-induced module.
The problem is that the subalgebra $\bar{\h}_X=<T,X^{\pm1}>$
has no one-dimensional representations and we need to
diminish $\overline{\HH}$ by switching
from $X^{\pm1}$ to other (non-invertible) generators.
It can be addressed as follows.

\subsubsection{\sf Alternative presentations}
Let us begin with the following remark.
Motivated the theory of the DAHA in the $q,t$\~case,
one can try to eliminate $\pi$
from the definition of nil-DAHA. It is doable, but less
useful than in the general case.

The claim is that $\overline{\HH}=<T,X^{\pm 1},Y>$. The remaining
elements $\pi, Y'$  are expressed as follows:
\begin{align}\label{piviaXY}
\pi=q^{1/2}X^{-1}YX-Y,\ Y'=T'\pi=T'(q^{1/2}X^{-1}YX-Y).
\end{align}
In terms of $T,X,Y$, the defining relations will
be
\begin{align}\label{relxyt}
&T(T+1)=0,\ \pi^2=(q^{1/2}X^{-1}YX-Y)^2=1,\\
&Y(X+X^{-1})=(q^{1/2}X^{-1}+q^{-1/2}X)Y,\notag\\
&X(Y+Y')=(q^{1/2}Y+q^{-1/2}Y')X.\notag
\end{align}

However, the algebra $\overline{\HH}$ does not
have the PBW-property in terms of $T,X,Y$ (which
is the key fact in the $q,t$\~case). Also, as
noted above, the
affine Hecke subalgebra $\bar{\h}_X=<T,X^{\pm1}>$
has no one-dimensional representation due to the
invertibility of $X$.
\smallskip

One may also consider $\overline{\HH}$
as the algebraic span
$\overline{\HH}=<T,\pi,\tilde{\pi}>$ over $\C[q^{\pm1/4}]$. The
defining relations in this presentation are as follows:
\begin{align}\label{dahapitilde}
&T(T+1)=0,\ \, \pi^2\,=\,1\,=\,\tilde{\pi}^2,\\
&\pi X\pi =q^{1/2}X^{-1} \for X\equal q^{1/4}\tilde{\pi}\pi,\notag\\
&\tilde{\pi} Y\tilde{\pi} =q^{-1/2}Y' \for Y\equal\pi T.\notag
\end{align}
The function ord$(H)$ becomes
$\max_i\{\hbox{ord}(N_i)\}$
for any {\em reduced} expression $H=\sum_i c_iN_i$, where
$N_i$ are products (monomials) of $T,\pi,\tilde{\pi}$ and
where ord$(N_i)$ is the number of $T$ in the word for $N_i$.
By reduced, we mean that the sum  $H=\sum_i c_iN_i$ must have
the least possible $\max_i\{\hbox{ord}(N_i) : c_i\neq 0\}$ among
all such expressions
for $H$ in terms of the products of $T,\pi,\tilde{\pi}$.

The embedding of $\overline{\HH}$ into the abstract algebra
defined in
(\ref{dahapitilde}) is straightforward; the
formulas for the images of $X,Y$ are provided there. It is obviously
an isomorphism. The element $\tilde{\pi}$ is actually from
(\ref{diamondef}):
\begin{align}\label{pitildemap}
&\tilde{\pi}\equal\tau_+(\pi)=q^{-1/4}X\pi =q^{1/4}\pi X^{-1}.
\end{align}
\smallskip

\subsubsection{\sf Transitional subalgebras}
The approach we need in order to address the spinor representation
is actually opposite to the previous remarks; we need to
eliminate $\pi,\tilde{\pi}, X$ from the list of generators.
Let us introduce two {\em proper} subalgebras of $\overline{\HH}$:
\begin{align}\label{tildeHH}
&\widetilde{\HH}_{\tilde{\pi},Y}\equal<T,\tilde{\pi},Y,Y'>,\
\widetilde{\HH}_{\pi,\tilde{X}}\equal<T,\pi,\tilde{X},\tilde{X}'>,
\end{align}
where
\begin{align}\label{tildex}
&\tilde{X}\equal \tilde{\pi} T'=q^{1/4}YX,\\
&\tilde{X}'\equal T\tilde{\pi}=q^{-1/4}X^{-1}Y'.
\notag
\end{align}
These elements are direct nil-counterparts of $\tilde{X}^{\pm1}$
used in (\ref{psidef}). The {\em anti-involution} $\psi$ defined there
acts in $\overline{\HH}$ as does
$\tau_+$; note that $\psi\,\tau_+\,\psi=\tau_+^{-1}$.

One has:
\begin{align}\label{tildehh}
&\psi(\widetilde{\HH}_{\tilde{\pi},Y})=
\widetilde{\HH}_{\pi,\tilde{X}},\where \psi: \\
&Y\mapsto \tilde{X}',Y'\mapsto \tilde{X},\,
\tilde{\pi}\mapsto \pi,\, T\mapsto T.
\end{align}
The following relations hold:
\begin{align}\label{ytilde}
&\tilde{X}T=0=T'\tilde{X}', T\tilde{X}=\tilde{X}'T'.
\end{align}
\smallskip

\comment{
Similar to the definitions $Y=\pi T,Y'=T'\pi$, one can
introduce the elements in $\widetilde{\HH}_{\tilde{\pi},Y}$:
\begin{align}\label{tildey}
&\tilde{Y}\equal\tau_+(Y)=\tilde{\pi} T=q^{-1/4}XY,\\
&\tilde{Y}'\equal\tau_+(Y')=T'\tilde{\pi}=q^{1/4}Y'X^{-1}.
\notag
\end{align}
They satisfy the relations:
\begin{align}\label{ytilde}
&\tilde{Y}T'=0=T\tilde{Y}',\  T'\tilde{Y}=\tilde{Y}'T
(=T' \tilde{\pi} T),\\
&\breve{\pi} \tilde{Y}\breve{\pi} =q^{-1/2}\tilde{Y}' \for
\breve{\pi}\equal X\pi X^{-1}.\notag
\end{align}
Note that $\breve{\pi}$ does not belong to
$\widetilde{\HH}_{\tilde{\pi},Y}$.
\smallskip
}

\subsubsection{\sf Defining relations}
We claim that $\widetilde{\HH}_{\tilde{\pi},Y}$, considered
as an abstract algebra, has the following defining relations:
\begin{align}\label{HHtildepi}
&T(T+1)=0,\ \tilde{\pi}^2=1,\ \tilde{\pi} Y\tilde{\pi} =q^{-1/2}Y',\\
&YY'=0=Y'Y,\ YT'=0=TY',\ T'Y=Y'T,\ T'\equal T+1\notag
\end{align}
Applying $\psi$, the defining relations of
$\widetilde{\HH}_{\pi,\tilde{X}}$ are then
\begin{align}\label{HHtildex}
&T(T+1)=0,\  \pi^2=1,\ \pi \tilde{X} \pi=q^{1/2}\tilde{X}',\\
&\tilde{X}\tilde{X}'=0=\tilde{X}'\tilde{X},\
\tilde{X}T=0=T'\tilde{X}',\ T\tilde{X}=\tilde{X}'T'.\notag
\end{align}

It is obvious that all these relations hold in $\overline{\HH}$.
What is less obvious is that the corresponding homomorphisms
\begin{align}\label{tildeinject}
\widetilde{\HH}_{\tilde{\pi},Y}\ \rightarrow \ \overline{\HH}
\ \leftarrow\ \widetilde{\HH}_{\pi,\tilde{X}}
\end{align}
are injective if these algebras are defined by the presentations
above, not as subalgebras.
This follows from Theorem \ref{SPINPBW}, ($i,ii$) below.
\smallskip

\rmk The standard way to verify the injectivity in
(\ref{tildeinject}) and to check similar facts is as follows.

(a) First of all, one must check that the abstract algebras
from (\ref{HHtildepi},\ref{HHtildex}) satisfy the claims of
Theorem \ref{SPINPBW} below, except for the uniqueness part there,
i.e., without verification of the linear independence
of the terms listed there. Then one defines the tilde-polynomial
representations of these algebras as proper induced modules
(see below).

(b) Assuming that the uniqueness holds, the formulas for the
generators acting in the tilde-representations can be calculated.
Then one verifies {\em directly} that these formulas really
give representations of the corresponding algebras, so these
representations can be defined explicitly without any reference
to the induction construction.

(c) Next, it is not difficult to check that the terms listed
in Theorem \ref{SPINPBW} really are linearly independent
{\em as operators}
acting in the corresponding tilde-polynomial representations
for generic $q$. This readily gives the uniqueness claims
from this theorem {\em for all} $q$ and, finally,
the injectivity of the maps in (\ref{tildeinject}).

{\sf The nil-case}.
Using the order in nil-Hecke algebras simplifies
the considerations versus the $q,t$\~case; it can be used
instead of the order of operators
needed in (c). This order readily results in the {\em existence}
of the PBW-type decompositions of the elements in $\overline{\HH}$.
The fact that the nil-algebras are limits of $\HH$ (where the PBW
Theorem holds) gives the {\em uniqueness}; the latter results in
the injectivity from (\ref{tildeinject}).
\medskip

\subsection{\bf PBW and filtrations}
The next step will be the PBW Theorems for the
tilde-algebras.

\begin{theorem} {\sf (tilde-PBW)}\label{SPINPBW}
(i) An arbitrary $H\in \widetilde{\HH}_{\tilde{\pi},Y}$
can be uniquely represented as a linear
combination of the terms
$$
(\tilde{X}')^l T Y^m,\ \tilde{X}^l \tilde{\pi}\,(Y')^m,\
\tilde{\pi}\, Y^m (m>0),\ \tilde{X}^l Y^m,\
 (\tilde{X}')^l(Y')^m (l+m>0),
$$
where $l\ge 0,m\ge 0$.

(ii)  An arbitrary $H\in \widetilde{\HH}_{\pi,\tilde{X}}$
can be uniquely represented as a linear
combination of the terms
$$
(Y')^l T' \tilde{X}^m,\ Y^l \pi (\tilde{X}')^m,\
\pi \tilde{X}^m (m>0),\
Y^l \tilde{X}^m,\ (Y')^l(\tilde{X}')^m (l+m>0),
$$
where $l\ge 0,m\ge 0$.
\end{theorem}

{\it Proof.}
Let us check ($i$). We can
move $Y$ and $Y'$ to the right (through
$T$ and $\tilde{\pi}$ modulo lower terms). Therefore, any
$H\in \widetilde{\HH}_{\tilde{\pi},Y}$ can be represented
as a linear combination of the terms
$\tilde{M}M$, $\tilde{M}\tilde{\pi}\,M$ and
 $\tilde{M}TM$ for the monomials in the form
$\tilde{M}=\tilde{X}^l$ or $\tilde{M}=(\tilde{X}')^l$ and
$M=Y^m$ or $M=(Y')^m$, where $l,m\ge 0$.
Further reductions are based on vanishing
properties of the products of the generators and induction with
respect to the degree defined by $l+m$ (i.e., we can disregard $T,\pi$ when
they appear in the expressions).

Since $\tilde{X}T=0$ and $TY'=0$, the terms
$(\tilde{X}')^l T Y^m$ are sufficient among
those with $T$ in the middle.
Next,
$$
\tilde{X}Y'=\tilde{\pi}T'T'\pi=\tilde{\pi}T'\pi=\tilde{\pi}Y', \
\tilde{X}'\tilde{\pi}=T\tilde{\pi}^2=T,
$$
which makes the terms $\tilde{X}^l(Y')^m$ $(m>0)$,
$(\tilde{X}')^l\tilde{\pi}\, Y^m$ and
$(\tilde{X}')^l\tilde{\pi}\,(Y')^m$ unnecessary modulo the
terms of lower degree. Now
\begin{align}\label{txpry}
&\tilde{X}'Y=T\tilde{\pi}\,Y=q^{-1/2}TY'\tilde{\pi}=0,\\
&Y'\tilde{X}=T'\pi \tilde{X}=q^{1/2}T'\tilde{X}'\pi=0. \notag
\end{align}
The first of these identities completes the existence part.
See (\ref{degreedef}) for more discussion concerning the degree.
The {\em uniqueness} formally follows from the {\em existence}
and the fact that $\overline{\HH}$ is a limit
of $\HH$, which satisfies the PBW Theorem.

Note that the {\em existence} claim can readily obtained from
the existence part of Proposition \ref{SOUBLEWEYL} below and
can be also deduced from the explicit
formulas for the generators acting in the tilde-polynomial
representations (see the comment above).

Claim ($ii$) is very much similar. Let us list the
key identities necessary in this case:
\begin{align}\label{yprtx}
&Y'\tilde{X}=Y\tilde{X}'= YT'\tilde{X}=Y'\pi \tilde{X}'=0,\\
&Y\tilde{X}'=-\pi\tilde{X}', \ Y\pi\tilde{X}=-\tilde{X},\
Y'\pi=T'.\notag
\end{align}
\sq
\smallskip

\subsubsection{\sf Using the order}
A more conceptual way for {\em normal ordering} of the
operators described in the theorem is based on the
following proposition.

\begin{proposition}\label{SOUBLEWEYL}
Let $\hbox{ord}(N)$ be the number of $Y,Y',T,T'$ in
the elements (words) $N$ listed in (i) or in (ii).
We set
$$
\f_m\equal \{\sum_i c_i N_i : \hbox{ord}(N_i)\le m\}.
$$
Then
$H_1 H_2\in \f_{m_1+m_2}$ for  $H_1\in \f_{m_1}$
and $H_2\in \f_{m_2}$. The
graded algebras {\rm gr}$\tilde{\HH}$=
$\oplus_{m=0}^\infty\,\f_{m+1}/\f_m$
of $\widetilde{\HH}_{\tilde{\pi},Y}$ and
$\widetilde{\HH}_{\pi,\tilde{X}}$ are correspondingly
\begin{align}\label{tildewtpi}
&\tilde{\w}_{\tilde{\pi},Y}\equal <Y,Y',\tilde{X},\tilde{X}',
\tilde{\pi}>, \\
&\tilde{\w}_{\pi,\tilde{X}}\equal <Y,Y',\tilde{X},\tilde{X}',
\pi>\, \label{tildepi}
\end{align}
subject to the corresponding Weyl-type relations
\begin{align*}
\pi^2=&1,\  \tilde{\pi}^2=1,\ \,
YY'=0=Y'Y,\ \tilde{X}\tilde{X}'=0=\tilde{X}'\tilde{X},\\
&Y\tilde{X}'=0=\tilde{X}Y',\ \
\tilde{\pi}\tilde{X}=\tilde{X}'\tilde{\pi},\
\pi Y= Y'\pi, \\
\pi \tilde{X}&=q^{1/2}\tilde{X}'\pi,\ \
\tilde{\pi}Y= q^{-1/2}Y'\tilde{\pi}, \ \
\tilde{X}Y=q^{1/2}Y\tilde{X}.
\end{align*}
In particular, it follows that  $N\not\in \f_{m-1}$ if
ord$(N)=m \,(m>0)$ for any element $N$ from the theorem,
i.e., these elements are exactly of order ord$(N)$ with
respect to the filtration $\{\f_i\}$.
\end{proposition}

{\it Proof.} The procedure for taking {\em gr} is as follows.
We rescale the elements $N$ from the theorem:
$N=h^{-ord(N)}N(h)$, for instance, $T=h^{-1}T(h)$.
Then we send $h\to 0$; for instance, $(h^{-1}T(h))^2=h^{-1}T(h)+1$
results in $T^2=0$ for $T=T(0)$ in gr$\tilde{\HH}$.
Thus $T$ and $T'$ coincide under {\em gr}, $\pi Y=T=
T'=Y'\pi$ and $\tilde{\pi}\tilde{X}=T'=T=\tilde{X}'\tilde{\pi}$.
Similarly, $\tilde{X}Y'=\tilde{\pi}T'T'\pi=0$ in
{\rm gr}$\tilde{\HH}$ and
\begin{align*}
&\tilde{X}Y=\tilde{X}\pi T=q^{1/2}\pi \tilde{X}'T'=
q^{1/2}\pi T \tilde{X}=q^{1/2}Y\tilde{X}.
\end{align*}

Conjugating by $\tilde{\pi}$ and $\pi$,
we obtain the remaining relations in {\rm gr}$\tilde{\HH}$:
\begin{align*}
\tilde{X}'Y=0=Y\tilde{X}',\ \tilde{X}'Y'=q^{1/2}Y'\tilde{X}'.
\end{align*}
We see that the operators $\tilde{X}, Y$ and
$\tilde{\pi}$ or $\pi$ can be normally ordered in $gr$.
The uniqueness, the fact that the relations from the proposition
are really defining, follows from the fact that $gr$ is a limit
of the corresponding $\tilde{\HH}$, where the PBW Theorem was
already checked. We actually repeat here the deduction of
the linear independence of the elements from ($i,ii$) from
the fact that nil-DAHA are limits of the general DAHA.\sq

Note that the $\f_m$ are finite-dimensional
vector spaces in contrast to the ord-filtration
for $\overline{\HH}$, where all such spaces are infinite
dimensional since ord$(X^{\pm 1})=0$.
\medskip

\setcounter{equation}{0}
\section{\sc Induced representations}
\subsection{\bf Tilde-polynomial modules}
Replacing $X^{\pm1}$ by $\tilde{X},\tilde{X}'$
significantly increases the list of modules
of polynomial type.

\subsubsection{\sf Tilde-induction}
Let us begin with $\tilde{\HH}_{\pi,\tilde{X}}$.
The tilde-polynomial representations for this algebra are
$$
\tilde{\mathscr{X}}_{\pi,\tilde{X}}^{\ep,\pm}\equal
\hbox{Ind}_{\tilde{\h}_{\tilde{X}}}^
{\tilde{\HH}_{\pi,\tilde{X}}}\, \C_{\tilde{X}}^{\ep,\de},
$$
where $\tilde{\h}_{\tilde{X}}\equal <T,\tilde{X},\tilde{X}'>$.
Here
$\C^{\ep,\de}_{\tilde{X}}$ for $\ep=0,1, \de=\pm$ is the
restriction of
the one-dimensional representation of
$\overline{\h}_{\tilde{X}}\equal <T,\tilde{X},\tilde{X}',\tilde{\pi}>$
defined as follows (cf. (\ref{odimr})):
\begin{align}\label{odimrx}
&T'(1)=\ep,\, T(1)=\ep-1,\ \tilde{\pi}\,(1)=\de,\\
&\tilde{X}(1)=\de\ep,\ \tilde{X}'(1)=\de(\ep-1).\notag
\end{align}
\smallskip
Here and in the sequel, $\de=\pm$ is understood as $\pm 1$
when applicable.

\begin{theorem}\label{TILDEPOL}
For $\ep=0$, the natural map
$$
\tilde{\mathscr{X}}_{\pi,\tilde{X}}^{0,\pm} \to
\overline{\mathscr{X}}_{\tilde{X}}^{0,\pm}\equal
\mathrm{Ind}_
{\overline{\h}_{\tilde{X}}}^
{\overline{\HH}}\, \C^{0,\pm}_{\tilde{X}}
$$
is an isomorphism. Thus $\overline{\HH}$ naturally acts
in the former module. The counterpart of this claim for the algebra
$\widetilde{\HH}_{\tilde{\pi},Y}$ (see \ref{odimry} below)
is for $\ep=1$:
$$
\tilde{\mathscr{X}}_{\tilde{\pi},Y}^{1,\pm} \to
\overline{\mathscr{X}}_{Y}^{1,\pm}\equal
\mathrm{Ind}_
{\overline{\h}_{Y}}^
{\overline{\HH}}\, \C^{1,\pm}_{\tilde{X}}.
$$
\end{theorem}

{\it Proof.} The first map is an embedding due to
the PBW theorem.
Thus we need only check that it is surjective.
It suffices to check that
the induced  $\overline{\HH}$\~module
$\overline{\mathscr{X}}_{\tilde{X}}^{0,\pm}$
remains irreducible upon the restriction to
$\tilde{\HH}_{\pi,\tilde{X}}$ for {\em generic} $q$.
Here we diagonalize $\tilde{X},\tilde{X}'$ in
$\overline{\mathscr{X}}_{\tilde{X}}^{0,\pm}$ and check that their
spectrum is {\em simple}. This is completely analogous to
the diagonalization of $Y,Y'$ in the standard polynomial
representation $\mathscr{X}$. Then we use that $\pi$ intertwines
$\tilde{X}$ and $\tilde{X}'$.

The adjustment of $\ep$ when inducing from
$Y,Y'$ versus $\tilde{X},\tilde{X}'$
really is necessary; see Proposition \ref{COLLAPSE} below.
\sq
\smallskip

\subsubsection{\sf Explicit identification}
Theorem \ref{TILDEPOL} can be
obtained directly via explicit formulas for the action of
the generators, i.e.,
without any reference to the irreducibility of the
induced polynomial-type representations for generic $q$. Later we
will need the explicit formulas anyway.

Let us consider the induced modules governed
by the PBW Theorem from ($ii$) where we take
$\ep=0, \de=\pm$. There are five types of monomials
$N$ listed in Theorem \ref{SPINPBW}, ($ii$).
Only the following  evaluations of such $N$  at
$1\in \tilde{\mathscr{X}}_{\pi,\tilde{X}}^{\ep,\de}$
do not vanish:
$$
Y^l \pi (\tilde{X}')^m(1),\
Y^l(1) (l>0),\ (Y')^l(\tilde{X}')^m(1).
$$
We use that $T'(1)=0, \tilde{X}(1)=\tilde{\pi}\,T'(1)=0$.
Since $Y^l(1)=Y^{l-1}\pi T(1)=-Y^{l-1}\pi$, the evaluations
$Y^l(1)$ with $l>0$ can be omitted (this case is
included in the previous one). Moreover, the relation
$T(1)=-1$ results in
\begin{align}\label{piofone}
&\pi(1)=-\pi T(1)=-Y(1),\
Y^l \pi (\tilde{X}')^m(1)=-(-\de)^m Y^{l+1}(1).
\end{align}

We see that the evaluation at $1$
naturally leads to the identification
\begin{align}\label{spindirect}
&\tilde{\mathscr{X}}_{\pi,\tilde{X}}^{0,\de}\cong
\{\sum_{l\ge 0} \, a_{l+1} Y^{l+1}(1)+ a_{-l}(Y')^l(1)\};
\end{align}
the coefficients are from $\C[q^{\pm 1/4}]$.
To prove the proposition in these two cases ($\de=\pm$),
we need to check that the latter space is invariant
under the application of $\tilde{\pi}$, which
is missing in $\tilde{\HH}_{\pi,\tilde{X}}$.
Recall
that $\tilde{\pi}(1)=\de=\pm 1$; see (\ref{odimrx}).

Using the relation $\tilde{\pi} Y\tilde{\pi}=q^{-1/2}Y'$
from (\ref{dahapitilde}),
\begin{align}\label{pitildepol}
\tilde{\pi}(\sum_{l\ge 0}& a_{l+1} Y^{l+1}+ a_{-l} (Y')^l)(1)\\
&=\de(\sum_{l\ge 0}\, a_{l+1}q^{-(l+1)/2}(Y')^{l+1}+
a_{-l}q^{l/2}Y^l)(1),\notag
\end{align}
which gives the required ``extra" formula.

\comment{
Let us also provide the formula for the action of $\pi$ in this
case. We have:
\begin{align}\label{piep0}
&\pi(1)=-Y,\, \pi(Y)=-1,\,
\pi(Y^2)=TY(1)=-Y-Y',\\
\pi(Y^m)&=-Y^{m-1}-(Y')^{m-1},\
\pi((Y')^m)=Y^{m+1}-Y^{m-1}-(Y')^{m-1}\notag
\end{align}
for $m\ge 1$.
\smallskip
}

For $\tilde{\mathscr{X}}_{\tilde{\pi},Y}^{1,\pm}$,
the demonstration of Theorem \ref{TILDEPOL} via
the explicit formulas for $\pi$, missing
in $\widetilde{\HH}_{\tilde{\pi},Y}$, is analogous.
Note that we take $\ep=1$ in this case; see
Proposition \ref{COLLAPSE} below.
\smallskip

\subsubsection{\sf Basic operators}
Continuing our explicit analysis
based on the identification from (\ref{spindirect}),
let us obtain the
formulas for the action of $T, \pi$ in this tilde-polynomial module
(which are guaranteed by the induction construction).

Generalizing relations $TY'=YT', Y'T=T'Y$ and
using that $YY'=Y'Y$,
we arrive at
$T(Y')^2=YT'Y'=Y(T+1)Y'=Y^2 T',\, \ldots\, ,$
\begin{align}\label{typrime}
T(Y')^m=Y^mT',\ (Y')^m T=T' Y^m \for m>0.
\end{align}
Thus,
\begin{align}\label{tsumpol}
T\,(\sum_{l\ge 0}& a_{l+1} Y^{l+1}+ a_{-l} (Y')^l)(1)\\
&=\sum_{l\ge 0}\, -a_{l+1}(Y^{l+1}+(Y')^{l+1})(1).\notag
\end{align}

Similarly, $\pi Y=T$ and $\pi(Y^{m+1})=\pi Y^{m+1}(1)=T(Y^{m})$ for
$m\ge 0$. Using that \ $\pi Y'=\pi T' \pi=Y\pi +1$,
\begin{align}
\pi (Y')^2=(Y&\pi+1)Y'=Y(Y\pi+1)+Y'=
Y^2\pi+Y+Y',\ \ \pi (Y')^3=\notag\\
=(Y^2\pi+Y&+Y')Y'=Y^2(Y\pi+1)+(Y')^2=
Y^3\pi+Y^2+(Y')^2, \ldots,\notag\\
&\ \pi (Y')^m\ =\ Y^m\pi+Y^{m-1}+(Y')^{m-1} \for m>1.\label{piyprime}
\end{align}
Due to $\pi(1)=-Y(1)=-Y$, we obtain finally that
\begin{align}\label{pisumpol}
\pi\,&\bigl(\sum_{l\ge 0} a_{l+1} Y^{l+1}+ a_{-l} (Y')^l\bigr)(1)\\
&=\sum_{l\ge 0}\, -a_{l+1}\,M_l(1)
+\sum_{l\ge 0}\, -a_{-l}\,(Y^{l+1}-M_{l-1})(1),\notag
\end{align}
where $M_m=Y^m+(Y')^m$ for $m>0$, $M_0=1, M_{-1}=0$.

This provides a complete and effective
description of the action of $\overline{\HH}$
in the tilde-polynomial representation
$\tilde{\mathscr{X}}_{\pi,\tilde{X}}^{0,\de}$
upon the identification from
(\ref{spindirect}). Indeed, $Y,Y'$ act by the multiplication
and the action of $X,X^{-1}$ can be calculated using
(\ref{pitildepol},\ref{piyprime}) and $X=q^{1/4}\tilde{\pi}\pi$.
One has:
\begin{align*}
&X(1) = -\de q^{-\frac{1}{4}}Y'(1)=-\de q^{-\frac{1}{4}}Y',\
X(Y)(1) = -\de q^\frac{1}{4}, \\
&X(Y')(1) = -\de q^\frac{1}{4}(q^{-1}(Y')^2 - 1)(1)
=\de q^\frac{1}{4}(q^{-1}(Y')^2 - 1).\notag
\end{align*}
and for $m \geq 2$:
\begin{align}\label{yactionx}
X(Y^m)(1) = -\de q^\frac{1}{4}\,
(\,&q^\frac{m-1}{2}Y^{m-1}+q^{-\frac{m-1}{2}}(Y')^{m-1}\,)(1),\\
X(Y')^m(1)= -\de q^\frac{1}{4}\,
(\,q^{-\frac{m+1}{2}}(Y')^{m+1}
-&q^\frac{m-1}{2}Y^{m-1}-q^{-\frac{m-1}{2}}(Y')^{m-1}\,)(1)\notag.
\end{align}
These formulas can be identified with those for the action
of the spinor-Dunkl operators in the spinor representation;
see (\ref{spinidentif}) below.

\smallskip
Let us briefly consider the subalgebra
$\tilde{\HH}_{\tilde{\pi},Y}$
following part ($i$) of Theorem \ref{SPINPBW}.
The tilde-polynomial representations for this algebra are
$$
\tilde{\mathscr{X}}_{\tilde{\pi},Y}^{1,\pm}\equal
\hbox{Ind}_{\tilde{\h}_{Y}}^{\tilde{\HH}_{\tilde{\pi},Y}}\,
\C^{1,\pm}_{Y},
$$
where $\tilde{\h}_{Y}\equal <T,Y,Y'>$.
Analogous to ($ii$),
$\C^{\ep,\de}_{Y}$ is the
restriction of the one-dimensional representations of
$\overline{\h}_{Y,\pi}\equal<T,Y,Y',\pi>$
given by
\begin{align}\label{odimry}
&T'(1)=\ep,\, T(1)=\ep-1,\ \pi(1)=\de,\\
&Y(1)=\de(\ep-1),\ Y'(1)=\de\ep,\ \ep=0,1,\ \de=\pm.\notag
\end{align}
to its subalgebra $\overline{\h}_{Y}\equal <T,Y,Y'>$.
\smallskip

Now $\ep=1$, which gives that $T(1)=0=Y(1)$,
$\pi(1)=\de=Y(1)$ and $\tilde{\pi}(1)=\tilde{X}(1)$
due to $T'(1)=1$.
Since we induce from the same
$\overline{\h}_{Y,\pi}$ as for the polynomial representation,
the resulting tilde-representation
{\em canonically} coincides with
$\overline{\mathscr{X}}=\C[q^{\pm 1/4}][X^{\pm 1}]$.

\begin{proposition}\label{PITILDEY}
The polynomial representation $\overline{\mathscr{X}}$ can be
naturally identified with
\begin{align}\label{spindirecty}
&\tilde{\mathscr{X}}_{\tilde{\pi},Y}^{1,+}\cong
\{\sum_{l\ge 0} \, a_{l+1} \tilde{X}^{l+1}(1)+
a_{-l}(\tilde{X}')^l(1)\}
\end{align}
for the action of the powers of $\tilde{X},\tilde{X}'$ on $1$ in the
tilde-polynomial representation.
Equivalently,  $\{\tilde{X}^{l+1}(1), (\tilde{X}')^l(1), l\ge 0\}$
is a basis of $\overline{\mathscr{X}}$, where $H(1)$ is understood
here as $H$ applied to $1$
in $\overline{\mathscr{X}}$.\sq
\end{proposition}
The final form of this proposition is presented in
Theorem \ref{EPOLTILDE} below.
\medskip

\subsection{\bf The core subalgebra}
It is now quite natural to consider the
intersection subalgebra
\begin{align}\label{tildeint}
\tilde{\HH}_{Y,\tilde{X}}\equal
\widetilde{\HH}_{\tilde{\pi},Y}\cap\widetilde{\HH}_{\pi,\tilde{X}}.
\end{align}
By construction, this subalgebra is preserved by the anti-involution
$\psi$ from  (\ref{tildehh}) sending
\begin{align*}
Y\mapsto \tilde{X}',Y'\mapsto \tilde{X},\,
\tilde{\pi}\mapsto \pi,\, T\mapsto T.
\end{align*}

\begin{theorem}\label{TILDEINTER}
(i) {\sf (Defining relations)} As an abstract algebra,
$\tilde{\HH}_{Y,\tilde{X}}$ is generated by
$\tilde{X},\tilde{X}',Y,Y',T$ subject to
the relations
\begin{align}\label{tildeintrel}
&TT'=0=T'T=\tilde{X}\tilde{X}'=\tilde{X}'\tilde{X}=YY'=Y'Y,\notag\\
&\tilde{X}T=0=T'\tilde{X}'=
YT'=TY'=\tilde{X}'Y=Y'\tilde{X},\notag\\
&T\tilde{X}=\tilde{X}'T',\ \  T'Y=Y'T,\ \
\tilde{X}Y'=-q^{1/2}Y\tilde{X}',\notag\\
&\tilde{X}Y=q^{1/2}Y(\tilde{X}+\tilde{X}'),\
\tilde{X}'Y'=q^{1/2}(Y'+Y)\tilde{X}'.
\end{align}

(ii) {\sf (PBW Theorem)}
An arbitrary element $H\in \tilde{\HH}_{Y,\tilde{X}}$
can be uniquely expressed as a linear combination of
the elements
\begin{align}\label{tildeintpbw}
&Y^l\tilde{X}^m,\ \ (Y')^l(\tilde{X}')^m (l+m>0), \notag\\
&Y^l(\tilde{X}')^m (lm>0),\ \, (Y')^l (\tilde{X}')^m T.
\end{align}
Here the ordering $\{Y,Y'\},\{\tilde{X},\tilde{X}'\},T$
can be changed to its arbitrary permutation (similar to the
$q,t$\~case).

(iii) {\sf (Tilde-Weyl algebra)}
For the order from Proposition \ref{SOUBLEWEYL} restricted
to $\tilde{\HH}_{Y,\tilde{X}}$ the corresponding
{\rm gr}$\tilde{\HH}_{Y,\tilde{X}}$ is algebraically
generated by $\C$ and
the images of following elements in $\tilde{\HH}_{Y,\tilde{X}}$
of order one
\begin{align}\label{tildord}
\tilde{X},\tilde{X}',\ Y,Y',\ T,\ P\equal
\tilde{X}Y'=-q^{1/2}Y\tilde{X}'
\end{align}
subject to the nil-nil-Weyl relations:
\begin{align}\label{tildeintweyl}
\tilde{X}\tilde{X}'&=YY'=0=\tilde{X}'\tilde{X}=Y'Y,\
T\tilde{X}=\tilde{X}'T,\ TY=Y'T,\notag\\
&T^2=0=\tilde{X}T=T\tilde{X}'=
YT=TY'=\tilde{X}'Y=Y'\tilde{X},\notag\\
&P^2=0=P\tilde{X}=\tilde{X}'P=
PY=Y'P=\tilde{X}Y'=Y\tilde{X}',\notag\\
TP=&\tilde{X}'Y',\ \ PT=\tilde{X}Y,\ \ PY'=q^{1/2}YP,\ \
\tilde{X}P=q^{1/2}P\tilde{X}',\notag\\
&\ \ \ \ \ \ \ \ \tilde{X}Y=q^{1/2}Y\tilde{X},\ \
\tilde{X}'Y'=q^{1/2}Y'\tilde{X}'.
\end{align}

(iv) {\sf (Induced modules)} Setting
\begin{align}\label{interindy}
&\tilde{\mathscr{X}}_{Y}^{\ep,\pm}\equal
\mathrm{Ind}_{\tilde{\h}_{Y}}^{\tilde{\HH}_{Y,\tilde{X}}}\,
\C^{\ep,\pm}_{Y},\\
\label{interindx}
&\tilde{\mathscr{X}}_{\tilde{X}}^{\ep,\pm}\equal
\mathrm{Ind}_{\tilde{\h}_{\tilde{X}}}^{\tilde{\HH}_{Y,\tilde{X}}}\,
\C^{\ep,\pm}_{\tilde{X}}
\end{align}
for $\C^{\ep,\pm}$ defined in (\ref{odimr},\ref{odimrx}),
the natural maps
\begin{align*}
&\tilde{\mathscr{X}}_{Y}^{1,\de}\, \to\,
\tilde{\mathscr{X}}_{\tilde{\pi},Y}^{1,\de},\ \
\tilde{\mathscr{X}}_{\tilde{X}}^{0,\de}\, \to\,
\tilde{\mathscr{X}}_{\pi,\tilde{X}}^{0,\de}
\end{align*}
are isomorphisms for $\de=\pm.$
\end{theorem}
{\it Proof.}
The vanishing conditions for the monomials from part ($i$)
are due to $TT'=0=T'T$ directly or (less directly) as
in (\ref{txpry}) and (\ref{yprtx}).
Next, the relations
\begin{align}\label{txypr}
&Y\tilde{X}'=\pi T^2\tilde{\pi}=-\pi T\tilde{\pi}\\
=&-Y\tilde{\pi}=-q^{-1/2}\tilde{\pi}Y'=-\pi \tilde{X}'=
-q^{-1/2}\tilde{X}\pi,\notag
\end{align}
combined with the analogous ones
\begin{align*}
&\tilde{X}Y'=\tilde{\pi}T'\pi=\tilde{X}\pi=\tilde{\pi}Y',
\end{align*}
readily result in
$$
Y\tilde{X}'=-q^{-1/2}\tilde{X}\pi=-q^{-1/2}\tilde{X}Y'.
$$

Now,
\begin{align*}
&\tilde{X}Y=\tilde{\pi}T'Y=
\tilde{\pi}Y'T=q^{1/2}Y\tilde{\pi}T=
q^{1/2}Y\tilde{\pi}(T'-1)\\
&=q^{1/2}Y\tilde{X}-q^{1/2}Y\tilde{\pi}=
q^{1/2}Y(\tilde{X}+\tilde{X}'),
\end{align*}
where we use (\ref{txypr}) and \ $\tilde{\pi}Y'\tilde{\pi}=q^{1/2}Y$
\ from (\ref{HHtildepi}).
 Similarly,
\begin{align*}
&\tilde{X}'Y'=T\tilde{\pi}Y'=
q^{1/2}T Y\tilde{\pi}=
q^{1/2}(T'-1)Y\tilde{\pi}\\
&=q^{1/2}Y'T\tilde{\pi}-q^{1/2}Y\tilde{\pi}=
q^{1/2}Y'\tilde{X}'-q^{1/2}Y\tilde{\pi}=
q^{1/2}(Y'+Y)\tilde{X}'.
\end{align*}

We see that the relations from ($i$) hold in $\overline{\HH}$.
These relations are sufficient to ensure that the elements listed in
($ii$) linearly generate $\tilde{\HH}_{Y,\tilde{X}}$,
which is sufficient to establish their
linear independence.

The passage to $gr$ is straightforward; express $P$ using
relations (\ref{txypr}) and use (\ref{tildeintrel}); we emphasize
that $P$
is {\em not} $\tilde{X}Y'$
in \,{\rm gr\,}$\tilde{\HH}_{Y,\tilde{X}}$\,, where $\tilde{X}Y'=0$.
Since $P$ is $\psi$\~invariant,
the following are sufficient to check the relations
where $P$ appears:
\begin{align*}
&TP=T\tilde{X}Y'=\tilde{X}'(T')^2Y'=\tilde{X}'T'Y'=\tilde{X}'Y' \and
\\
&(Y\tilde{X}')Y'=q^{1/2}Y(Y'\tilde{X}'+Y\tilde{X})=
q^{1/2}Y^2\tilde{X})=q^{1/2}Y(Y\tilde{X}').
\end{align*}
Both identities hold in  $\tilde{\HH}_{\tilde{\pi},Y}$, i.e.,
before taking $gr$.

Finally, the claims from $(iv)$ follow
from $(i,ii)$ and the previous considerations for
$\tilde{\HH}_{\tilde{\pi},Y}$ and
$\tilde{\HH}_{\pi,\tilde{X}}$.
\sq
\smallskip

The last part of the theorem is
the most logically transparent way to introduce the
tilde-polynomial representations. The {\em intermediate}
subalgebras $\tilde{\HH}_{\tilde{\pi},Y}$ and
$\tilde{\HH}_{\pi,\tilde{X}}$ are transitional, as they are
convenient for calculating the formulas for the action of
$\pi,\tilde{\pi}$ and $X^{\pm 1}$.

We note that one can introduce a different order $\ \tilde{ord}\ $
in $\,\tilde{\HH}_{Y,\tilde{X}}\,$, not the
order $\, ord\, $ from $\tilde{\HH}$ used in $(iii)$,
by making $\tilde{ord}(P)=2$. Namely, the elements of tilde-order
no greater than $1$
will be $1,T,Y,Y',\tilde{X},\tilde{X}'$ and their linear
combinations. They will algebraically generate
$\tilde{gr}\tilde{\HH}_{Y,\tilde{X}}$ subject to
the relations from $(i)$ where we set $T^2=0$ and replace $T'$
everywhere by $T$.
We use that these relations are homogeneous for $\tilde{ord}$
upon the formal identification of $T$ and $T'$; for instance,
the relations
\begin{align*}
\tilde{X}Y\ =\ q^{1/2}Y(\tilde{X}+\tilde{X}'),\ \ \,
\tilde{X}'Y'\ =\ q^{1/2}(Y'+Y)\tilde{X}'
\end{align*}
remain unchanged after
taking $\,\tilde{gr}\,$.
Such graded algebra does not seem (immediately) connected with
the ``natural" nil-nil-Weyl algebra from $(iii)$ (which can be
generalized to higher ranks).
\smallskip

In contrast to the intermediate subalgebras and the
main nil-DAHA $\overline{\HH}$,
the algebra $\tilde{\HH}_{Y,\tilde{X}}$ is {\em bi-graded}.
We set
\begin{align}\label{degreedef}
&\hbox{deg}_{\tilde{x},y}(T)=0,\  \hbox{deg}_{\tilde{x},y}(AB)=
\hbox{deg}_{\tilde{x},y}(A)+\hbox{deg}_{\tilde{x},y}(B),\\
\hbox{deg}_y&(Y^m)=m=\hbox{deg}_y((Y')^m),\
\hbox{deg}_{\tilde{x}}(\tilde{X}^m)=m=
\hbox{deg}_{\tilde{x}}((\tilde{X}')^m),\notag\\
\hbox{deg}_y&(\tilde{X}^m)=\,0\,=\hbox{deg}_y((\tilde{X}')^m),\
\hbox{deg}_{\tilde{x}}(Y^m)=\,0\,=\hbox{deg}_{\tilde{x}}((Y')^m),
\notag
\end{align}
where $A,B\in \tilde{\HH}_{Y,\tilde{X}}$; $q$ is of degree zero.

Note that the subalgebras $\tilde{\HH}_{\tilde{\pi},Y}$ and
$\tilde{\HH}_{\pi,\tilde{X}}$ possess
$\hbox{deg}_y$ and $\hbox{deg}_{\tilde{x}}$, respectively (but
not both). This bi-grading plays a significantly more important
role in the paper than the order in $(iii)$.
\smallskip

\subsubsection{\sf Tilde-intertwiners}
The modules from
(\ref{interindy}) are graded and even bi-graded upon the following
modification of their definitions.
Let $\tilde{\mathscr{X}}_{\tilde{X}}^{\ep,u}$
be induced from  the one-dimensional
$\tilde{\h}_{\tilde{X}}$\~modules $\C^{\ep,u}_{\tilde{X}}$
introduced as follows:
\begin{align}\label{onetxu}
&\tilde{X}(1)=0=T'(1), \tilde{X}'(1)=-u \for \ep=0,\\
&\tilde{X}'(1)=0=T(1), \ \tilde{X}(1)=u \for \ep=1.\notag
\end{align}
The degree $\hbox{deg}_y$ is present here.
If $u$ is considered an independent variable such
that $\hbox{deg}_{\tilde{x}}(u)=1$, then the latter
degree can be extended to
$\tilde{\mathscr{X}}_{\tilde{X}}^{\ep,u}.$
Such bi-grading will not be used in the paper;
we will mainly treat $u$ as a constant.
However, $\hbox{deg}_y$ will be needed.

Analogously, the
$\tilde{\h}_{Y}$\~modules $\C^{\ep,u}_{Y}$ are
\begin{align}\label{oneyu}
&Y'(1)=0=T'(1), Y(1)=-u \for \ep=0,\\
&Y(1)=0=T(1), \ \, Y'(1)\, =\, u \, \for \ep=1.\notag
\end{align}
The corresponding induced modules are
denoted by $\tilde{\mathscr{X}}_{Y}^{\ep,u}$;
the grading here is $\hbox{deg}_{\tilde{x}}.$

These starting one-dimensional modules can be naturally extended to
one-dimensional $\tilde{\HH}_{Y,\tilde{X}}$\~modules,
which will be denoted by $\tilde{\C}^{\ep,u}_{\tilde{X}}$ and
$\tilde{\C}^{\ep,u}_{Y}$. The action of the remaining
operators is zero, namely, $Y(1)=0=Y'(1)$ for (\ref{onetxu}) and
$\tilde{X}(1)=0=\tilde{X}'(1)$ for (\ref{oneyu}).

The following proposition explains why only one choice
of $\ep$ is possible for each type of the induction;
see Theorem \ref{TILDEINTER},($iv$).

\begin{proposition}\label{COLLAPSE}
(i) The modules
\begin{align*}
&\tilde{\mathscr{X}}_{Y}^{0,u} \and
 \tilde{\mathscr{X}}_{\tilde{X}}^{1,u}
\end{align*}
are one-dimensional. They are equal to
$\tilde{\C}^{0,u}_{\tilde{X}}$ and
$\tilde{\C}^{1,u}_{Y}$, respectively.

(ii) Assuming that $q$ is generic (i.e., not a root of unity),
the only nonzero irreducible finite-dimensional
$\tilde{\HH}_{Y,\tilde{X}}$\~modules
are:
$$\tilde{\C}^{1,u}_{Y},\ \tilde{\C}^{0,u}_{Y} \hbox{\ \ or\ \ }
\tilde{\C}^{1,u}_{\tilde{X}},\ \tilde{\C}^{0,u}_{\tilde{X}}
\for u\neq 0.
$$
\end{proposition}
{\it Proof}. Let us start with
$\tilde{\mathscr{X}}_{\tilde{X}}^{1,u}$. Since $\ep=1$,
$$
T'(1)=1,\ T(1)=0=\tilde{X}'(1),\ \tilde{X}(1)=u.
$$
Using that $Y'\tilde{X}=0$, one of the defining relations
from part ($i$) of the theorem, we obtain that $Y'(1)=0$. Also,
$YT'=0$, which gives that  $Y(1)=YT'(1)=0$. Thus the induced
module is one-dimensional in this case, where only
$\tilde{X}$ and $T'$ are non-zero among
$Y,Y',\tilde{X},\tilde{X}',T,T'.$

Similarly, if $\ep=0$ and we induce from $Y,Y'$, then:
$$
T'(1)=0,\, T(1)=-1,\ Y'(1)=0,\ Y(1)=-u,
$$
which results in $\tilde{X}(1)=0$ due to $\tilde{X}T=0$
and in $\tilde{X}'(1)=0$ due to $\tilde{X}'Y=0$.
Part ($i$) is verified.

\begin{lemma}\label{TILDEINTW}
(i) Let us assume that $Y(e)=\La e$ for $\La\neq 0$.
Then $Y'(e)=0=\tilde{X}'(e)$. If
$\tilde{e}\equal\tilde{X}(e)\neq 0$, then
$$
Y(\tilde{e})=q^{1/2}\La\tilde{e},\
Y'(\tilde{e})=0=\tilde{X}'(\tilde{e}).
$$
Let $e'\equal T'(e)\neq 0$. Then
$
Y'(e')=\La e', \
Y(e')=0=T(e').
$
If $T'(e)= 0$ then $\tilde{X}(e)=0$ and
the module $\tilde{\HH}_{Y,\tilde{X}}(e)$ equals
$\tilde{\C}^{0,u}_{Y}$ for $u=-\La'$.

(ii) Now we assume that $Y'(e')=\La' e'$ for $\La'\neq 0$,
which automatically results in $Y(e')=0=T(e')$,
and set $\tilde{e}'\equal\tilde{X}(e')$.
Then
$$
Y(\tilde{e}')=q^{-1/2}\La' \tilde{e}',\
Y'(\tilde{e}')=0=\tilde{X}'(\tilde{e}').
$$
If here $\tilde{e}'=0$, then $\tilde{X}'(e')=0$ and
the module $\tilde{\HH}_{Y,\tilde{X}}(e')$ equals
$\tilde{\C}^{1,u}_{Y}$ for $u=q^{-1/2}\La'$.

(iii) The statements from (i) hold correspondingly for
$\tilde{e}\equal Y(e)$ and $e'\equal T(e)$ assuming that
$\tilde{X}(e)=\La e$ for $\La\neq 0$\,;\,
the latter implies that $Y'(e)=0=\tilde{X}'(e)$. In this case,
$$
\tilde{X}(\tilde{e})=q^{1/2}\La \tilde{e},\
\tilde{X}'(e')=\La e';\  T(e)=0 \Rightarrow
Y(e)=0.
$$
Following (ii), let
$\tilde{X}'(e')=\La' e'$ for $\La'\neq 0$ and
$\tilde{e}'\equal Y(e')$. Then
$\tilde{X}(\tilde{e}')=q^{-1/2}\La'\tilde{e}'.$
\end{lemma}
{\it Proof.}
First of all, $Y'Y(e)=0=\La Y'(e)$ and $\tilde{X}'(e)=0$ due
to the identity $\tilde{X}'Y=0$. For
$\tilde{e}\equal\tilde{X}(e),$
$$
Y(\tilde{e})=Y\tilde{X}(e)=(q^{1/2}\tilde{X}Y-Y\tilde{X}')(e)=
q^{1/2}\tilde{X}Y(e)=q^{1/2}\La\tilde{e}.
$$
Then automatically, $Y'(\tilde{e})=0=\tilde{X}'(\tilde{e}).$

Setting $e'\equal T'(e),$ we arrive at:
$$
Y'(e')=Y'T'(e)=Y'(T+1)(e)=T'Y(e)=\La e',\
Y(e')=0=T(e').
$$
Now, if $T'(e)= 0$ then $0=\tilde{X}T(e)=-\tilde{X}(e)$, since
$\tilde{X}T=0$ and $T(e)=-e$.

($ii$)
One has $YY'(e')=0=\La' Y(e)$ and  $TY'(e')=\La' T(e')=YT'(e')=0$.
Using (\ref{tildeintrel}) from Theorem \ref{TILDEINTER},
\begin{align}\label{ytxprime}
&Y\tilde{X}'(e')=-q^{-1/2}\tilde{X}Y'(e')=-q^{-1/2}\La'\tilde{X}(e')=
-q^{-1/2}\La'\tilde{e}'\\
&Y\tilde{X}'(e')=q^{1/2}\tilde{X}Y(e')-Y\tilde{X}(e')=
0-Y\tilde{X}(e')=-Y(\tilde{e}').\notag
\end{align}
If $\tilde{e}'=0$ then $\tilde{X}(e')=0$, which results
in
$$
\tilde{X}'T'(e')=0=\tilde{X}'(e')\hbox{\ due to \ }
T'(e')=e'.
$$

Note that $\tilde{X}=\tilde{\pi}T'$ is proportional
to the product $\Pi(T+1)$ of the intertwining
operators used in (\ref{PiTplus1}). Thus, this
claim and those from ($i$) formalize the technique of
intertwiners in the tilde-setting.
\smallskip

The counterparts of ($i$) and ($ii$) for $\tilde{X}$ instead
of $Y$ are straightforward. Let us demonstrate that
$\tilde{X}(\tilde{e}')=q^{-1/2}\La'\tilde{e}'$ if
for $\tilde{e}'=Y(e')$ if $\tilde{X}'(e')=\La' e'$.
We need the following straightforward corollaries of the
defining relation from (\ref{tildeintrel}):
\begin{align}\label{txyq}
\tilde{X}Y=q^{1/2}Y\tilde{X}-\tilde{X}Y', \
Y\tilde{X}= q^{-1/2}\tilde{X}(Y+Y').
\end{align}
Then, following (\ref{ytxprime}),
\begin{align*}
&\tilde{X}Y'(e')=-q^{1/2}Y\tilde{X}'(e')=-q^{1/2}\La'Y(e')=
-q^{1/2}\La'\tilde{e}',\\
&\tilde{X}Y'(e')=q^{1/2}Y\tilde{X}(e')-\tilde{X}Y(e')=
-\tilde{X}Y(e')=-\tilde{X}(\tilde{e}')
\end{align*}
and we can equate the outputs.
Similarly, if $\tilde{X}(e)=\La e$ for $\La\neq 0$,
then $\tilde{X}'(e)=0=Y'(e)$ and
(\ref{txyq}) gives that
$$
\tilde{X}(\tilde{e})=\tilde{X}Y(e)=
q^{1/2}Y\tilde{X}(e)-\tilde{X}Y'(e)=q^{1/2}\La\tilde{e}.
$$
\sq

The lemma readily results in the description
of the finite-dimensional irreducible representations
from Proposition \ref{COLLAPSE} ($ii$).
Applying the intertwiners, we will eventually
make $e'=0$ or $\tilde{e}'=0$, provided $q$ is
not a root of unity.
Thus these modules can be only as claimed.
\sq

\subsubsection{\sf Dual induced modules}
Since the induced modules from part ($i$)
of Proposition \ref{COLLAPSE} are trivial, another approach
must be used for establishing a link to
the $\overline{\HH}$\~modules
\begin{align*}
&\overline{\mathscr{X}}_{\tilde{X}}^{1,\de},\
\overline{\mathscr{X}}_{Y}^{0,\de}.
\end{align*}
These modules are isomorphic
to $\C[q^{1/4}][X^{\pm 1}]$ as linear spaces;
cf. Theorem \ref{TILDEPOL}. However
they cannot be defined directly within
the tilde-theory. One must proceed as
follows.

We set
$\hbox{Hom}_{\hbox{\tiny\,deg}}(\v)\equal
\lim_{m\to \infty}\, \hbox{Hom}(\v/\v_m\,,\,\C)$
for a vector space $\v$ with the filtration
$\v=\v_0\supset \v_1\supset \v_2\supset \cdots\ ,$ corresponding
to a given degree function: $\v_m=\{v\in \v : \hbox{deg}(v)\ge m\}$.
We set $\v_m^\vee\equal\hbox{Hom}(\v/\v_{m+1}\,,\,\C)$.
We will use the inner product $\lan f,g\ran'$ defined
in (\ref{innerpw}) with corresponding anti-involution
$\psi$.

\begin{proposition}\label{DUALTILDE}
Let us define the following
$\tilde{\HH}_{Y,\tilde{X}}$\~modules:
\begin{align}\label{interindych}
&(\tilde{\mathscr{X}}_{Y}^{\ep,u})^\vee \equal
\mathrm{Hom}_{\hbox{\tiny \rm\,deg\,}_{\tilde{x}}}\Bigl(
\mathrm{Ind}_{\tilde{\h}_{Y}}^{\tilde{\HH}_{Y,\tilde{X}}}\,
\C^{\ep,u}_{Y}\,,\ \C_q\Bigr),\\
\label{interindxch}
&(\tilde{\mathscr{X}}_{\tilde{X}}^{\ep,u})^\vee \equal
\mathrm{Hom}_{\hbox{\tiny \rm\,deg}_y}\Bigl(
\mathrm{Ind}_{\tilde{\h}_{\tilde{X}}}^{\tilde{\HH}_{Y,\tilde{X}}}\,
\C^{\ep,u}_{\tilde{X}}\,,\ \C_q\Bigr),
\end{align}
where the action of $\tilde{\HH}_{Y,\tilde{X}}$ in these spaces
is via the anti-involution $\psi$:
$H(f(P))\equal f(\psi(H)(P))$; $\C_q\equal\C[q^{\pm1/4}]$.

(i) Let $\v_m^\vee$
be $(\tilde{\mathscr{X}}_{Y}^{1,u})_m^\vee$
or $(\tilde{\mathscr{X}}_{\tilde{X}}^{0,u})_m^\vee$
with respect to $\hbox{deg}_{\tilde{x}}$ or
$\hbox{deg}_y$, $\v^\vee$ their inductive limit. Then
\begin{align}\label{dualydec}
&\v=\v_0\supset \v_1\supset \v_2 \supset \v_3\supset \, \cdots\, , \\
&\C_q\cong \v_1^\vee\subset \v_2^\vee\subset \v_3^\vee \
\cdots\
\subset \v^\vee\notag
\end{align}
are actually filtrations of
$\tilde{\HH}_{Y,\tilde{X}}$\~modules.

(ii) The pairing $\lan f,g\ran'$
induces the embeddings
\begin{align}\label{overhom}
&\overline{\mathscr{X}}_{\tilde{X}}^{1,\de}\,\hookrightarrow\,
\mathrm{Hom}\Bigl(
\overline{\mathscr{X}}_{Y}^{1,\de}\,,\, \C_q\Bigr),\ \
\overline{\mathscr{X}}_{Y}^{0,\de}\, \hookrightarrow\,
\mathrm{Hom}\Bigl(
\overline{\mathscr{X}}_{\tilde{X}}^{0,\de}\,,\,\C_q\Bigr),
\end{align}
where $\de=\pm$,
which result in the isomorphisms
of $\tilde{\HH}_{Y,\tilde{X}}$\~modules
\begin{align}\label{dualyx}
&(\tilde{\mathscr{X}}_{Y}^{1,\de})^\vee\, \cong \,
\overline{\mathscr{X}}_{\tilde{X}}^{1,\de},\ \
(\tilde{\mathscr{X}}_{\tilde{X}}^{0,\de})^\vee\, \cong \,
\overline{\mathscr{X}}_{Y}^{0,\de}.
\end{align}
\end{proposition}
{\it Proof.}
It suffices to check that the original descending
filtration
\begin{align}\label{directydec}
\v=\v_0\supset \v_1\supset \v_2 \supset \v_3 \ \cdots\
\end{align}
is that in terms of $\tilde{\HH}_{Y,\tilde{X}}$\~modules;
here $\v$ is equal to $\tilde{\mathscr{X}}_{Y}^{1,u}$
or $\tilde{\mathscr{X}}_{\tilde{X}}^{0,u}.$  This readily
follows from the properties of $deg_y$ and $deg_{\tilde{x}}$.
For (ii), we need to check that the embeddings in (\ref{overhom})
are compatible with the corresponding $deg$. It results from ($ii$);
see also Theorem \ref{DECOMIND} below.
\sq

\subsubsection{\sf Decomposition}
The deg-filtrations in
\begin{align}\label{yxinduced}
&\tilde{\mathscr{X}}_{Y}^{1,\de},\ \,
(\tilde{\mathscr{X}}_{Y}^{1,\de})^\vee,\ \ \
\tilde{\mathscr{X}}_{\tilde{X}}^{0,\de},\ \,
(\tilde{\mathscr{X}}_{\tilde{X}}^{0,\de})^\vee
\end{align}
give partial decomposition of the
$\overline{\HH}$\~modules
$\overline{\mathscr{X}}_{\tilde{X}}^{\ep,\de}$ and
$\overline{\mathscr{X}}_{Y}^{\ep,\de}$
upon the restriction to $\tilde{\HH}_{Y,\tilde{X}}$.

Recall that these filtrations $ \{\v_m (m\ge 0)\}$ and
$ \{\v_m^\vee (m\ge 1)\}$
are correspondingly either {\em descending co-finite-dimensional}
for the induced modules in (\ref{yxinduced})
or {\em ascending finite-dimensional} for their
dual counterparts.

For instance, in the check-case,
the modules $\v^\vee_1$ from the proposition are
$\tilde{\C}_{\tilde{X}}^{1,\de}$ or
$\tilde{\C}_{Y}^{0,\de}$, i.e., they are one-dimensional.
The quotients $\v^\vee_{m+1}/\v^\vee_{m}$ are two-dimensional
for $m\ge 1$. Since the irreducible
$\tilde{\HH}_{Y,\tilde{X}}$\~modules
are all one-dimensional by part ($ii$) of
Proposition \ref{COLLAPSE}, the constituents in this filtration
are not irreducible. The (complete) Jordan-Hoelder filtration
is actually a counterpart of the filtration in the
polynomial representation arising from the construction of
nonsymmetric Macdonald polynomials.

We use the identifications $A \mapsto A(1)$:
\begin{align}\label{dualyxid}
&\tilde{\mathscr{X}}_{Y}^{1,u}\, \cong\,
\C_q[\tilde{X},\tilde{X}'],\ \
\tilde{\mathscr{X}}_{\tilde{X}}^{0,u}\, \cong\,
\C_q[Y,Y'], \ \C_q=\C[q^{\pm 1/4}].
\end{align}
See (\ref{spindirect}) and (\ref{spindirecty}).

\begin{theorem}\label{DECOMIND}
For $u\neq 0$, let $\{\v_n : n\ge 0\}$ be the filtration of
$\tilde{\HH}_{Y,\tilde{X}}$\~submodules in
\begin{align*}
&\tilde{\mathscr{X}}_{Y}^{1,u} \hbox{\ \ \, or\ \ \ }
\tilde{\mathscr{X}}_{\tilde{X}}^{0,u}
\end{align*}
defined in Proposition \ref{DUALTILDE}, (ii).
For $n\ge 1$, we extend this filtration by adding the modules
$$
\v_{-n}=\C_q (\tilde{X}^n+(\tilde{X}')^n)+\v_{n+1} \hbox{\ \, or\ \ }
\v_{-n}=\C_q (Y^n+(Y')^n)+\v_{n+1}.
$$
Then $\{\v_n : n\in \Z\}$
is a filtration of $\tilde{\HH}_{Y,\tilde{X}}$\~modules.
The corresponding quotients of this filtration
are the one-dimensional modules
\begin{align}\label{filty1}
&\tilde{\C}_{Y}^{1,u},\  \tilde{\C}_{Y}^{0,q^{-1/2}u},\
\tilde{\C}_{Y}^{1,q^{-1/2}u},\ \tilde{\C}_{Y}^{0,q^{-1}u},\ \ldots\,
\for \tilde{\mathscr{X}}_{Y}^{1,u},\\
&\tilde{\C}_{\tilde{X}}^{0,u},\  \tilde{\C}_{\tilde{X}}^{1,q^{-1/2}u},\
\tilde{\C}_{\tilde{X}}^{0,q^{-1/2}u},\
\tilde{\C}_{\tilde{X}}^{1,q^{-1}u},\ \ldots\
\for \tilde{\mathscr{X}}_{\tilde{X}}^{0,u}.\label{filtx0}
\end{align}
\end{theorem}
{\it Proof.} Let us consider only the case of (\ref{filty1}).
We can apply the intertwiner $\tilde{X}$
from part ($ii$) of Lemma \ref{TILDEINTW}
to $e'=e_0=1$; indeed, $T'(1)=1$ and $Y'(1)=u$. Automatically,
$T(1)=0=Y(1)$. The element $e_1=\tilde{e}'=\tilde{X}(1)$ will satisfy
the conditions of part ($ii$) in the Lemma and we can now apply $T'$,
which will result in
$$
e_{-1}=T'\tilde{X}(1)=(T+1)\tilde{X}(1)=(\tilde{X}+\tilde{X}')(1).
$$
The next intertwiner in this chain will be $\tilde{X}$ again,
resulting in $e_2=\tilde{X}^2(1)$. Then
$e_{-2}=(\tilde{X}^2+(\tilde{X}')^2)(1)$ and so on:
\begin{align}\label{enminusn}
e_{n}=\tilde{X}^n(1),\ \ e_{-n}=(\tilde{X}^n+(\tilde{X}')^n)(1)
\for n\ge 1.
\end{align}

They all will be $Y,Y'$\~eigenvectors. Since the elements
$e_{-n}=(\tilde{X}^n+(\tilde{X}')^n)(1)$ for $n>0$ satisfy
$T(e_{-n})=0$, the module
$\tilde{\HH}_{Y,\tilde{X}}(e_{-n})$ divided by $\v_{n+1}$
is one-dimensional.

The consecutive quotients are exactly those  claimed in the theorem;
the description of the $Y,Y'$\~spectrum of $e_{\pm n}$
readily follows from the lemma.
\sq
\smallskip

Actually, our construction is a
counterpart of the {\em chain of intertwiners} in
the bar-polynomial representation $\overline{\mathscr{X}}$
from (\ref{PiTplus1}). Note that $u=1$ for the polynomial
representation: $Y'(1)=T'\pi(1)=1$. Also $\tilde{X}(1)
=\tilde{\pi} T'(1)=\Pi(1)$, i.e., the first intertwiner
we used is actually $\Pi$ from (\ref{PiTplus1}).
However, now the chain of intertwiners solves the problem
of the decomposition of the bar-polynomial representation
with respect to the action of $\tilde{\HH}_{Y,\tilde{X}}$,
i.e., contains different information.

Using the dualization, one obtains the ascending
filtrations in the remaining modules from
(\ref{yxinduced}).
\smallskip

Interestingly, using intertwiners in the nil-case can be done
within the subalgebra $\tilde{\HH}_{Y,\tilde{X}}$, i.e.,
the classical (invertible) intertwiners $\tilde{\pi}=\Pi$ and
$\pi$ can be avoided; the latter is a counterpart
of $\tilde{\pi}=\Pi$ for the modules induced from
$\tilde{X},\tilde{X}'$. It is directly
related to Proposition \ref{PITILDEY}, which can be
now stated in its final form. We come to the following
simple but fundamental theorem.

\begin{theorem}\label{EPOLTILDE}
Let us consider the elements $\{\,e_{n} : n\in \Z\,\}$ from
(\ref{enminusn}) in the polynomial representation
$\overline{\mathscr{X}}$. Then
\begin{align}\label{epoltildeformula}
\overline{E}_n &= q^{n^2/4}e_n,\ n \in \Z.
\end{align}
More explicitly,
\begin{align}\label{epoltildeformulaex}
\overline{E}_n &= q^{n^2/4}\tilde{X}^n(1),\ n \geq 0,\\
\overline{E}_{-n} &= q^{n^2/4}(\tilde{X}^n+(\tilde{X}')^n)(1),\
n > 0.
\end{align}
\end{theorem}
{\it Proof.} The $e_n$ were constructed using the intertwiners
$T'$ and $\tilde{X}$. Notice that $\tilde{X}(1) = q^{-1/4}\Pi(1)$
and that
\begin{align}
\tilde{X}T' &= q^{-1/4}\Pi(T')^2 = q^{-1/4}\Pi T'.
\end{align}
Hence we may just as well have used the intertwiners $T'$ and
$q^{-1/4}\Pi$, which are essentially those used to construct the
$\overline{E}_n$. Comparing with (\ref{PiTplus1}), one obtains
(\ref{epoltildeformula}).
\sq
\smallskip

\subsection{\bf Shapovalov forms}
Following \cite{C102}, let us interpret the
inner product $\lan f,g\ran'$ in the bar-polynomial
representation, which was defined in (\ref{innerpw}),
as a Shapovalov-type bilinear form.
\smallskip

\subsubsection{\sf Shapovalov trace}
The {\em Shapovalov trace} on the core
subalgebra $\tilde{\HH}_{Y,\tilde{X}}$
is given by
\begin{align}\label{shapow}
&
\lan (a(\tilde{X}')^{m'+1}+b\tilde{X}^{m})\,
(T')^\varsigma\, (c Y^{n+1}+d(Y')^{n'})\ran_\psi'\ \equal\ bd,
\end{align}
where $ \varsigma=0,1, \ m,m',n,n'\ge 0$.
The coefficients $a,b,c,d$ are
from $\C_q=\C[q^{\pm 1/4}]$. Recall that
$$
\tilde{X}=\psi(Y')=q^{1/4}YX,\  \tilde{X}'=\psi(Y)=q^{-1/4}X^{-1}Y'.
$$
We use here the PBW part of Theorem \ref{TILDEINTER}.

Generally, the Shapovalov trace $\lan\cdot\ran_\al$ requires a
subalgebra $\a$ in $\tilde{\HH}_{Y,\tilde{X}}$ or in
$\overline{\HH}$,
its character $\chi:\a\to \C$ and an anti-involution
$\al$. The defining relations are
$\lan\al(A)H\ran_\al=0=\lan HA\ran_\al$\, whenever \, $\chi(A)=0$
for any $H$.

Assuming that $\tilde{\HH}_{Y,\tilde{X}}$ is linearly generated
by the products from $\al(\a)\cdot\a$ (or $\overline{\HH}$ is
such) the dimension of the linear space of the Shapovalov traces
(called {\em DAHA-coinvariants} in \cite{C102})
is no greater than one.
It is obvious, as is the relation
$\lan H^\al \ran_\al=\lan H\ran_\al$ for
$H\in \tilde{\HH}_{Y,\tilde{X}}.$

Note that the definition from (\ref{shapow})
is adjusted to the concrete induced module, which is
the representation
$\tilde{\mathscr{X}}_{\tilde{\pi},Y}^{1,+}$ lifted to
the bar-polynomial representation of $\overline{\HH}$.
Namely, the algebra is $\a=<Y,T,\pi>$ and $\chi(H)=H(1)$ are
as follows:
$$
Y(1)=T(1)=0,\ \pi(1)=1 \and  Y'(1)=T'(1)=1.
$$
Generally, there are four possible choices for  $\chi$:
one can combine $T(1)=\ep=0,1$ with $\pi(1)=\de=\pm 1$.

The functional $\lan\cdot\ran_\psi'$
naturally maps through the tilde-module
$$
\tilde{\mathscr{X}}_{Y,\tilde{X}}^{1,+},
\hbox{\ which is naturally identified with \ }\overline{\mathscr{X}},
$$
so it can be extended to $\overline{\HH}$.
More explicitly, this extension
is uniquely determined by the relations
\begin{align}\label{shapfull}
&\lan\, A\,\pi^\varsigma\,(a Y^{n+1}+b (Y')^{n'})\,\ran_\psi'\ =\
b\lan\, A\, \ran_\psi'\\
&=\ \lan\, (a(\tilde{X}')^{n+1}+b\tilde{X}^{m'})\,
\tilde{\pi}^\varsigma\, A \,\ran_\psi'\notag
\end{align}
for $A\in\overline{\HH}$.

In particular, $\pi^\varsigma$ or $\tilde{\pi}^\varsigma$ can
be inserted instead of $(T')^\varsigma$ in (\ref{shapow});
recall that $\tilde{\pi}=\psi(\pi)$.
The formula $X=q^{1/4}\tilde{\pi}\pi$
readily results in $\lan X\ran_\psi'=q^{1/4}$.
Also $T'X^{-1}=XT$, which gives that $\lan X^{-1}\ran_\psi'=0$.
Moreover,
$$
T'X^{-n}=X^nT-X^{n-2} \for n\ge 2 \and
\lan X^{-n}\ran_\psi'=-\lan X^{n-2}\ran_\psi'.
$$
Combining this relation with
\begin{align*}
&\lan X^n\ran_\psi'=
\lan q^{n/4}(\tilde{\pi}\pi)^n\ran_\psi'
=\lan q^{n/4}(\tilde{\pi}\pi)^{n-1}\ran_\psi'=
\lan q^{n/4}q^{(n-1)/4}X^{1-n}\ran_\psi',
\end{align*}
we arrive at the formulas for the traces
of arbitrary $X$\~monomials:
\begin{align}\label{shapformulax}
&\lan X^n \ran_\psi'=-q^{n/2-1/4}\lan X^{n-3} \ran_\psi' \for
n\ge 1,\ \lan X^{-n}\ran_\psi'=\lan X^{n-2}\ran_\psi',\\
\where &\lan 1\ran_\psi'=1,\, \lan X \ran_\psi'=q^{1/4},\,
\lan X^2 \ran_\psi'=0=\lan X^{-1} \ran_\psi',\,
\lan X^{-2} \ran_\psi'=-1.
\notag
\end{align}
\smallskip

\subsubsection{\sf Shapovalov pairing}
The next step is the pairing
$\lan A, B\ran_\psi'\equal \lan A^\psi\, B\ran_\psi'$
for $A,B\in \overline{\HH}$. It is symmetric and induces the
anti-involution $\psi$ due to $\lan A^\psi \ran_\psi'=\lan A\ran_\psi'$.
Its $q,t$\~generalization is an important particular case of the
theory of Shapovalov-type forms from \cite{C102}.

\begin{theorem}\label{SHAPOL}
The pairing above naturally maps through
$\overline{\mathscr{X}}\times \overline{\mathscr{X}}$, so we come
to the following definition:
\begin{align}\label{shpoformpol}
\lan f,g\ran'_{\hbox{\tiny alg}}\equal \lan f(X),g(X)\ran_\psi'=
\lan g,f\ran'_{\hbox{\tiny alg}}.
\end{align}
The corresponding anti-involution is $\psi$ by construction.
For any $q$ (including the roots of unity), the following
holds:
$$
\lan f,g\ran'_{\hbox{\tiny alg}}=\lan f,g\ran' \for
f,g\in \overline{\mathscr{X}}.
$$
\end{theorem}
{\it Proof.} We use the irreducibility of
$\overline{\mathscr{X}}$ for generic $q$ and the uniqueness
of the bi-linear form in this representation associated with
$\psi$ (equivalently, the PBW claim from Theorem \ref{TILDEINTER}).
\sq

Note that the Shapovalov construction automatically holds for
arbitrary $q$. It is instructional to compare the coefficients
of the (analytic) expansion of $\overline{\mu}\tga'$ in terms
of $X^{-n}$ from (\ref{mugaone}) with the entirely algebraic
formulas for $\lan X^n \ran_\psi'$ from
(\ref{shapformulax}). The theorem states that they must
coincide; they really do.

\subsubsection{\sf The tilde-case}
We will use the identification of $\overline{\mathscr{X}}$ from
Proposition \ref{PITILDEY}:
\begin{align}\label{spindirectx}
&\tilde{\mathscr{X}}_{\tilde{\pi},Y}^{1,+}=
\tilde{\mathscr{X}}\equal
\{\sum_{l\ge 0} \, a_{l+1} \tilde{X}^{l+1}(1)+
a_{-l}(\tilde{X}')^l(1)\}.
\end{align}
Here the coefficients are from $\C[q^{\pm 1/4}]$. For the sake of
definiteness, only the case of $\ep=1$ and $\de=+$
will be considered here.

We now set
\begin{align}\label{innerpwt}
&\lan f,g\ran_\circ=\lan fg\tilde{\mu}_\circ \ran,\
\lan f,g\ran'=\lan fg \tilde{\mu\ga}' \ran
\end{align}
for $f,g\in \tilde{\mathscr{X}}$, where the
kernels of these pairings are those from
(\ref{mutildemuoneh}) and (\ref{mugaone}) upon the
substitution:
$X\mapsto \tilde{X}, X^{-1}\mapsto \tilde{X}'$. Recall that
$\tilde{X}\tilde{X}'=0=\tilde{X}'\tilde{X}$.
Explicitly,
\begin{align}\label{mutildex}
&\tilde{\mu}_\circ\ =\
\sum_{n=0}^\infty (-1)^nq^{\frac{n^2-n}{2}}
(\tilde{X}^{2n}-(\tilde{X}')^{2n-2}),\\
&\tilde{\mu\ga}'\ = \
\sum_{n=0}^\infty q^{n(n+2)/12}(\tilde{X}^{n+2}-(\tilde{X}')^{n}),
\where n\neq 2 \hbox{ mod } 3.\label{gatildex}
\end{align}

We claim that these inner products induce the
anti-involutions $\diamond$ and $\psi$, respectively, for the action
of $\overline{\HH}$
in $\tilde{\mathscr{X}}_{Y,\tilde{X}}^{1,+}$.
Recall that the anti-involution $\psi$ preserves the core
subalgebra $\tilde{\HH}_{Y,\tilde{X}}$ from
(\ref{tildeint}).

The definitions of the Shapovalov trace and
the pairing $\lan A,B\ran_\psi'$ in $\overline{\HH}$
remain unchanged. The difference is that now we
restrict this abstract pairing to ``functions" of
$\tilde{X}$.

\begin{theorem}\label{SHAPOLT}
The Shapovalov-type form
\begin{align}\label{shpoformpolt}
\lan f,g\ran'_{\hbox{\tiny alg}}\equal
\lan \psi(f)g\ran'_\psi=
\lan g,f\ran'_{\hbox{\tiny alg}}, \hbox{\, where\, }
f\in \tilde{\mathscr{X}}\ni g,
\end{align}
is associated to the anti-involution $\psi$. For any $q$,
$$
\lan f,g\ran'_{\hbox{\tiny alg}}=\lan f,g\ran' \for
f,g\in \tilde{\mathscr{X}}.
$$
\end{theorem}
\smallskip

\subsubsection{\sf Induction from tilde-X}
Let us briefly discuss the changes needed for
the modules induced from $\tilde{X}$:
\begin{align}\label{spindirectnew}
&\tilde{\mathscr{X}}_{\pi,\tilde{X}}^{0,-}=
\tilde{\mathscr{Y}}\equal
\{\sum_{l\ge 0} \, a_{l+1} Y^{l+1}(1)+ a_{-l}(Y')^l(1)\}.
\end{align}
See (\ref{spindirect}). We set
\begin{align*}
&\lan f(Y),g(Y)\ran_\circ=\lan f(Y)g(Y)\tilde{\mu}_\circ(Y) \ran,\
\lan f(Y),g(Y)\ran'=\lan fg \tilde{\mu\ga}'(Y) \ran
\end{align*}
for $f(Y),g(Y)\in \tilde{\mathscr{Y}}$, where the
kernels of these pairings are from
(\ref{mutildemuoneh}) and (\ref{mugaone}) upon the
substitution:
$X\mapsto Y', X^{-1}\mapsto Y$.
Explicitly,
\begin{align*}
&\tilde{\mu}_\circ(Y)\ =\
\sum_{n=0}^\infty (-1)^nq^{\frac{n^2-n}{2}}((Y')^{2n}-Y^{2n-2}),\\
&\tilde{\mu\ga}'(Y)\ = \
\sum_{n=0}^\infty q^{n(n+2)/12}((Y')^{n+2}-Y^{n}),\where
n\neq 2 \hbox{ mod } 3.
\end{align*}

We claim that these inner products induce the anti-involutions
$\diamond$ and $\psi$ for the action
of $\overline{\HH}$
in $\tilde{\mathscr{Y}}$ supplied with the
action of $\overline{\HH}$ via
$\tilde{\mathscr{X}}_{\pi,\tilde{X}}^{1,+}$.

Concerning the Shapovalov-type construction, we need
to reverse the order of operators in (\ref{shapow}).
The Shapovalov trace now reads:
\begin{align}\label{shapowx}
&
\lan\, (a(Y')^{m'}+bY^{m+1})\,
(T')^\varsigma\, (c\tilde{X}^{n}+ d(\tilde{X}')^{n'+1})\,
\ran_\psi^\backprime\ =\ ac,
\end{align}
where $\varsigma=0,1$, $m,m',n,n'\ge 0$. Its extension to the
whole $\overline{\HH}$ is determined from
\begin{align}\label{shapowext}
&\lan\, (a(Y')^{m'}+bY^{m+1})\,\pi^\varsigma\,A\,
\ran_\psi^\backprime\ =\
a\lan\, A \,\ran_\psi^\backprime\\
&=\ \lan\, A\, \tilde{\pi}^\varsigma\,
(a\tilde{X}^{n}+ b(\tilde{X}')^{n'+1})\,
\ran_\psi^\backprime.\notag
\end{align}
For instance, the relations from (\ref{shapowx}) are
satisfied with $\pi^\varsigma$ or $\tilde{\pi}^\varsigma$
instead of $(T')^\varsigma$ ($\varsigma=0,1$).

Theorem \ref{SHAPOLT} must be now stated in terms of
the restriction of the pairing $\lan A,B\ran_\psi^\backprime$ to
the ``functions" of $Y$:
$$
\lan f(Y),g(Y)\ran^\backprime_{\psi}=
\lan f(Y),g(Y)\ran' \for
f(Y),g(Y)\in \tilde{\mathscr{Y}}.
$$
\medskip

\setcounter{equation}{0}
\section{\sc Spinors and the
\texorpdfstring{$q$}{q}
-Toda theory}

In contrast to QMBP, the $q$\~Toda operators (\ref{qToda})
are nonsymmetric;
the corresponding Dunkl operators were not
expected to exist. The formulas from \cite{C102}
were really surprising. They provided an
exact counterpart of the $q,t$\~representation
$\l=Y+Y^{-1}$ (upon the restriction
to the symmetric functions) for $\l$ from (\ref{Lopertaor}),
but in the spinor sense. The introduction of the
spinor-Dunkl operators made it possible to use
DAHA methods at full potential, including the
theory of the $q$\~Whittaker functions.
This construction can be readily extended to arbitrary
root systems. The fact that the fundamental operators
are of first order dramatically simplifies
the theory.
We will begin with the definition of the spinors
following \cite{C102}.

\subsection{\bf The spinors}\label{sect:spinors}
Generally, $W$\~spinors are needed in the DAHA theory.
In the $A_1$\~case,
we will call them simply {\em spinors}. In this case,
they really are connected with spinors
from the theory of the Dirac operator (and with
super-algebras). Under the rational degeneration,
the Dunkl operator for $A_1$ becomes
the square root of the (radial part of the)
Laplace operator, i.e., the {\em Dirac operator}.
However, this direct relation (and using super-variables)
is a special feature of the root system $A_1$.

Let us mention that there are works in which DAHA is coupled with
the Clifford algebra. This approach provides a formula for the Dirac
operator in terms of Dunkl operators for any root systems;
however, it is a different approach.

For practical calculations with spinors, the language of
$\Z_2$\~graded algebras can be used in the $A_1$\~case.
However, we prefer
to do it in a way that does not rely on the special
symmetry of the $A_1$\~case and can be transferred
to $W$\~spinors for arbitrary root systems.
\smallskip

\subsubsection{\sf The definitions}
The {\em spinors\,} are simply pairs $\{f_1,f_2\}$ of
elements (functions) from a space $\f$ with an action
of  $s$; the addition or multiplication
(if applicable) of spinors is componentwise.
The space of spinors will be denoted by $\hat{\f}$.

The involution $s$ on spinors is defined as
$s\{f_1,f_2\}=\{f_2,f_1\}$, so it does not involve the
action of $s$ in $\f$. There is a ``natural"
embedding $\rho:\f\to \hat{\f}$ mapping
$f\mapsto f^\rho=\{f,s(f)\}$ and also the diagonal
embedding $\delta: \f\to \hat{\f}$ sending
$f\mapsto f^\de=\{f,f\}$. Accordingly,
for an arbitrary operator $A$ acting in $\f$,
$A^\rho=\{A,s(A)\}, A^\de=\{A,A\}$. The images
$f^\rho$ of $f\in\f$ are called {\em functions} (in contrast to
{\em spinors}) or {\em principal spinors} (like for adeles).

For instance, for $\f=\mathscr{X}$,
\begin{align*}
&X^\rho:\{f_1,f_2\}\mapsto\{Xf_1,X^{-1}f_2\},\
&&\Gamma^\rho:\{f_1,f_2\}\mapsto\{\Gamma(f_1),
\Gamma^{-1}(f_2)\},\\
&X^\de:\{f_1,f_2\}\mapsto\{Xf_1,Xf_2\},\
&&\Gamma^\de:\{f_1,f_2\}\mapsto\{\Gamma(f_1),\Gamma(f_2)\},
\end{align*}
where, recall, $\Gamma(X)=q^{1/2}X$. We simply put
$$
X^\rho=\{X,X^{-1}\},\, \Ga^\rho=\{\Ga,\Ga^{-1}\},
X^\de=\{X,X\},\, \Ga^\de=\{\Ga,\Ga\}.
$$
Obviously, $s^\rho=s=s^\de$.

If a function $f\in \f$ or an operator $A$ acting in $\f$
has no super-index $\de$, then
they will be treated as $f^\rho, A^\rho$. I.e.,
by default, functions
and operators are embedded
into $\hat{\f}$ and the algebra of spinor operators
using $\rho$.
\smallskip

\subsubsection{\sf Spinor operators}
If the operator $A$ is explicitly expressed as $\{A_1,A_2\}$,
then $A_1$ and $A_2$ must be applied to the
corresponding components of $f=\{f_1,f_2\}$. In the calculations
below, $A_i$ may contain $s$. Then $A_i$ must be presented
as $A_i'\cdot s$, where $A_i'$ contains no $s$\,; i.e.,
in practice $s$ must be placed on the right. In the operators in
$\mathscr{X}$ we will consider, the commutation relations
between $s$ and $X, \Ga$ must be used when moving $s$. Then
the component
$i$ of $Af$ will be $A_i'(f_{3-i})$, i.e., $s$ placed
on the right means the switch to the other component
before applying $A_i'$.

For instance, $\{\Ga s,\,s-1\}(\{f_1,f_2\})=$
$\{\Ga(f_2), f_1-f_2\}$.

We will frequently use the
vertical mode for spinors:
$$
\{f_1,f_2\}=
\left\{\begin{array}{c}f_1 \\f_2\end{array}\right\},\
\{A_1,A_2\}=
\left\{\begin{array}{c}A_1 \\A_2\end{array}\right\}.
$$

\comment{
\begin{align*}
&X^\rho:\left\{\begin{array}{c}f_{1} \\f_{2}
\end{array}\right\}
\mapsto\left\{\begin{array}{c}Xf_1 \\X^{-1}f_2
\end{array}\right\}
&&\Gamma^\rho:\left\{\begin{array}{c}f_{1} \\f_{2}
\end{array}\right\}
\mapsto\left\{\begin{array}{c}\Gamma(f_1) \\\Gamma^{-1}(f_2)
\end{array}\right\}
\\
&X^\de:\left\{\begin{array}{c}f_{1} \\f_{2}
\end{array}\right\}
\mapsto\left\{\begin{array}{c}Xf_1 \\Xf_2
\end{array}\right\},\
&&\Gamma^\de:\left\{\begin{array}{c}f_{1} \\f_{2}
\end{array}\right\}
\mapsto\left\{\begin{array}{c}\Gamma(f_{1})
\\ \Gamma(f_{2})\end{array}\right\},
\end{align*}
}

\subsubsection{\sf The q-Toda operators via DAHA}
The $q$\~Toda {\em spinor\,} operator is the following
{\em symmetric\,} (i.e., $s$\~invariant) difference
{\em spinor\,} operator
\begin{align}\label{Lspin}
\widehat{\t}\ =\
\{\Gamma^{-1}+(1-X^{-2})\Gamma,\,\Gamma^{-1}+(1-X^{-2})\Gamma\}.
\end{align}
Its first component is the operator $R\!\!E(\mathcal{L})$
from Section \ref{sect:whitlim}; we will use the notation
and definitions
from that section.

We claim that $\widehat{\t}$
can be represented as $\hat{Y}+\hat{Y}^{-1}$
upon the restriction to {\em symmetric spinors}, i.e.,
to $\{f,f\}\in \hat{\f}$. The construction of the
{\em spinor-difference Dunkl operator}
$\hat{Y}$ goes as follows.
\smallskip

Let us introduce the $\de$\~counterpart of $\hbox{\ae}$
from (\ref{remain}):
\begin{align}\label{remainde}
&\hbox{\ae}^\de(\mathcal{A})\equal
(q^{kx}\Ga_{k}^{-1})^\de\,\mathcal{A}\,(\Ga_{k}q^{-kx})^\de,\\
\label{aede}
&\hbox{\ae}^\de:\, X\mapsto \tilde{t}^{\ -1/2}X,\,
\Gamma \mapsto \tilde{t}^{\ -1/2}\Gamma,\, s\mapsto s
\end{align}
for the
{\em spinor constant\,}
$\tilde{t}^{\ 1/2}\equal\{t^{1/2},t^{-1/2}\}$.
Spinor constants are actually diagonal matrices, which may not
commute with $s$ but commutes with $\Ga$ and $X$.
The {\em spinor $R\!\!E$-construction\,} is:
$$
R\!\!E^\de:\ A\mapsto \lim_{t\to 0}\, \hbox{\ae}^\de(A).
$$
It is of course very different from the procedure
$R\!\!E^\rho$ from Section \ref{sect:whitlim}. The spinor-Dunkl
operators are $\hat{Y}=R\!\!E^\de(Y),\, \hat{Y}\,'=
R\!\!E^\de(Y^{-1})$.
They are inverse to each other: $\hat{Y}\hat{Y}\,'=1$.
The following  theorem is from \cite{C102}.

\begin{theorem}\label{MainToda}
The map
\begin{align*}
&Y^{\pm 1}\mapsto \hat{Y}^{\pm 1},\
\breve{\pi}\mapsto R\!\!E^\de(XT),\\
&T\mapsto \hat{T}=R\!\!E^\de(t^{1/2}T),\
T'\mapsto \hat{T}'=R\!\!E^\de(t^{1/2}T^{-1})\,
\end{align*}
can be extended to
a representation of the algebra $\,\overline{\HH}^{\,\vph}\,$
in the space $\,\hat{\mathscr{X}}\,$ of spinors over
$\mathscr{X}=\C[q^{\pm 1/4}][X^{\pm 1}]$.
Correspondingly,
\begin{align*}
&X\mapsto R\!\!E^\de(t^{1/2}X)=R\!\!E^\de(XT)\circ \hat{T}',\\
&X'\mapsto R\!\!E^\de(t^{1/2}X^{-1})=\hat{T}
\circ R\!\!E^\de(XT).
\end{align*}
The commutativity
of\, $T$  and\, $Y+Y^{-1}$\, in $\overline{\HH}^{\,\vph}$ results
in the $s$\~invariance of $\ \hat{Y}+\hat{Y}^{-1}\,$
and the $s$\~invariance of this operator upon its restriction
to the space of $s$\~invariant spinors; this operator is
$\hat{\t}$  from (\ref{Lspin}).\sq
\end{theorem}

It is clear from the construction that all
hat-operators preserve the space of spinors
for the Laurent polynomials in terms of $X^{\pm 1}$.
Upon multiplication by the Gaussian,
it contains the {\em spinor polynomial representation},
isomorphic to the Fourier image of the bar-polynomial
representation times the Gaussian; see
Section \ref{sect:QHermite}, formula (\ref{sitausi})
and Theorem \ref{SPIN-polyn} below.
The reproducing kernel of the isomorphism between
these two modules
inducing $\si:\overline{\HH}\to$ $\overline{\HH}^{\,\vph}$
at the operator level is given by the {\em spinor (nonsymmetric)
$q$\~Whittaker function}.
\medskip

\subsection{\bf Spinor Whittaker function}
\subsubsection{\sf Basic operators}
Let us reproduce some formulas and calculations from
\cite{C102}, which will be needed to extend them
to the case $|q|>1$ and, mainly, to establish the connection
with the core-subalgebra. We will show below that the final
formulas and certain steps in caclulating the limits are
directly related to the core-induced modules.

We begin with $\hat{Y}=R\!\!E^{\de}(Y)=
\lim_{t\to 0}\hbox{\ae}^\de(Y)$.
Using formulas (\ref{aede}):
\begin{eqnarray*}
\hbox{\ae}^\de(Y)
&=&s\cdot(\tilde{t}^{\ -1/2}\Gamma)\cdot
\left(t^{1/2}s+\frac{t^{1/2}-t^{-1/2}}
{\tilde{t}^{-1}X^{2}-1}\cdot (s-1)\right)\\
&=&t^{1/2}\tilde{t}^{\ 1/2}\Gamma^{-1}
+\tilde{t}^{\ 1/2}\Gamma^{-1}\cdot\frac{t^{1/2}-t^{-1/2}}
{\tilde{t}X^{-2}-1}\cdot(1-s)\\
&=&\left\{\begin{array}{c}
t\Gamma^{-1}+\Gamma^{-1}\frac{t-1}{tX^{-2}-1}\cdot(1-s) \\
\Gamma+\Gamma\frac{1-t^{-1}}{t^{-1}X^{2}-1}\cdot(1-s)
\end{array}\right\}\\
\xrightarrow{t\to 0}\ \hat{Y}&=&
\left\{\begin{array}{c}
\Gamma^{-1}\cdot(1-s) \\
\Gamma-\Gamma\cdot X^{-2}\cdot (1-s)\end{array}\right\}.
\end{eqnarray*}
Recall that $\tilde{t}^{\ 1/2}= \{t^{1/2},t^{-1/2}\}$.
A little bit more involved calculation is needed for
$\hat{Y}\,'=R\!\!E^{\de}(Y^{-1})$:

\begin{eqnarray*}
\hbox{\ae}^\de(Y^{-1})
&=&\left(t^{-1/2}s+\frac{t^{-1/2}-t^{1/2}}{\tilde{t}X^{-2}-1}
\cdot (s-1)\right)\cdot
(\tilde{t}^{\ 1/2}\Gamma^{-1}s)\\
&=&\left(\frac{t^{\ -1/2}\tilde{t}X^{-2}-t^{1/2}}
{\tilde{t}X^{-2}-1}\cdot s
-\frac{t^{\ -1/2}-t^{1/2}}{\tilde{t}X^{-2}-1}\right)
\cdot(\tilde{t}^{\ 1/2}\Gamma^{-1}s)\\
&=&\frac{t^{\ -1/2}\tilde{t}X^{-2}-t^{1/2}}
{\tilde{t}X^{-2}-1}\tilde{t}^{\ -1/2}\Gamma
-\frac{t^{-1/2}-t^{1/2}}{\tilde{t}X^{-2}-1}
\tilde{t}^{\ 1/2}\Gamma^{-1}s\\
&=&\left\{\begin{array}{c}
\frac{X^{-2}-1}{tX^{-2}-1}\Gamma-\frac{1-t}{tX^{-2}-1}
\Gamma^{-1}s\\
\frac{t^{-1}X^{2}-t}{t^{-1}X^{2}-1}
\Gamma^{-1}-\frac{t^{-1}-1}{t^{\,-1}X^{2}-1}\Gamma s
\end{array}\right\}\xrightarrow{t\to 0}\\
\hat{Y}\,'&=& \left\{\begin{array}{c}
(1-X^{-2})\Gamma+\Gamma^{-1}s \\
\Gamma^{-1}-\frac{1}{X^{2}}\Gamma s\end{array}\right\}\\
&=& \left\{\begin{array}{c}1-X^{-2} \\1\end{array}\right\}
\Gamma+\left\{\begin{array}{c}1 \\-X^{2}\end{array}\right\}
\Gamma^{-1}s.
\end{eqnarray*}
Automatically, $\hat{Y}\hat{Y}\,'=1$. Now, as we claimed,
\begin{eqnarray*}
&&R\!\!E^{\de}(Y+Y^{-1})=\lim_{t\to 0}\hbox{\ae}^\de(Y+Y^{-1})\\
&=&\left\{\begin{array}{c}\Gamma^{-1}(1-s)+(1-X^{-2})
\Gamma+\Gamma^{-1}s
 \\\Gamma-\Gamma\frac{1}{X^{2}}(1-s)+\Gamma^{-1}
 -\frac{1}{X^{2}}\Gamma s\end{array}\right\}\\
&=&\left\{\begin{array}{c}\Gamma^{-1}+(1-X^{-2})\Gamma \\
\Gamma^{-1}+(1-X^{-2})\Gamma\end{array}\right\}
\bigl(\text{ mod }\bigl(\cdot\bigr)(s-1)\bigr).
\end{eqnarray*}
\smallskip

For $X$ and $X^{-1}$, we have
\begin{align}{\label{eqn:hatX}}
\hat{X}=R\!\!E^{\de}(t^{1/2}X)
=\lim_{t\to 0}\hbox{\ae}^\de(t^{1/2}X)
=\lim_{t\to 0}t^{1/2}\tilde{t}^{\ -1/2}X
=\left\{\begin{array}{c}X \\0\end{array}\right\},\\
\hat{X}\,'=R\!\!E^{\de}(t^{1/2}X^{-1})
=\lim_{t\to 0}\hbox{\ae}^\de(t^{1/2}X^{-1})
=\lim_{t\to 0}t^{1/2}\tilde{t}^{\ -1/2}X^{-1}
=\left\{\begin{array}{c}0 \\X\end{array}\right\}.\notag
\end{align}
Obviously, $\hat{X}\hat{X}\,'=0$.
Next,
\begin{align*}
&\hat{T}=R\!\!E^{\de}(t^{1/2}T)=
\lim_{t\to 0}\hbox{\ae}^\de(t^{1/2}T)
=\left\{\begin{array}{c}0 \\s-1\end{array}\right\},\\
&\hat{T}\,'=R\!\!E^{\de}(t^{1/2}T^{-1})=\lim_{t\to 0}
\hbox{\ae}^\de(t^{1/2}T^{-1})
=\left\{\begin{array}{c} 1 \\ s \end{array}\right\}.
\end{align*}
See \cite{C102} for more formulas and
explicit verifications of the basic relations.
\comment{
\begin{align}
&\hat{T}'=\hat{T}+1,\
\hat{T}\hat{T}\,'\,=\,0\,=\,\hat{T}\,'\hat{T},\ \,
\hat{T}\,'\hat{X}\,'=0=\hat{X}\hat{T},
{\label{eqn:TY1}}\\
&\hat{T}\hat{Y}-\hat{Y}^{-1}\hat{T}\,=\,-\hat{Y},
\ \hat{T}\hat{Y}^{-1}-\hat{Y}\hat{T}\,=\,\hat{Y},
{\label{eqn:TY2}}\\
&\hat{T}\hat{X}-\hat{X}\,'\hat{T}=\hat{X}\,',\
\hat{T}\hat{X}\,'-\hat{X}\hat{T}=-\hat{X}\,',
\ \hat{X}+\hat{X'}=X^{\de}.{\label{eqn:TX}}
\end{align}
Relations \eqref{eqn:TY2} imply that

The relation
\begin{align}\label{hatTYcom}
\hat{T}(\hat{Y}+\hat{Y}^{-1})=(\hat{Y}+\hat{Y}^{-1})\hat{T}
\end{align}
proves that the spinor operator $\hat{Y}+\hat{Y}^{-1}$ is
symmetric (recall that $\hat{Y}\,'=\hat{Y}$).}

\smallskip
\subsubsection{\sf Using the components}
Explicitly, the action of $\hat{Y}$ and $\hat{Y}\,'$
on the spinors is as follows:
\begin{align*}
&\hat{Y}(\left\{\begin{array}{c}f_{1} \\f_{2}
\end{array}\right\})
=\left\{\begin{array}{c}\Gamma^{-1}(f_{1}-f_{2}) \\
\Gamma(f_{2})-\Gamma(\frac{f_{2}-f_{1}}{X^{2}})
\end{array}\right\},
\\
&\hat{Y}\,'(\left\{\begin{array}{c}f_{1} \\f_{2}
\end{array}\right\})
=\left\{\begin{array}{c}(1-X^{-2})
\Gamma(f_{1})+\Gamma^{-1}(f_{2}) \\
\Gamma^{-1}(f_{2})-\frac{1}{X^{2}}\Gamma(f_{1})
\end{array}\right\}.
\end{align*}
It is simple but not immediate to check the relation
$\hat{Y}\hat{Y}\,'= 1$ and other identities for
$\hat{Y}^{\pm1}$
using the component formulas.
The explicit formulas for $\hat{T}$ and $\hat{T}\,'$
are:
\begin{align}\label{Tspinors}
&\hat{T}(\left\{\begin{array}{c}f_{1} \\f_{2}
\end{array}\right\})
=\left\{\begin{array}{c}0\\f_{1}-f_{2}
\end{array}\right\},
&\hat{T}\,'(\left\{\begin{array}{c}f_{1} \\f_{2}
\end{array}\right\})
=\left\{\begin{array}{c}f_{1}\\
f_{1}\end{array}\right\}.
\end{align}
\comment{
It readily gives  \eqref{eqn:TY1}, \eqref{eqn:TY2}.
Generally, there is no need to establish and check
the formulas for $\hat{X}$ and $\hat{X}\,'$ (although
they are simple). Indeed,
\begin{align*}
&\hat{X}\ =\ R\!\!E^\de(XT)\cdot \hat{T}\,',\
\hat{X}\,'\ =\ \hat{T}
\cdot R\!\!E^\de(XT).
\end{align*}
Thus we need only to know $\hat{\pi}\equal
R\!\!E^\de(XT)$.}

Finally,
\begin{eqnarray*}
\hbox{\ae}^\de(XT)
&=&(\tilde{t}^{-1/2}X)(t^{1/2}s
+\frac{t^{1/2}-t^{-1/2}}{\tilde{t}^{-1}X^{2}-1}(s-1))\\
&=&\tilde{t}^{-1/2}t^{1/2}Xs
+\frac{X(\tilde{t}^{-1/2}t^{1/2}-\tilde{t}^{-1/2}t^{-1/2})}
{\tilde{t}^{-1}X^{2}-1}(s-1)\\
&=&\left\{\begin{array}{c}Xs \\tX^{-1}s\end{array}\right\}
+\left\{\begin{array}{c}\frac{X(1-t^{-1})}
{t^{-1}X^{2}-1}(s-1) \\
\frac{X^{-1}(t-1)}{tX^{-2}-1}(s-1)\end{array}\right\}.
\end{eqnarray*}
Taking the limit $t\to 0$,
\begin{align*}
\hat{\pi}=\left\{\begin{array}{c}Xs \\0\end{array}\right\}
+\left\{\begin{array}{c}-X^{-1}(s-1) \\X^{-1}(s-1)
\end{array}\right\}
=\left\{\begin{array}{c}Xs-X^{-1}(s-1) \\X^{-1}(s-1)
\end{array}\right\}.
\end{align*}
Using the components,
\begin{align}\label{picompo}
\hat{\pi}:\ \left\{\begin{array}{c}f_{1} \\f_{2}
\end{array}\right\}\mapsto
\left\{\begin{array}{c}Xf_{2}+\frac{f_{1}-f_{2}}{X} \\
\frac{f_{1}-f_{2}}{X}\end{array}\right\}.
\end{align}
Check directly that  $\hat{\pi}^{2}=id$.

This formula completes the ``component presentation" of the
{\em hat-module} of $\overline{\HH}^{\,\vph}$
from Theorem \ref{MainToda}:
$$
T, \breve{\pi}, Y\ \mapsto\
 \hat{T}, \hat{\pi}, \hat{Y}.
$$

\comment{
The extension of this Theorem to arbitrary (reduced)
root systems is straightforward as well as the justification;
we will address it (and the applications) in further paper(s).
The formulas for the $\overline{Y}$\~operators are of course
getting more involved. The justifications
in the spinor $q$\~Toda theory (including global Whittaker
functions) are entirely based on the DAHA-theory.
We calculate and check practically everything
explicitly in this work mainly to demonstrate the practical
aspects of the technique of spinors (and because of novelty
of this topic).}

\subsubsection{\sf The main formula}
Let us apply the procedure $R\!\!E^\de$ to
the {\em global difference spherical function}
$G$ from (\ref{gxla}). It was denoted by
$\e_q(x,\la)$ in \cite{C5}, Section 5 (arbitrary
reduced root systems). See \cite{Sto}
for the $C^\vee C$\~case. We need to use
Proposition \ref{tildee}
to express the conjugated $E$\~polynomials in terms
of $E$\~polynomials.  Generally, this relation
requires the action of $T_{w_0}$.

Following \cite{C102}, we arrive at the
spinor (nonsymmetric) generalization
of $\w$ from (\ref{Whitsym})\,:
\begin{align}
\Om(X,\La)= (\tga'(X)\tga'(\La))^{-1}\,\Bigl(1+
\sum_{m=1}^{\infty} q^{m^{2}/4}\,
\bigl(&\frac{\overline{E}_{-m}(\La)}
{\prod_{s=1}^{m} (1-q^{s})}
\left\{\begin{array}{c}X^m \\ q^mX^m\end{array}\right\}\notag\\
+\,&\frac{\overline{E}_{m}(\La)}
{\prod_{s=1}^{m-1} (1-q^{s})}
\left\{\begin{array}{c} 0 \\ X^m\end{array}\right\}
\bigr)\Bigr).\label{spinwhito}
\end{align}
Using the Pieri rules from (\ref{pienilp-}), it can be
presented as follows:
\begin{align}\label{spinwhit}
\Om= (\tga'(X)\tga'(\La))^{-1}\,
\sum_{m=0}^{\infty} \frac{q^{m^{2}/4}}
{\prod_{s=1}^{m} (1-q^{s})}
\left\{\begin{array}{c}X^m \overline{E}_{-m}(\La) \\
X^m \La^{-1}\overline{E}_{m+1}(\La)\end{array}\right\}.
\end{align}

Either of these two presentations readily gives that the
spinor symmetrization of $\Om$ is $\{\w,\,\w\}$
for the symmetric $q$\~Whittaker function $\w$.
We need
to apply the symmetrizer $T'=T+1$ to $\Om$,
equivalently, make the second component
equal to the first one;
see (\ref{Tspinors}).
Note that $\La$ is a non-spinor variable.

The spinor function $\Om$
intertwines the bar-representation of  $\overline{\HH}$ and
the hat-representation of  $\overline{\HH}^{\,\vph}$. Namely,
\begin{align}\label{omident}
&\hat{Y}(\Om)=\La^{-1}(\Om),\ \
\hat{X}(\Om)=\overline{Y}'_\La (\Om),\ \
\hat{X}\,'(\Om)=\overline{Y}_\La (\Om),\\
&\hat{\pi}(\Om)\ =\ \pi_\La(\Om),\  \,
\hat{T}(\Om)\ =\  \overline{T}_\La(\Om),\,\
\hat{T'}(\Om)\ =\  \overline{T'}_\La(\Om),
\label{Omegainter}
\end{align}
where $\overline{Y}\,'_\La,\, \overline{Y}_\La,\,\pi_\La,\,
\overline{T}_\La$\,
act on the argument $\La$; the other operators are
$X$\~operators. Recall that $\hat{\pi}$ is the action
of $\breve{\pi}=XT$ in the spinor representation.

These (and other
related identities) follow from the general
theory for arbitrary reduced root systems (at least, in
the twisted case). However, in the rank one case
(and for $A_n$), one can use the Pieri
rules from (\ref{pienilp+}),(\ref{pienilp-})
and formulas (\ref{Ynilp+}),
(\ref{Ynilp-}) for the direct verification.
See \cite{C102} for explicit calculations.
For instance,
\begin{align}\label{hatYga}
\tga'(X)\,\hat{Y}\,\tga'(X)^{-1}\left\{\begin{array}{c}f_1\\
f_2\end{array}\right\}=q^{\frac{1}{4}}
\left\{\begin{array}{c} X^{-1}\Ga^{-1}(f_1-f_2) \\
X\Ga(f_2)+q^{-1}\frac{\Ga(f_1-f_2)}{X}\end{array}\right\}.
\end{align}

Formulas (\ref{hatYga}), (\ref{picompo}) and (\ref{Tspinors})
were used in \cite{C102} to introduce the {\em spinor-polynomial}
representation.

\begin{theorem}\label{SPIN-polyn}
The space
$$
\mathscr{X}_{spin}\equal\C\oplus\ \bigl
(\oplus_{m=1}^{\ \infty} (\C\{X^m,0\}\oplus
\C\{0,X^m\})\bigr).
$$
is an irreducible
$\overline{\HH}^{\,\vph}$\~submodule of the space of
spinors over $\C[X^{\pm1}]$ supplied with the twisted
action:
$$
\overline{\HH}^{\,\vph}\ni A\mapsto q^{-x^2}\,\hat{A}\,q^{x^2}.
$$
Equivalently, $\mathscr{X}_{spin}$
is invariant and irreducible under the
action of operators $\hat{T},
\hat{\pi}$  and $q^{-x^2}\,\hat{Y}\,q^{x^2}$. \sq
\end{theorem}
\medskip

\subsection{\bf Algebraic theory}
\subsubsection{\sf Relation to tilde-modules}
It is important that the $\overline{\HH}^{\,\vph}$\~module
$\mathscr{X}_{spin}$ can be identified with the induced
$\overline{\HH}$\~module
$\tilde{\mathscr{X}}_{\pi,\tilde{X}}^{\,0,-}$
(where $\ep=0$ and $\de=-1$). The identification goes as follows.

\begin{theorem}\label{TILDESPIN}
Let us define the $\C$\~linear map
$\,\chi:\tilde{\mathscr{X}}_{\pi,\tilde{X}}^{0,-}
\rightarrow \mathscr{X}_{spin}$ by
\begin{align}\label{spinidentif}
1\mapsto \{1,1\},\ Y^m \mapsto \{X^m,0\},\
(Y')^m \mapsto \{0,X^m\},\ q\mapsto q^{-1}.
\end{align}
for $m>0$. It induces the following isomorphism on operators:
\begin{align}\label{chiso}
&\chi:\overline{\HH} \rightarrow \overline{\HH}^{\,\vph},\,
T \mapsto
-T',\ \pi \mapsto -\breve{\pi},\ X^{\pm 1}\mapsto Y^{\pm 1},\
q \mapsto q^{-1},\\
&\chi(T')=1-T'=-T,\ \chi(Y)=\breve{\pi}T'=X,\
\chi(Y')=T\breve{\pi}=X'.\notag
\end{align}
\end{theorem}
{\it Proof.}
Formula (\ref{hatYga}) gives that
$q^{-x^2}\widehat{Y}q^{x^2}(\{1,1\})=q^{1/4}\{0,X\}$,
\begin{align}\label{hatYgaex}
q^{-x^2}\widehat{Y}q^{x^2}&\left\{\begin{array}{c}X^m\\
0\end{array}\right\}
=q^\frac{1}{4}\left\{\begin{array}{c}X^{-1}q^{-\frac{m}{2}}X^m\\
q^{-1}X^{-1}q^\frac{m}{2}X^m\end{array}\right\}\\
=q^{-\frac{1}{4}}(
q^{-\frac{m-1}{2}}&\left\{\begin{array}{c}X^{m-1}\\
0\end{array}\right\}+q^\frac{m-1}{2}\left\{\begin{array}{c}0\\
X^{m-1}\end{array}\right\}),\notag\\
q^{-x^2}\widehat{Y}q^{x^2}&\left\{\begin{array}{c}0\\
X^m\end{array}\right\}
=q^\frac{1}{4}
\left\{\begin{array}{c}-X^{-1}q^{-\frac{m}{2}}X^m\\
Xq^\frac{m}{2}X^m-
q^{-1}X^{-1}q^\frac{m}{2}X^m\end{array}\right\}\notag\\
=q^{-\frac{1}{4}}\,(
q^\frac{m+1}{2}&\left\{\begin{array}{c}0\\X^{m+1}\end{array}\right\}
-q^{-\frac{m-1}{2}}\left\{\begin{array}{c}X^{m-1}\\
0\end{array}\right\}\notag\\
&\quad -q^\frac{m-1}{2}
\left\{\begin{array}{c}0\\X^{m-1}\end{array}\right\}),
\hbox{\ \  where } m \geq 1.\notag
\end{align}
Compare with (\ref{yactionx}); the action of
$Y \in \overline{\HH}^\varphi$ coincides with that of
$X \in \overline{\HH}$ upon the identification given by
(\ref{spinidentif}).

Next, the operator $\widehat{T}' = \{1,s\}$
coincides with $-T$ defined by (\ref{tsumpol}):
\begin{align}
-T\,(\sum_{l\ge 0}& a_{l+1} Y^{l+1}+ a_{-l} (Y')^l)(1)\\
&=\sum_{l\ge 0}\, a_{l+1}(Y^{l+1}+(Y')^{l+1})(1).\notag
\end{align}

Finally,
the component presentation of $\widehat{\pi}$ from
(\ref{picompo}) results in
\begin{align}
&\widehat{\pi}(\left\{\begin{array}{c}1\\1\end{array}\right\})
=\left\{\begin{array}{c}X\\0\end{array}\right\},\
\widehat{\pi}(\left\{\begin{array}{c}X^m\\0\end{array}\right\})
= \left\{\begin{array}{c}X^{m-1}\\X^{m-1}\end{array}\right\},\\
&\widehat{\pi}(\left\{\begin{array}{c}0\\X^m\end{array}\right\})
=\notag
\left\{\begin{array}{c}X^{m+1}\\0\end{array}\right\}
-\left\{\begin{array}{c}X^{m-1}\\X^{m-1}\end{array}\right\}.
\end{align}
This coincides with $-\pi \in \overline{\HH}$ defined by
(\ref{pisumpol}). \sq

We see that $\chi$ is essentially an involution and can be
used equally well to go from $\overline{\HH}^{\,\vph}\,$ to
$\overline{\HH}$. However we prefer to use $\chi^{-1}$
for the inverse map. Note that $\chi$ is nothing but
$$
\nu_{-\,}\vep:\, \overline{\HH} \rightarrow
\overline{\HH}^{\,\vph\dag}\rightarrow \overline{\HH}^{\,\vph},
$$
where $\nu_{\pm\,}$ are defined in (\ref{nudag}),
$\, \vep$ is the bar-restriction of that from (\ref{vepanti}),
naturally sending
\begin{align}\label{vepbar}
&\vep: \overline{\HH} \rightarrow \overline{\HH}^{\,\vph\dag}:\ \,
T\mapsto (T^\dag)',\, X\mapsto Y,\, Y\mapsto X,\, q\mapsto q^{-1}.
\end{align}
\smallskip

\subsubsection{\sf Algebraic Whittaker function}
Let us apply the identification $\chi^{-1}$ from
(\ref{spinidentif}) to $\Om(X,\Lambda)$.
We obtain
\begin{align}\label{omalg}
\Om_\mathrm{alg}(Y,Y';\La) \equal
&q^{-\la^2}(1+
\sum_{m=1}^\infty q^{-\frac{m^2}{4}}\left(
\frac{(\overline{E}_{m}(\La))^\ast}
{\prod_{s=1}^{m-1} (1-q^{-s})}(Y')^m\right.\\
&+\notag\left.\frac{(\overline{E}_{-m}(\La))^\ast}
{\prod_{s=1}^{m}(1-q^{-s})}
(Y^m + q^{-m}(Y')^m)\right)).
\end{align}
Here it is convenient to use the Gaussian
$q^{-\la^2}$ instead of $\widetilde{\gamma}'(\La)$;
it is understood as a formal symbol satisfying the standard
relations.
We drop the Gaussian for the spinor variables because it is already
incorporated into the representation and the conjugation
by this function is included in $\chi$.

The coefficients of the summation in (\ref{omalg})
belong to
$
\overline{\mathscr{X}}^\dag_\La\otimes
\tilde{\mathscr{X}}_{\pi,\tilde{X}}^{\,0,-}.
$
Note $\dag$, which appears due to the conjugation of
$\overline{E}_m(\La)$;\
$\overline{\mathscr{X}}^\dag_\La $ is a
module over $\overline{\HH}^\dag$ (not over
 $\overline{\HH}$).

Applying $\chi^{-1}$ to
(\ref{omident}), the relations satisfied by
$\Omega_\mathrm{alg}$ are
\begin{align}\label{omalgeqn}
\chi^{-1}(H^\vph)(\Om_\mathrm{alg})
&= \eta(H_\La)(\Om_\mathrm{alg}), \where H\in \overline{\HH}.
\end{align}
Here $\vph$ is from (\ref{vphnil}) and $\eta$ is from (\ref{etaiso}).
Using that
$$
\chi^{-1}=\vep\nu_{-\,}=\nu_{-\,}\vep \and
\nu_{-\,}\vep\vph\eta=\nu_{-}\ast\eta=\nu_{-\,}\diamond,
$$
we can rewrite (\ref{omalgeqn}) as follows:
\begin{align}\label{omalgeqnn}
(H^\diamond)(\Om_\mathrm{alg})
&= \nu_{-\,}(H_\La)(\Om_\mathrm{alg}), \where H\in \overline{\HH}.\
\end{align}

Explicitly,
\begin{align}\label{omalgident}
&X(\Om_\mathrm{alg})=\La(\Om_\mathrm{alg}),\ \ \ \
-\pi(\Om_\mathrm{alg})=\pi_\La(\Om_\mathrm{alg}),\\
-&T(\Om_\mathrm{alg})=\overline{T}^\dag_\La(\Om_\mathrm{alg}),\ \
\notag
-T'(\Om_\mathrm{alg})=(\overline{T}^\dag_\La-1)(\Om_\mathrm{alg}),\\
&Y(\Om_\mathrm{alg})=\overline{T}^\dag_\La\pi(\Om_\mathrm{alg}),\
\notag\
Y'(\Om_\mathrm{alg})=\pi(\overline{T}^\dag_\La-1)(\Om_\mathrm{alg}).
\end{align}

\comment{
\subsubsection{\sf An explicit verification}
These formulas can be of course checked directly.
The second formula in (\ref{omalgident}) is the most interesting
example. We will use the $\ast$-variants of
the Pieri rules from (\ref{pienilp+}), namely
\begin{align}\label{starpieri}
\La\overline{E}_{-n}^\ast
&=\overline{E}_{-(n+1)}^\ast-\overline{E}_{n+1}^\ast \ (n\geq 0),\\
\La\overline{E}_n^\ast
&=(1-q^{-(n-1)})\overline{E}_{n-1}^\ast\notag
+q^{n-1}\overline{E}_{-(n-1)}^\ast \ (n\geq 1).
\end{align}

First, use (\ref{starpieri}) to rewrite $\Om_\mathrm{alg}$ as
\begin{align}\label{omalgpieri}
\Om_\mathrm{alg}(Y,Y';\La)&=q^{-\la^2}(1+
\sum_{m=1}^\infty q^{-\frac{m^2}{4}}\times\\
\left(\frac{\La(\overline{E}_{m+1}(\La))^\ast}
{\prod_{s=1}^m (1-q^{-s})}\right.&(Y')^m\notag
+\left.\frac{(\overline{E}_{-m}(\La))^\ast}
{\prod_{s=1}^{m}(1-q^{-s})}Y^m\right))\\
&=\notag q^{-\la^2}
\sum_{m=0}^\infty \frac{q^{-m^2/4}}
{\prod_{s=1}^{m} (1-q^{-s})}\mathcal{F}_m,
\end{align}
where $\mathcal{F}_0 \equal 1,\ \mathcal{F}_m \equal
\La(\overline{E}_{m+1}(\La))^\ast (Y')^m
+(\overline{E}_{-m}(\La))^\ast Y^m,\ m\geq 1$.
Note that $\pi_\La(q^{-\la^2}) = q^{-1/4}\La q^{-\la^2}$ and
$\pi_\La(\overline{E}_{m}^\ast)=
q^{-m/2}\La\overline{E}_{1-m}^\ast$.
Hence,
\begin{align}\label{pilaFga}
\pi_\La(\mathcal{F}_m q^{-\la^2}) &=
q^{-\frac{1}{4}}\La q^{-\la^2}(
q^{-\frac{m}{2}}\overline{E}_{-m}^\ast (Y')^m
+q^\frac{m}{2}\La\overline{E}_{m+1}^\ast Y^m),\\
\pi_\La(\mathcal{F}_0 q^{-\la^2})&=\notag
q^{-\frac{1}{4}}\La q^{-\la^2}.
\end{align}

Now let us apply $-\pi$ using (\ref{pisumpol}). We find
\begin{align}\label{minuspiF}
-\pi(\mathcal{F}_m) &= \La\overline{E}_{m+1}^\ast(Y^{m+1}-M_{m-1})
+\overline{E}_{-m}^\ast M_{m-1},
\end{align}
where $M_{-1} \equal 0,\ M_0 \equal 1,\
M_m \equal Y^m+(Y')^m,\ m\geq 1$.
Thus,
\begin{align}\label{minuspiomalg}
-\pi(\Om_\mathrm{alg}) &= q^{-\la^2}(q^{-1/4}\La
+\sum_{m=1}^\infty (q^{-(m+1)^2/4}\times\\
\frac{(-\La\overline{E}_{m+2}^\ast
+\overline{E}_{-(m+1)}^\ast)}{\prod_{s=1}^{m+1}(1-q^{-s})}
&(Y^m+(Y')^m)+q^{-(m-1)^2/4}\notag
\frac{\La\overline{E}_m^\ast}{\prod_{s=1}^{m-1}(1-q^{-s})}Y^m).
\end{align}

Using (\ref{starpieri}), we simplify:
\begin{align}
-\La\overline{E}_{m+2}^\ast +\overline{E}_{-(m+1)}^\ast &=\notag
-(1-q^{-(m+1)})\overline{E}_{m+1}^\ast
+(1-q^{-(m+1)})\overline{E}_{-(m+1)}^\ast\\
&=\notag (1-q^{-(m+1)})\La\overline{E}_{-m}^\ast.
\end{align}
Hence
\begin{align}\label{minuspiomalgg}
-\pi(\Om_\mathrm{alg}) &= q^{-\la^2}(q^{-1/4}\La
+\sum_{m=1}^\infty q^{-(m+1)^2/4}\times\\
\frac{\La\overline{E}_{-m}^\ast}
{\prod_{s=1}^{m}(1-q^{-s})}&(Y')^m\notag
+\notag q^{-(m-1)^2/4}
\frac{\La^2\overline{E}_{m+1}^\ast}
{\prod_{s=1}^{m}(1-q^{-s})}Y^m),\\
&=\notag q^{-1/4}\La q^{-\la^2}(1
+\sum_{m=1}^\infty\frac{q^{-m^2/4}}
{\prod_{s=1}^{m}(1-q^{-s})}\times\\
&\notag(q^{-m/2}\overline{E}_{-m}^\ast(Y')^m
+q^{m/2}\La\overline{E}_{m+1}^\ast Y^m))
\end{align}
Comparing this with (\ref{omalgpieri}) and (\ref{pilaFga})
shows that $\pi_\La(\Om_\mathrm{alg})=-\pi(\Om_\mathrm{alg})$.
}

\begin{theorem}\label{TILDEOM}
(i) Let us define the transform $\B:
\overline{\mathscr{X}}_{\La}q^{-\la^2}\, \rightarrow\,
\tilde{\mathscr{X}}_{\pi,\tilde{X}}^{0,-}$ by
\begin{align}\label{transformalg}
\B\bigl(f(\La)q^{-\la^2}\bigr)(Y,Y')\equal
(\!(\,\Omega_{\mathrm{alg}}(Y,Y';\La)\, ,\,
f(\La)q^{-\la^2}\, )\!)_\circ
\end{align}
for $(\!(\ ,\ )\!)_\circ$ from (\ref{barinner}).
Then it induces
\,$\nu_{-\,}\eta:\,\overline{\HH}\,\to\,
\overline{\HH}^{\vph}$ on operators:
\begin{align}\label{transforma}
\B(H(f(\La)q^{-\la^2}))
\ =\ (\nu_{-\,}\eta(H))\bigl(\B(f(\La)q^{-\la^2})\bigr).
\end{align}

(ii)
Using $M_m^Y=Y^m+(Y')^m \for m> 0$
and $M^Y_0=1$,
\begin{align}\label{btranpol}
\B(\,
\overline{E}_{-m}(\La) q^{-\la^2}\,)\ &=\
q^{-\frac{m^2}{4}}M^Y_{m}(1) \for m\ge 0,\\
\B(\,\overline{E}_{+m}(\La)q^{-\la^2}\,)\ &=\
q^{-\frac{m^2}{4}}Y^{m}(1) \for m\ge 1.\notag
\end{align}
We use here $H(1)$ in $\tilde{\mathscr{X}}_{\pi,\tilde{X}}^{0,-}$;
see the identification from (\ref{spindirect}).
\end{theorem}

\subsubsection{\sf Justification}
Theorem \ref{TILDEOM} is essentially a variant of
Theorem \ref{FOURE}, ($ii$) adjusted to the nil-setting.
It induces $\nu_{-}\eta$ on operators by construction. However,
part ($ii$), concerning the explicit formulas,
is not quite obvious.
The following
transformation of the expression in (\ref{omalg})
is needed.

\begin{lemma}\label{LEMOMALG}
In terms of the polynomials $\{\overline{E}_{m}^\dag\}$
from (\ref{astepol}),
\begin{align}\label{omalgdag}
\Om_\mathrm{alg}(Y,Y';\La)=
q^{-\la^2}\Bigr(
&\sum_{m=1}^\infty q^{-\frac{(m-1)^2}{4}}M^Y_{m-1}(1)
\frac{\overline{E}_{1-m}^\dag(\La)}
{\prod_{s=1}^{m-1} (1-q^{-s})}\\
+&\sum_{m=0}^\infty q^{-\frac{(1+m)^2}{4}}Y^{m+1}(1)
\frac{\overline{E}_{1+m}^\dag(\La)}
{\prod_{s=1}^{m} (1-q^{-s})}\Bigr).\notag
\end{align}
\end{lemma}
{\it Proof.} To avoid misunderstanding, let us
calculate a couple of first terms in the right-hand side.
The term with $m=1$ in the first summation is $1$.
The term with $m=0$ in the second summation is
$q^{-1/4}Y(1)\,E_1^\dag(\La)=q^{-1/4}\La Y$,
since $\overline{E}_1^\dag(\La)=\La$. There is another appearance
of $Y$ due to $M_1^Y$; the corresponding term is
$$
q^{-\frac{1}{4}}Y(1)
\frac{\overline{E}_{-1}^\dag(\La)}{1-q^{-1}}, \where
\overline{E}_{-1}^\dag(\La)=\La^{-1}+q^{-1}\La.
$$
Thus the total coefficient of $Y$ is
$$
q^{-1/4}\bigl((1+\frac{q^{-1}}{1-q^{-1}})\La+\frac{1}{1-q^{-1}}
\La^{-1}\bigr)=q^{-1/4}\frac{\La+\La^{-1}}{1-q^{-1}}.
$$
This is exactly the coefficient of $Y$ from
the original formula (\ref{omalg}), which is
$q^{-1/4}\overline{E}_{-1}(\La)^\ast/(1-q^{-1})$.

The general verification is based on the
formulas
$$\overline{E}_{m}(\La)^\ast=
\La^{-1}\overline{E}_{1-m}^\dag,\
X(\Om_\mathrm{alg})=\La(\Om_\mathrm{alg}).
$$
One has:
\begin{align}\label{omalgg}
&\Om_\mathrm{alg}(Y,Y';\La)=
\La X^{-1}(\Om_\mathrm{alg}(Y,Y';\La))=\\
&q^{-\la^2}\Bigl(\La X^{-1}(1)+
\sum_{m=1}^\infty q^{-\frac{m^2}{4}}\left(
\frac{\overline{E}_{1-m}^\dag(\La)}
{\prod_{s=1}^{m-1} (1-q^{-s})}X^{-1}((Y')^m)\right.\\
&+\notag\left.\frac{\overline{E}_{1+m}^\dag(\La)}
{\prod_{s=1}^{m}(1-q^{-s})}
X^{-1}(Y^m + q^{-m}(Y')^m)\right)(1)\Bigr).
\end{align}
Now we can use the formulas (\ref{yactionx}) for $\de=-1$:
\begin{align}\label{yactionxx}
&X^{-1}(Y') = (-\de) q^{\frac{1}{4}},\
X^{-1}(1) = (-\de) q^{-\frac{1}{4}}Y(1),\\
X^{-1}&((Y')^m) = (-\de) q^{-\frac{1}{4}}\,
q^{+m/2}(\,Y^{m-1}+(Y')^{m-1}\,)(1),\notag\\
X^{-1}&(Y^m+q^{-m}(Y')^m)= (-\de) q^{-\frac{1}{4}}q^{-m/2}\,
Y^{m+1}(1)\  (m\geq 2).\notag
\end{align}
It concludes the verification of the lemma.
\sq
\smallskip

According to the general approach from
(\ref{reproker}), the series
from (\ref{omalgdag}) is exactly the reproducing
kernel for $\B$ from the theorem for
the basis $\{f_n=\overline{E}_n(\La) q^{-\la^2}\}$ and
basis $\{f_n'=\overline{E}^\dag_n(\La) q^{-\la^2}\}$ orthogonal
to each other with respect to $(\!(\,,\,)\!)_\circ.$
Thus,
\begin{align}\label{btransom}
&\Om_\mathrm{alg}(Y,Y';\La)=
\sum_{m\in \Z}
\frac{\overline{E}^\dag_m(\La)\,q^{-\la^2}\cdot
\B(\overline{E}_m(\La)\,q^{-\la^2})}
{(\!(\,\overline{E}_m^\dag(\La) q^{-\la^2}\ ,\
\overline{E}_m(\La) q^{-\la^2}\,)\!)_\circ},\\
&(\!(\,\Om_\mathrm{alg}(Y,Y';\La)\,,\,
q^{-\la^2}\overline{E}_m(\La)\,)\!)=
\B(q^{-\la^2}\overline{E}_m).
\end{align}
See (\ref{innerepol}) for the values of
$(\!(\,\overline{E}_m^\dag q^{-\la^2}\,,\,
\overline{E}_m q^{-\la^2}\,)\!)_\circ=
(\!(\,\overline{E}_m^\dag\,,\,
\overline{E}_m(\La)\,)\!)_\circ.$
\sq
\smallskip

It is instructional to check directly the following
corollary of the theorem. Since,
$\tau_+^{-1}(Y)_\La=q^{1/4}X^{-1}Y$ and
$\tau_+^{-1}(Y')_\La=q^{-1/4}Y'X$
diagonalize $\overline{E}_n(\La)q^{-\la^2}$ the operators
\begin{align}\label{yniceop}
&\y=\eta\tau_+^{-1}(Y)=q^{-1/4}X\pi T^{-1}=q^{1/4}\pi TX=
q^{1/4}YX,\\
&\y^{\,\prime}=
\eta\tau_+^{-1}(Y')=q^{-1/4}X^{-1}Y',\notag
\end{align}
must diagonalize the images of
$\overline{E}_n(\La)q^{-\la^2}$ under
$\B$. Here $\nu_{-}$ acts trivially on $X,Y$.

More exactly, the eigenvalues must be preserved
too. Due to (\ref{Ynilp+}) and (\ref{Ynilp-}):
\begin{align}\label{Ynilp+n}
&\overline{Y}\,(\overline{E}_{n})
=\left\{\begin{array}{ccc}q^{-|n|/2}\overline{E}_{n}, &  & n>0, \\
0, &  & n\le 0.\end{array}\right\},\\
&\overline{Y}'\,(\overline{E}_{n})
=\left\{\begin{array}{ccc}q^{-|n|/2}\overline{E}_{n},
&  & n\le 0, \\
0, &  & n>0.\end{array}\right\}.\label{Ynilp-n}
\end{align}
So we must have:
\begin{align}\label{Ynilp+nn}
&\y\,(Y^n)=q^{-n/2}Y^n\, (n>0),\ \ \ \y\,(M_n^Y)=0\, (n\ge 0),\\
&\y^{\,\prime}\,(M_n^Y)=q^{-n/2}Y^n\, (n\ge 0),\
\y^{\,\prime}\,(Y^n)=0\, (n>0).
\label{Ynilp-nn}
\end{align}
This readily follows from (\ref{yactionxx}) for $\y^{\,\prime}$.
For $\y$, use these relations multiplied by $X$:
\begin{align}\label{yactionxxx}
& X(1)=(-\de) q^{-\frac{1}{4}}Y',\
X(Y)=(-\de) q^{\frac{1}{4}},\\
& X(\,Y^{m-1}+(Y')^{m-1}\,)=
(-\de) q^{\frac{1}{4}}q^{-m/2}(Y')^m,\notag\\
&X(Y^{m+1})=(-\de) q^{\frac{1}{4}}q^{m/2}(Y^m+q^{-m}(Y')^m)\,
\  (m\geq 2).\notag
\end{align}
\smallskip
More conceptually, follow Lemma \ref{TILDEINTW}.
\smallskip

\subsubsection{\sf Omitting the conjugation}
According to the general formula from
(\ref{reprokerpro}), if we use the
pairing $\lan f, g\ran_\circ$ from (\ref{innerinvw})
corresponding to the anti-involution $\diamond$,
the corresponding transform will represent
$\nu_{-\,}$, i.e., it will essentially be an isomorphism
of $\overline{\HH}$ modules. Let us demonstrate this now.
We will use directly the presentation
of $\Om_{\mathrm{alg}}$ from
(\ref{omalg}); there is no need for (\ref{omalgdag}).

Now $\{f_n=\overline{E}^\dag_n(\La) q^{\la^2}\}$ and
$\{f_n'=(\overline{E}_n(\La) q^{\la^2})^*\}$ are orthogonal
to each other with respect to $\lan\,,\,\ran_\circ.$
We come to the following counterpart of
Theorem \ref{TILDEOM}.

\begin{theorem}\label{TILDEOMX}
(i) Let us define $\F:
\overline{\mathscr{X}}^\dag_{\La}q^{+\la^2}\, \rightarrow\,
\tilde{\mathscr{X}}_{\pi,\tilde{X}}^{0,-}$ by
\begin{align}\label{transfalgx}
\F\bigl(f(\La)q^{+\la^2}\bigr)(Y,Y')\equal
\lan\,\Omega_{\mathrm{alg}}(Y,Y';\La)\, ,\,
f(\La)q^{-\la^2}\, \ran_\circ
\end{align}
for $\lan\ ,\ \ran_\circ$ from (\ref{innerinvw}).
Then
\begin{align}\label{transformax}
&\F(H(f(\La)q^{+\la^2}))
\ =\ \nu_{-\,}(H)\bigl(\F(f(\La)q^{+\la^2})\bigr),\\
\label{btranpolx}
&\F(\,
\overline{E}_{-m}^\dag(\La) q^{+\la^2}\,)\ =\
q^{-\frac{m^2}{4}}(Y^{m}+q^{-m}(Y')^m)(1) \for m\ge 0,\\
&\F(\,\overline{E}^\dag_{+m}(\La)q^{+\la^2}\,)\ =\
q^{-\frac{m^2}{4}}(Y')^{m}(1) \for m\ge 1.\notag
\end{align}
\sq
\end{theorem}

Similar to the explicit calculation above, let us check
that the eigenfunctions really go to the corresponding
eigenfunctions under $\F$. The polynomials
$\overline{E}_n^\dag(\La)=\La (\overline{E}_{1-n})^*$
solve the following eigenvalue problem:
\begin{align}\label{Ynilp+dag}
&(\overline{Y}^\dag)'\,(\overline{E}^\dag_{n})
=\left\{\begin{array}{ccc}q^{|n|/2}\overline{E}^\dag_{n}, &  & n>0, \\
0, &  & n\le 0.\end{array}\right\},\\
&\overline{Y}^\dag\,(\overline{E}^\dag_{n})
=\left\{\begin{array}{ccc}q^{|n|/2}\overline{E}^\dag_{n},
&  & n\le 0, \\
0, &  & n>0.\end{array}\right\}.\label{Ynilp-dag}
\end{align}
Correspondingly,
\begin{align}\label{ynicedag}
&\tilde{\y}=\nu_{-\,}\tau_+(Y^\dag)=q^{-1/4}XY,\
&\tilde{\y}^{\,\prime}=
\nu_{-\,}\tau_+(Y')^\dag=q^{+1/4}Y'X^{-1},
\end{align}
and we must have for $n\ge 0$:
\begin{align}\label{Ynilp+dagg}
&\tilde{\y}^{\,\prime}\,((Y')^{n+1})=q^{(n+1)/2}(Y')^{n+1},\ \ \
\ \ \
\tilde{\y}^{\,\prime}\,(Y^{n}+q^{-n}(Y')^n)=0,\\
&\tilde{\y}\,(Y^{n}+q^{-n}(Y')^n)=q^{n/2}(Y^{n}+q^{-n}(Y')^n), \
\tilde{\y}\,((Y')^{n+1})=0.
\label{Ynilp-dagg}
\end{align}
This is straightforward.

Generally, the map we constructed
is closely connected with those that can be obtained directly
from Proposition \ref{PITILDEY} and especially
Theorem \ref{EPOLTILDE}. As we
see, the real source of the function $\Omega$ appears in
the fact that the core subalgebra provides natural creation
and annihilation operators for the $\overline{E}$\~polynomials.
\medskip

\setcounter{equation}{0}
\section{\sc The case of large \texorpdfstring{$q$}{q}}
Let us address  the spinor limit of the ``second" function
$G^\checkmark$ defined for $|q|>1$ by formula (\ref{gxlaprime}).
Its theory is somewhat more involved than that for $|q|<1$.

\subsection{\bf Omega-check function}
The procedure
is based on the delta-version
of $R\!\!E^\checkmark$ from (\ref{reprime})
{\em twisted} as follows.

\subsubsection{\sf Twisted RE-procedure}
Setting $\tilde{\Ga_{k}}=\Ga_k^\de,\  \tilde{q^{kx}}=
(q^{kx})^\de$, let
\begin{align}\label{reprimediag}
\tilde{\hbox{\ae}}(\a)&\equal
\left\{\begin{array}{c}1 \\t^{-1}
\end{array}\right\}
\tilde{q^{-kx}}\tilde{\Ga_{k}}^{-1}\,\a\,
\tilde{\Ga_{k}}\tilde{q^{kx}}
\left\{\begin{array}{c}1 \\t
\end{array}\right\}, \\
\tilde{R\!\!E}(\a)\equal \lim_{k\to \infty}
\tilde{\hbox{\ae}}(\a),\ \
&\tilde{R\!\!E}(F)=\lim_{k\to \infty}
\left\{\begin{array}{c}1 \\t^{-1}
\end{array}\right\}\,
\tilde{q^{-kx}}\tilde{\Ga_{k}}^{-1}(F) \notag
\end{align}
for the operators $\a$ and functions $F$.
The spinor constant here is
directly related to the spinor constant
$\tilde{t}^{\ 1/2}= \{t^{1/2},t^{-1/2}\}$
used in the case $|q|<1$; we will return
to this point below.

\subsubsection{\sf Definition}
For the function $G^\checkmark$ defined in
(\ref{gxlaprime}),
\begin{align}\label{spinwhitonew}
\Om^\checkmark(X,\La)\equal  \tilde{R\!\!E}(G^\checkmark)&\\
=(\tga(X)\tga(\La))^{-1}\,\Bigl(
\left\{\begin{array}{c}1 \\ 0\end{array}\right\}+
\sum_{m=1}^{\infty} q^{-m^{2}/4}\,
\bigl(\frac{\overline{E}_{-m}^\dag(\La)}
{\prod_{s=1}^{m} (1-q^{-s})}&
\left\{\begin{array}{c}X^{-m} \\ 0\end{array}\right\}\notag\\
+\,\frac{\overline{E}_{m}^\dag(\La)}
{\prod_{s=1}^{m-1} (1-q^{-s})}&
\left\{\begin{array}{c} X^{2-m} \\ X^{-m}\end{array}\right\}
\bigr)\Bigr).\notag
\end{align}

Collecting the terms with the same products
in the denominator (the norms),
\begin{align}\label{spinwhitnew}
&\Om^\checkmark(X,\La) \ =\ (\tga(X)\tga(\La))^{-1}\\
\times\sum_{m=0}^{\infty} \frac{q^{-m^{2}/4}}
{\prod_{s=1}^{m} (1-q^{-s})}&
\left\{\begin{array}{c} X^{-m}\overline{E}_{-m}^\dag(\La)+
X^{1-m}q^{-m/2-1/4}\overline{E}_{m+1}^\dag(\La)\\
X^{-m-1}q^{-m/2-1/4}\overline{E}_{m+1}^\dag(\La)\end{array}\right\}.
\notag
\end{align}

With respect to $X^{-m}$,
\begin{align}\label{spinwhitneww}
&\Om^\checkmark \ =\ (\tga(X)\tga(\La))^{-1}\Bigl(
q^{-1/4}\,\overline{E}_{1}^\dag
\left\{\begin{array}{c} X \\ 0\end{array}\right\}\\
+&\sum_{m=0}^{\infty} \frac{q^{-m^{2}/4}}
{\prod_{s=1}^{m} (1-q^{-s})}
\left\{\begin{array}{c} X^{-m}\bigl(\overline{E}_{-m}^\dag+\,
\frac{q^{-m-1}}{1-q^{-m-1}}
\overline{E}_{m+2}^\dag\bigr)\\
X^{-m}(1-q^{-m})\overline{E}_{m}^\dag\end{array}\right\}\Bigr).
\notag
\end{align}

\subsubsection{\sf Symmetrization}
We claim that the spinor-symmetrization of
$\Om^\checkmark$ is $\{\w^\checkmark,\, 0\}$ for the symmetric
global Whittaker function from
(\ref{Whitsymprime}):
\begin{align}\label{Whitsymcheck}
&\w^{\,\checkmark}(X,\La)\ =\ (\tga(X)\tga(\La))^{-1}\\
\times\sum_{m=0}^{\infty}&q^{-\frac{m^{2}}{4}}\,
\overline{P}_{m}(\La;q^{-1})\,
X^{-m}\,\prod_{s=1}^{m}\frac{1}{1-q^{-s}}\,,\notag
\end{align}

The symmetrization
is now the application of $\widehat{T}^\checkmark$;
the latter is the $T$\~symmetrizer understood
as $t^{-1}+t^{-1/2}T$ under the limit  $t\to \infty$. Since
\begin{align}\label{tcheckf}
\widehat{T}^\checkmark=\tilde{R\!\!E}(t^{-1/2}T)=
\left\{\begin{array}{c} (1-X^2)s+1\\ 0\end{array}\right\}
\end{align}
(see below),
we need to check the identity:
\begin{align*}
\overline{E}_{-m}^\dag+(1-q^{-m})\overline{E}_{m}^\dag
=\overline{P}_m^\dag,\where
\overline{P}_m^\dag=P_m\mid{}_{t\to\infty}=
\overline{P}_m\,\mid{}_{q\mapsto q^{-1}}
\end{align*}
for the Rogers polynomials $P_n\, (n\ge 0)$.
Indeed,
\begin{align*}
P_n=E_{-n}+\frac{t-tq^n}{1-tq^n}E_{n}
\hbox{\ \  due to (\ref{pviae})},
\end{align*}
which readily gives the required limit.

\subsubsection{\sf Discussion}
Generally, we can conjugate the limiting procedures
by any {\em constant} spinors
without changing its action on diagonal operators.
For instance, the Toda operator
$\tilde{R\!\!E}(\l)$ will automatically
coincide with that from (\ref{qTodaprime}). The multiplication
by  $\{1,t^{-1}\}$ is necessary to ensure the
existence of $\Om^\checkmark$ and the corresponding operators.

Theorem \ref{CONJUGG} (see also
part ($iii$) of the following theorem) clarifies the
appearance of this multiplier in full.
Let us demonstrate it. This theorem gives that
$$
\tilde{G}=t^{-1/2}\La^{-1}\pi G \and
G^\checkmark=t^{1/2}\La\pi G^\ast.
$$
The operation
$R\!\!E^\de=\hbox{lim}_{t\to 0}\hbox{\ae}^\de$
for  $\hbox{\ae}^\de=(q^{kx}\Ga_k^{-1})^\de$
from  (\ref{aede}) can be applied to $G$.
Therefore $q^{-kx}\Ga_k^{-1}$ can be applied to its
conjugation $G^\ast$ followed by the limit $t\to\infty$.
Hence the operation
\begin{align}\label{reconj}
(q^{-kx}\Ga_k^{-1})^\de(t^{-1/2}\La^{-1}\pi) =
\La^{-1}\pi\left\{\begin{array}{c}1 \\ t^{-1}\end{array}\right\}
(q^{-kx}\Ga_k^{-1})^\de
\end{align}
can be applied to $G^\checkmark$ followed by the
limit $t\to\infty$. Up to $\La^{-1}\pi$, this is exactly the twisted
(spinor) $R\!\!E$\~procedure.

Recall that here $\pi=s\Ga$,
$\pi \Ga_k^\de = \Ga_k^\de\pi$, since $\Ga^\de$ commute
with $s$ and $\Ga$, and
$$
(q^{-kx})^\de \pi\ =\
\pi\left\{\begin{array}{c}t^{1/2} \\ t^{-1/2}\end{array}\right\}
(q^{-kx})^\de.
$$
\smallskip

\subsection{\bf Main theorem}
Combining (\ref{nildahay}) and (\ref{HHbardag}) , let
$\overline{\HH}^{\,\vph\dag}$ be the span
$\C[q^{\pm 1/4}]<T^\dag,\breve{\pi}, Y^{\pm 1}>$ subject to
\begin{align}\label{HHbardagvph}
&T^\dag (T^\dag)' = 0,\ \breve{\pi}^2 = 1,\
\breve{\pi} Y \breve{\pi} =
q^{-1/2}Y^{-1},\\
&T^\dag Y^{-1}=Y(T^\dag)'
\for (T^\dag)'\equal T^\dag-1.
\notag
\end{align}
We set $X\equal\breve{\pi}(T^\dag)'$ and
$X'\equal T^\dag\breve{\pi}$; then $T^\dag X=X'(T^\dag)'$
and $(T^\dag)'X'=0=XT^\dag$.

The algebra $\,\overline{\HH}^{\,\vph\dag}\,$ is the image of
the algebra
$\,\overline{\HH}^\dag\,$ under the
anti-isomorphism $\varphi$ sending
$$
T^\dag\mapsto T^\dag,\,
\pi\mapsto\breve{\pi},\, X\mapsto Y^{-1},\,
Y\mapsto X', \, Y'\mapsto X.
$$

We will use the notations:
\begin{align}\label{omseries}
\Om_{\hbox{\rm\tiny series}}=\tga'(X)\tga'(\La)\,\Om,\ \
\Om^\checkmark_{\hbox{\rm\tiny series}}=
\tga(X)\tga(\La)\,\Om^\checkmark.
\end{align}

\begin{theorem}\label{RELIMCHECK}
(i) The operators
\begin{align*}
&\widehat{T}^\checkmark=\tilde{R\!\!E}(t^{-1/2}T),\ \ \
(\widehat{T}')^\checkmark=\tilde{R\!\!E}(t^{-1/2}T^{-1}),\ \ \,
\widehat{XT}^\checkmark=\tilde{R\!\!E}(XT),\\
&\widehat{X}^\checkmark=\tilde{R\!\!E}(t^{-1/2}X),\
(\widehat{X}')^\checkmark=\tilde{R\!\!E}(t^{-1/2}X^{-1}),\\
&\widehat{Y}^\checkmark=\tilde{R\!\!E}(Y),\ \ \ \ \ \,
(\widehat{Y}^{-1})^\checkmark=\tilde{R\!\!E}(Y^{-1}),
\end{align*}
satisfy the relations in $\overline{\HH}^{\vph\dag}$
from (\ref{HHbardagvph}).
In particular,
$$(\widehat{T}^\checkmark)^2=\widehat{T}^\checkmark,\
(\widehat{T}')^\checkmark=\widehat{T}^\checkmark-1,\
\widehat{Y}^\checkmark (\widehat{Y}^{-1})^\checkmark=1,\
(\widehat{XT}^\checkmark)^2=1.
$$

(ii)  Similar to (\ref{Omegainter}),
$\Om^\checkmark=\tilde{R\!\!E}(G^\checkmark)$
satisfies the relations
\begin{align*}
&\hat{Y}^\checkmark(\Om^\checkmark)=\La^{-1}(\Om^\checkmark),\
\hat{X}^\checkmark(\Om^\checkmark)=
\overline{Y}\,'_\La (\Om^\checkmark),\
(\hat{X}')^\checkmark\,(\Om^\checkmark)=
\overline{Y}_\La (\Om^\checkmark),\notag\\
&\hat{XT}^\checkmark(\Om^\checkmark)=\pi_\La(\Om^\checkmark),\
\hat{T}^\checkmark(\Om^\checkmark)=
\overline{T}_\La(\Om^\checkmark), \
(\hat{T}')^\checkmark(\Om^\checkmark)=
(\overline{T}')_\La(\Om^\checkmark),
\end{align*}
where $\overline{Y}\,'_\La,\, \overline{Y}_\La,\,\pi_\La,\,
\overline{T}_\La$\,
act on the argument $\La$; the other operators are
$X$\~operators.

(iii) Following Theorem \ref{CONJUGG},
let $\tilde{\Om}=(\Om^\checkmark)^\ast$, where
$(q^{1/4})^\ast=q^{-1/4}$,
$$
\ (\overline{E}_m^\dag(\La))^\ast=
\La^{-1}\overline{E}_{1-m},\ \
\{X^m,X^n\}^\ast=\{X^{-m},X^{-n}\} \for m,n\in \N.
$$
It extends the standard conjugation considered in
(\ref{astepol}) from functions to spinors.
Also, for $\pi$ and $\tau_+^{-1}(\pi)$
$=q^{1/4}X^{-1}\pi$,
we need their (natural) extensions to spinors:
\begin{align*}
\pi_{\hbox{\tiny\rm spin}}\left\{\begin{array}{c}
X^m\\X^n\end{array}\right\}=\left\{\begin{array}{c}
q^{-\frac{n}{2}}\, X^{n}\\
q^{+\frac{m}{2}}\, X^{m}\end{array}\right\},\ \,
\Pi_{\hbox{\tiny\rm spin}}\left\{\begin{array}{c}
X^m\\X^n\end{array}\right\}=\left\{\begin{array}{c}
q^{\frac{1}{4}-\frac{n}{2}}\, X^{n-1}\\
q^{\frac{1}{4}+\frac{m}{2}}\, X^{m+1}
\end{array}\right\}.
\end{align*}
Then
\begin{align}\label{conjgfnctom}
&\Om(X,\La)\ =\ \pi_{\hbox{\tiny\rm spin}}
(\La\tilde{\Om}(X,\La)),\\
&\Om(X,\La)_{\hbox{\rm\tiny series}}\ =\
\Pi_{\hbox{\tiny\rm spin}}
(\La\tilde{\Om}_{\hbox{\rm\tiny series}}(X,\La)).\notag
\end{align}
\end{theorem}

\subsubsection{\sf Proof}
The relations from $(i)$ hold by construction.
Recall that $\widehat{XT}^\checkmark$ corresponds
to the element $\breve{\pi}$ in $\overline{\HH}^{\vph\dag}$.

Similarly, the relations for
$\Om^\checkmark(X,\La)$ result from the definition of this
function as the limit of $G^\checkmark$; thus it is not
necessary to verify them. Nevertheless, it is instructional
to check them directly using the following dag-Pieri rules:
\begin{align}\label{piericheckdag}
&X\overline{E}_n^\dag\ =\ \overline{E}_{n+1}^\dag-
\overline{E}_{1-n}^\dag \ &(n>0),\\
&X\overline{E}_{-n}^\dag\ =\ (1-q^{-n})\overline{E}_{1-n}^\dag+
q^{-n}\overline{E}_{1+n}^\dag \ &(n\ge 0),\notag\\
&X^{-1}\overline{E}_{n}^\dag\ =\ (1-q^{1-n})\overline{E}_{n-1}^\dag+
\overline{E}_{1-n}^\dag \ &(n> 0),\notag \\
&X^{-1}\overline{E}_{-n}^\dag\ =\ \overline{E}_{-n-1}^\dag+
q^{-n-1}\overline{E}_{1-n}^\dag \ &(n\ge 0).\notag
\end{align}
For instance, we can rewrite (\ref{spinwhitneww})
in a convenient compact way:
\begin{align}\label{spinwhitnewww}
&\Om^\checkmark(X,\La) \ =\ (\tga(X)\tga(\La))^{-1}\Bigl(
q^{-1/4}\overline{E}^\dag_{1}(\La)
\left\{\begin{array}{c} X\\0\end{array}\right\}\\
+\sum_{m=0}^{\infty} &\frac{q^{-m^{2}/4}}
{\prod_{s=1}^{m} (1-q^{-s})}
\left\{\begin{array}{c} X^{-m}
(1-q^{-m-1})^{-1}\La\overline{E}_{-m-1}^\dag(\La)\\
X^{-m}(1-q^{-m})\overline{E}_{m}^\dag(\La)\end{array}\right\}\Bigr).
\notag
\end{align}

This formula is the best to check claim ($iii$).
Actually, this claim follows
from (\ref{reconj}), but we prefer to give a direct verification
as follows.

Upon the conjugation of
$\Om^\checkmark(X,\La)_{\hbox{\rm\tiny series}}$ and
the multiplication by $\La$:
\begin{align*}
\La\Bigl(& q^{-1/4}\overline{E}^\dag_{1}
\left\{\begin{array}{c}X\\0\end{array}\right\}\Bigr)^\ast
= q^{1/4}\left\{\begin{array}{c}X^{-1}\\0\end{array}\right\}\,,\\
\La\Bigl(&\sum_{m=0}^{\infty} \frac{q^{-m^{2}/4}}
{\prod_{s=1}^{m} (1-q^{-s})}
\left\{\begin{array}{c} X^{-m}
(1-q^{-m-1})^{-1}\La\overline{E}_{-m-1}^\dag\\
X^{-m}(1-q^{-m})\overline{E}_{m}^\dag\end{array}\right\}
\Bigr)^*\\
=&\sum_{m=0}^{\infty} \frac{q^{m^{2}/4}}
{\prod_{s=1}^{m} (1-q^{s})}
\La\left\{\begin{array}{c} X^{m}
(1-q^{m+1})^{-1}\La^{-1}(\La^{-1}\overline{E}_{2+m})\\
X^{m}(1-q^{m})(\La^{-1}\overline{E}_{1-m})\end{array}\right\}
\\
=&\left\{\begin{array}{c}(1-q)^{-1}\La^{-1}\overline{E}_{2}\\
0\end{array}\right\}+\sum_{m=1}^{\infty}
\left\{\begin{array}{c}
\frac{q^{(m+1)^{2}/4-m/2-1/4}}{\prod_{s=1}^{m+1} (1-q^{s})}
X^{m}\La^{-1}\overline{E}_{2+m}\\
\frac{q^{(m-1)^{2}/4+m/2-1/4}}{\prod_{s=1}^{m-1} (1-q^{s})}
X^{m}\overline{E}_{1-m}\end{array}\right\}.
\end{align*}
Finally, applying $\Pi_{\hbox{\tiny\rm spin}}$ and using that
$$
\Pi_{\hbox{\tiny\rm spin}}\left\{\begin{array}{c} A X^{m}\\
B X^{m}\end{array}\right\} =
\left\{\begin{array}{c} B q^{-m/2+1/4}X^{m-1}\\
A q^{m/2+1/4}X^{m+1}\end{array}\right\},
$$
where $A,B$ do not depend on $X$,
we arrive at $\Om(X,\La)_{\hbox{\rm\tiny series}}$;
claim ($iii$) is verified.

Using ($iii$), claims
from ($ii$) can be reduced to those for
$G$; see (\ref{omident}). Generally,
\begin{align*}
&A(\Om)=B_\La(\Om) \Leftrightarrow A^\vee(\Om^\checkmark)=
B_\La^\vee(\Om^\checkmark),\where\\
&A=\pi\eta(A^\vee)\pi,\
B_\La=\La\,(\eta(B^\vee))_\La\,\La^{-1},
\end{align*}
where $A,B$ are elements from the corresponding
nil-DAHA algebras, and $\eta$ is from
(\ref{etaction}). Equivalently,
\begin{align*}
&A^\vee=\pi\,\eta(A)\,\pi,\ \
B_\La^\vee=\La\,(\eta(B))_\La\,\La^{-1},
\end{align*}
since $\eta(\pi)=\pi,\eta(X)=X^{-1}$.

Here $\pi$ must be eliminated from the
formulas for $A$ and $A^\vee$, since it
does not belong to $\overline{\HH}^\vph$ or
$\overline{\HH}^{\vph\dag}$ ($\breve{\pi}$ does).
Indeed, otherwise $T$ would be invertible.
One can use $\pi$  in the
intermediate calculation (subject to all standard
identities) but the final result must not
contain it. Let us demonstrate how this works
for the basic operators.

The following identities give the required:
\begin{align*}
&\eta(Y)=\pi Y^{-1}\pi=q^{1/2}X^{-1}YX\ \Rightarrow\
Y^\vee=Y^{-1},\  (Y_\La)^\vee=q^{1/2}Y^{-1},\\
&\eta(X)=X'\ \Rightarrow\ X^\vee=q^{-1/2}X',\ \ \eta(\pi)=\pi
\ \Rightarrow\ (\pi_\La)^\vee=
\La\pi_\La\La^{-1},\\
&(\breve{\pi})^\vee\ =\ (XT)^\vee\ =\
\pi TXTY^{-1}\ =\ (\pi T)(XT)Y^{-1}\ \, =\ \,Y\breve{\pi}Y^{-1},\\
&\eta(T)=T'\,\ \Rightarrow\,\
T^\vee=\pi T'\pi=YT'Y^{-1},\ \ (T_\La)^\vee=
\La T'_\La\La^{-1}.
\end{align*}
For instance,
$Y^{-1}(\Om^\checkmark)=(\La^{-1})^\vee(\Om^\checkmark)=
\La(\Om^\checkmark)$. Multiplying by
$Y\La$, we obtain that $Y(\Om^\checkmark)$
$=\La^{-1}(\Om^\checkmark)$. As another example,
the relations $T(\Om)=T_\La(\Om)$ and
$T'(\Om)=T_\La'(\Om)$ give that
\begin{align*}
&YT'Y^{-1}(\Om^\checkmark)=
\La T'_\La \La^{-1}(\Om^\checkmark)\ \Rightarrow\
T'_{\La}(\Om^\checkmark)\\
=&YT'Y^{-1}(\Om^\checkmark)=
\La T'_\La \La^{-1}(\Om^\checkmark)\ =\ \,T'(\Om^\checkmark).
\end{align*}

It is equally possible to obtain these identities using the
limiting procedure from $G^\checkmark$ or via explicit
formulas for the basic operators.

\subsection{\bf Basic operators (large q)}
\subsubsection{\sf T-check operator}
Let us first obtain formula (\ref{tcheckf}) for
the image $\widehat{T}^\checkmark$ of the operator
\begin{align}\label{tformfull}
t^{-1/2}T=s+\frac{1-t^{-1}}{1-X^2}\,(s-1)=
\frac{X^2-t^{-1}}{X^2-1}\,s-
\frac{1-t^{-1}}{X^2-1}
\end{align}
under $\tilde{R\!\!E}$. Before taking the limit
$t\to\infty$,
\begin{align}\label{aetdiag}
&\tilde{\hbox{\ae}}(t^{-1/2}T)=\\
\left\{\begin{array}{c}1 \\t^{-1}
\end{array}\right\}
\tilde{q^{-kx}}\tilde{\Ga_{k}}^{-1}&\,
(\frac{X^2-t^{-1}}{X^2-1}\,s-
\frac{1-t^{-1}}{X^2-1})\,
\tilde{\Ga_{k}}\tilde{q^{kx}}
\left\{\begin{array}{c}1 \\t\end{array}\right\}.
\notag
\end{align}
Here $\tilde{\Ga_k},\tilde{q^{kx}}$ are
delta-spinors, so they commute with $s$.

The only non-trivial commutation relations
we need for this calculation and below are
\begin{align}\label{commutde}
&\tilde{q^{-kx}}\,\Ga\,\tilde{q^{kx}}=
\left\{\begin{array}{c}t^{1/2} \\t^{-1/2}\end{array}\right\}\Ga,\ \,
\tilde{\Ga_k^{-1}}\,X\,\tilde{\Ga_k}=
\left\{\begin{array}{c}t^{-1/2} \\t^{1/2}\end{array}\right\}X,
\notag\\
&\left\{\begin{array}{c}1 \\ t^{-1}\end{array}\right\}\,s\,
\left\{\begin{array}{c}1 \\ t\end{array}\right\}\ =\
\left\{\begin{array}{c}t \\ t^{-1}\end{array}\right\}\,s\ =\
s\,\left\{\begin{array}{c}t^{-1} \\ t \end{array}\right\},\\
&\,\Ga\,=\,
\left\{\begin{array}{c}\Ga \\ \Ga^{-1}\end{array}\right\},\ \
X\,=\,\left\{\begin{array}{c}X \\ X^{-1}\end{array}\right\},\ \
s\,=\,\left\{\begin{array}{c}s \\ s\end{array}\right\}.\notag
\end{align}
When $t\to \infty$, we readily obtain that
\begin{align}\label{tcheckff}
\widehat{T}^\checkmark=
\left\{\begin{array}{c} (1-X^2)s+1\\ 0\end{array}\right\},\
(\widehat{T}^\checkmark)^2=\widehat{T}^\checkmark,
(\widehat{T}')^\checkmark=\widehat{T}^\checkmark-1.
\end{align}
\smallskip

\subsubsection{\sf Y-check operators}
Let us present the calculation
of the check-spinor Dunkl operators
in detail. One has:
\begin{align}\label{aeydiag}
&\tilde{\hbox{\ae}}(Y)=\\
\left\{\begin{array}{c}1 \\t^{-1}
\end{array}\right\}&
\tilde{q^{-kx}}\tilde{\Ga_{k}}^{-1}\,
s\Ga(\frac{t^{1/2}X^2-t^{-1/2}}{X^2-1}\,s-
\frac{t^{1/2}-t^{-1/2}}{X^2-1})\,
\tilde{\Ga_{k}}\tilde{q^{kx}}
\left\{\begin{array}{c}1 \\t\end{array}\right\}. \notag
\end{align}

Performing the conjugations,
\begin{align}\label{aeydiagg}
\tilde{\hbox{\ae}}(Y)\hspace*{3in} & \\
=s\Ga\left\{\begin{array}{c}t^{-1/2} \\t^{1/2}\end{array}\right\}
\times \Bigl(\
\frac{t^{1/2}
\left\{\begin{array}{c}t^{-1} \\t\end{array}\right\}
\left\{\begin{array}{c}X^2 \\ X^{-2}\end{array}\right\}-
t^{-1/2}
}
{\left\{\begin{array}{c}t^{-1} \\t\end{array}\right\}
\left\{\begin{array}{c}X^2 \\ X^{-2}\end{array}\right\}
-1} \, {\left\{\begin{array}{c}t \\t^{-1}\end{array}\right\}}\, s&
\notag\\
-\frac{t^{1/2}-t^{-1/2}}
{\left\{\begin{array}{c}t^{-1} \\t\end{array}\right\}
\left\{\begin{array}{c}X^2 \\ X^{-2}\end{array}\right\}
-1}\,\Bigr)&
\notag\\
=s\Ga\Bigl(\frac{
\left\{\begin{array}{c}1 \\t\end{array}\right\}
\left\{\begin{array}{c}X^2 \\ X^{-2}\end{array}\right\}-
\left\{\begin{array}{c}1 \\ t^{-1}\end{array}\right\}
}
{\left\{\begin{array}{c}t^{-1} \\t\end{array}\right\}
\left\{\begin{array}{c}X^2 \\ X^{-2}\end{array}\right\}
-1} \, s -
\frac{
\left\{\begin{array}{c}1 \\t \end{array}\right\}-
\left\{\begin{array}{c}t^{-1} \\ 1\end{array}\right\}
}
{\left\{\begin{array}{c}t^{-1} \\t\end{array}\right\}
\left\{\begin{array}{c}X^2 \\ X^{-2}\end{array}\right\}
-1}\Bigr)&.\notag
\end{align}

\subsubsection{\sf The symmetrization}
Similarly, using that
$$
T^{-1}=T-(t^{1/2}-t^{-1/2})=
\frac{t^{1/2}X^2-t^{-1/2}}{X^2-1}\,s-
(t^{1/2}-t^{-1/2})\frac{X^2}{X^2-1},
$$
\begin{align}\label{aeyminusdiag}
&\tilde{\hbox{\ae}}(Y^{-1})=\\
\Bigl(\frac{
\left\{\begin{array}{c}1 \\t\end{array}\right\}
X^2-
\left\{\begin{array}{c}1 \\ t^{-1}\end{array}\right\}
}
{\left\{\begin{array}{c}t^{-1} \\t\end{array}\right\}
X^2
-1}& \, s -
\frac{
\left\{\begin{array}{c}1 \\t \end{array}\right\}X^2-
\left\{\begin{array}{c}t^{-1} \\ 1\end{array}\right\}X^2
}
{\left\{\begin{array}{c}t^{-1} \\t\end{array}\right\}X^2
-1}\Bigr)\,\Ga^{-1}s.\notag
\end{align}

Finally, in the limit $t\to\infty$,
\begin{align}\label{ycheck}
&\hat{Y}^\checkmark=
\Ga^{-1}\,\left\{\begin{array}{c}1 \\ 1-X^2\end{array}\right\}
+s\Ga\,\left\{\begin{array}{c}1 \\-X^2\end{array}\right\},\\
\label{yminuscheck}
&(\hat{Y}^{-1})^\checkmark=
\left\{\begin{array}{c}1-X^2 \\ 1\end{array}\right\}\,\Ga
+\left\{\begin{array}{c}X^2 \\-1\end{array}\right\}\,\Ga^{-1}s.
\end{align}
It is instructional to check that the product of these two
operators is really $1$; the formulas we give are adjusted
to make this check almost immediate.

Since $T^\dag(Y+Y^{-1})
=(Y+Y^{-1})T^\dag$ in $\overline{\HH}^{\vph\dag}$
it is granted that the operator
$\hat{Y}^\checkmark+(\hat{Y}^\checkmark)^{-1}$
preserves the space of spinors with zero second component,
i.e., those fixed by $\widehat{T}^\checkmark$.
Upon the restriction to this space,
\begin{align}\label{lspincheck}
&\hat{Y}^\checkmark+(\hat{Y}^\checkmark)^{-1}=
\left\{\begin{array}{c}\t^\checkmark \\ 0\end{array}\right\}
\for \t^\checkmark \ =\ (1-X^2)\Ga+\Ga^{-1}.
\end{align}
The latter is exactly the check-$q$\~Toda operator from
(\ref{qTodaprime}).

\subsubsection{\sf Pi-check operator}
To complete the list of basic operators, let us
provide the formula $\hat{\pi}^\checkmark$, which is
$\tilde{R\!\!E}(XT)$.
It is the image of
$\breve{\pi}$ in the
hat-dag-polynomial representation
of $\overline{\HH}^{\vph\dag}$ we are considering now.
One has:
\begin{align}\label{pihatdag}
\hat{\pi}^\checkmark =
\left\{\begin{array}{c} X \\ X^{-1}\end{array}\right\}&
\left\{\begin{array}{c} (1-X^2)s+1 \\ s-X^2\end{array}\right\}=
\left\{\begin{array}{c} X(1-X^2) \\ X^{-1}\end{array}\right\}s
+\left\{\begin{array}{c} X \\ -X\end{array}\right\},\notag\\
&\hat{\pi}^\checkmark
(\left\{\begin{array}{c} f_1 \\ f_2\end{array}\right\})=
\left\{\begin{array}{c} X(1-X^2)f_2+Xf_1 \\
X^{-1}f_1-X f_2\end{array}\right\}.
\end{align}
Recall that $\{As,Bs\}(\{f_1,f_2\})=
\{A(f_2),B(f_1)\}$, i.e., $s$ at the end of the component
means that the remaining operator will be applied to the
other component.

It is instructional to check directly that
$(\hat{\pi}^\checkmark)^2=1$ directly using the
last formula from (\ref{pihatdag}). We
we will leave it as an exercise.

It is interesting to examine the operator interpretation
of the relations
between $\pi^\checkmark, T^\checkmark$
and $X^\checkmark$. The following must hold due to the
general theory:

\begin{align}\label{xhatdagg}
(\hat{\pi}^\checkmark)(\hat{T}')^\checkmark=\hat{X}^\checkmark
\equal \tilde{R\!\!E}(t^{-1/2}X)=
\left\{\begin{array}{c} 0 \\ X^{-1}\end{array}\right\},\\
(\hat{T}^\checkmark)(\hat{\pi}^\checkmark)=(\hat{X}')^\checkmark
\equal \tilde{R\!\!E}(t^{-1/2}X^{-1})=
\left\{\begin{array}{c} X^{-1} \\ 0\end{array}\right\}.\notag
\end{align}

Here the calculation of operators $\hat{X}^\checkmark,\
(\hat{X}')^\checkmark$ themselves is simple; use
(\ref{commutde}) and the fact that diagonal constant
spinors commute with any diagonal operators. A direct verification
of (\ref{xhatdagg}) is a bit more
involved. For instance, let us check the second formula.
The calculation goes as follows:

\begin{align*}
&(\hat{T}^\checkmark)(\hat{\pi}^\checkmark)\\
=&\Bigl(
\left\{\begin{array}{c} 1-X^2\\ 0\end{array}\right\}\,s
+\left\{\begin{array}{c} 1 \\ 0\end{array}\right\}\Bigr)
\Bigl(
s\,\left\{\begin{array}{c} X^{-1}\\ X(1-X^2)\end{array}\right\}
+\left\{\begin{array}{c} X \\ -X\end{array}\right\}\Bigr)\\
=&
\left\{\begin{array}{c} 1-X^2\\ 0\end{array}\right\}
\left\{\begin{array}{c} X^{-1} \\ X(1-X^2)\end{array}\right\}
+\left\{\begin{array}{c} 1\\ 0\end{array}\right\}
\left\{\begin{array}{c} X(1-X^2) \\ X^{-1}\end{array}\right\}\,s\\
&+\left\{\begin{array}{c} 1\\ 0\end{array}\right\}
\left\{\begin{array}{c} X\\ -X\end{array}\right\}
+\left\{\begin{array}{c} 1-X^2\\ 0\end{array}\right\}
\left\{\begin{array}{c} -X \\ X\end{array}\right\}\,s\ =\
\left\{\begin{array}{c} X^{-1} \\ 0\end{array}\right\}.
\end{align*}

\subsection{\bf The eigenvalue problem}
It is interesting to understand how the symmetries
of $\Om^\checkmark(X,\La)$ can be obtained directly,
without using the general theory. For arbitrary root
systems no direct methods are known. For $A_n$, the Pieri
formulas are essentially sufficient, although the
calculations are very involved for $n>1$. Even for
$n=1$ it is not immediate.

\subsubsection{\sf Preparations}
Let us consider the main
symmetry, which is the eigenvalue problem. It reads:
$\hat{Y}^\checkmark(\Om^\checkmark)=
\La^{-1}\Om^\checkmark$. We will leave the other
relations from part ($ii$) of Theorem \ref{RELIMCHECK}
as an exercise.

Let use the formula for $\hat{Y}^\checkmark$,
\begin{align}\label{ymincheck}
\hat{Y}^\checkmark&=
\left\{\begin{array}{c}1\\ 1-qX^2\end{array}\right\}\,\Ga^{-1}
+\left\{\begin{array}{c}-q^{-1}X^2 \\1\end{array}\right\}\,
\Ga^{-1}s,
\notag\\
\hat{Y}^\checkmark
(\left\{\begin{array}{c} f_1 \\ f_2\end{array}\right\})&=
\left\{\begin{array}{c}\Ga^{-1}(f_1)-q^{-1}X^2\Ga^{-1}(f_2)
\\ (1-qX^2)\Ga(f_2)+\Ga(f_1).
\end{array}\right\}\,
\end{align}
We will need both presentations
of $\Om^\checkmark(X,\La)$,
namely, that from (\ref{spinwhitnewww}):
\begin{align}\label{spinwhitnewa}
&\Om^\checkmark(X,\La) \ =\ (\tga(X)\tga(\La))^{-1}\Bigl(
q^{-1/4}\overline{E}^\dag_{1}(\La)
\left\{\begin{array}{c} X\\0\end{array}\right\}\\
+\sum_{m=0}^{\infty} &\frac{q^{-m^{2}/4}}
{\prod_{s=1}^{m} (1-q^{-s})}
\left\{\begin{array}{c} X^{-m}
(1-q^{-m-1})^{-1}\La\overline{E}_{-m-1}^\dag(\La)\\
X^{-m}(1-q^{-m})\overline{E}_{m}^\dag(\La)\end{array}\right\}\Bigr).
\notag
\end{align}
and that from (\ref{spinwhitneww}):
\begin{align}\label{spinwhitnewb}
&\Om^\checkmark \ =\ (\tga(X)\tga(\La))^{-1}\Bigl(
q^{-1/4}\,\overline{E}_{1}^\dag
\left\{\begin{array}{c} X \\ 0\end{array}\right\}\\
+&\sum_{m=0}^{\infty} \frac{q^{-m^{2}/4}}
{\prod_{s=1}^{m} (1-q^{-s})}
\left\{\begin{array}{c} X^{-m}\bigl(\overline{E}_{-m}^\dag+\,
\frac{q^{-m-1}}{1-q^{-m-1}}
\overline{E}_{m+2}^\dag\bigr)\\
X^{-m}(1-q^{-m})\overline{E}_{m}^\dag\end{array}\right\}\Bigr).
\notag
\end{align}
Recall that
$\Ga(\tga(X)^{-1})=q^{-1/4}X^{-1}\tga(X)^{-1}$ and
$$
X=\{X,X^{-1}\},\ \
\Ga^{-1}(\{X^{-m},X^{-m}\})= \{q^{m/2}X^{-m},q^{-m/2}X^{-m}\}.
$$

The following notations will be convenient
\begin{align}\label{epolnormal}
\overline{E}_m^\ddag=\overline{E}_m^\dag
\prod_{s=1}^{|\tilde{n}|-1}(1-q^{-s})^{-1}\,,\ \,
\left\{\begin{array}{c}
|\tilde{n}|=n \for n>0\\|\tilde{n}|=n+1 \hbox{\, otherwise}
\end{array}\right\};
\end{align}
see (\ref{tilden}).

\subsubsection{\sf Explicit verification}
Applying (\ref{ymincheck}) to (\ref{spinwhitnewb}),
one obtains:
\begin{align}\label{spinwhitnec}
\hat{Y}^\checkmark(\Om^\checkmark) =
(\tga(X)\tga(\La))^{-1}\ \ \ \ \ \ \ \ \ \ \ \ \ \ \
\ \ \ \ \ \ \ \ \ &\notag\\
\times\Bigl(
\overline{E}^\ddag_{1}
\left\{\begin{array}{c}q^{-1}X^2\\
1\end{array}\right\}
-\left\{\begin{array}{c}0\\
q^{-1/4}(X^{-1}-qX)\end{array}\right\}&\\
+\sum_{m=0}^{\infty}
\left\{\begin{array}{c}
q^{-(m-1)^2/4}X^{1-m}(
\overline{E}_{-m}^\ddag+q^{-m-1}\overline{E}_{m+2}^\ddag
-q^{-1}X^2\overline{E}_{m}^\ddag)
\\
q^{-(m+1)^2/4}X^{-1-m}(
(1-qX^2)\overline{E}_{m}^\ddag+\overline{E}_{-m}^\ddag
+q^{-m-1}\overline{E}_{m+2}^\ddag)
\end{array}\right\}&\notag\Bigr).
\end{align}

Until the end of this calculation,
we will ignore certain terms involving powers
of $X$ greater than or equal to $-1$ (but
not all such terms).
Correspondingly, we will use $\approx$
instead of $=$.

Since
$q^{-(m-1)^2/4}q^{-m}\overline{E}_{m+2}^\ddag=
q^{-(m+1)^2/4}\overline{E}_{m+2}^\ddag$,
only the terms $\overline{E}_{-m}^\ddag$ will
contribute to the first component; the other terms
will cancel each other in the final summation.

The same kind of cancelation will occur in the second
component. Only $\overline{E}_{m}^\ddag$
and $\overline{E}_{-m}^\ddag$ will really contribute
to the final summation.

Thus,
\begin{align}\label{spinwhitnecc}
&\hat{Y}^\checkmark(\Om^\checkmark)\ \approx\ \\
(\tga(X)\tga(\La))^{-1}\,
&\sum_{m=0}^{\infty}\,
\left\{\begin{array}{c}
q^{-(m-1)^2/4}X^{1-m}\,
\overline{E}_{-m}^\ddag
\\
q^{-(m+1)^2/4}X^{-1-m}(
\overline{E}_{m}^\ddag+\overline{E}_{-m}^\ddag)
\end{array}\right\}\,.\notag
\end{align}

However, this expression coincides with
\begin{align}\label{spinwhitnewaa}
&\La^{-1}\Om^\checkmark(X,\La) \ \approx\ \\
(\tga(X)\tga(\La))^{-1}\,
&\sum_{m=0}^{\infty} q^{-m^{2}/4}
\left\{\begin{array}{c} X^{-m}
\overline{E}_{-m-1}^\ddag(\La)\\
\,X^{-m}\La^{-1}\overline{E}_{m}^\ddag(\La)
\end{array}\right\}
\notag
\end{align}
due to the Pieri formula
$$
\La^{-1}\overline{E}^\ddag_m=
\overline{E}^\ddag_{m-1}+\overline{E}^\ddag_{1-m}.
$$
It is of course with the reservation that we
ignore certain powers of $X$. It holds
without this reservation, which is easy to check.
Generally, controlling
the boundary terms in this and similar calculations
is a nontrivial combinatorial problem even for $A_n\,(n>1)$.
The $q,t$\~setting and spherical polynomials manage this
difficulty automatically.

\setcounter{equation}{0}
\section{\sc Appendix: Givental-Lee theory}
\label{sec:qKth}
\subsection{\sf The J-series}
Let $G$ be a complex semi\-simple algebraic group and let
$\mathfrak{B}=G/B$ be the flag variety, where $B$ is a Borel
subgroup. For $G=SL(n+1)$, Givental and Lee show in \cite{GiL}
that quantum $K$\~theory on $\mathfrak{B}$ produces
formal series solutions to the $q$\~Toda eigenvalue problem.
A conjectured generalization to arbitrary semisimple $G$ is stated
in \cite{GiL}, and for simply-laced $G$ this has been proved by
Braverman and Finkelberg \cite{BF}.

Let us consider the case of $G=SL(n+1)$ in more detail and establish
the connection to the $q$\~Whittaker functions considered in this paper.
In the notation of \cite{GiL}, the $q$\~Toda difference operator
is given by
\begin{align}
\widehat{H}&=q^{\partial/\partial t_0}+
q^{\partial/\partial t_1}(1-Q_1)+\cdots+
q^{\partial/\partial t_n}(1-Q_n),
\end{align}
where $Q_i=e^{t_{i-1}-t_i}$ act as multiplication operators and
$q^{\partial/\partial t_j} : t_i \mapsto t_i+\delta_{ij} \ln q$
are translations. We write $Q=(Q_1,\ldots,Q_n)$.
A $q$\~Whittaker function is a solution $I=I(Q,\Lambda,q)$
to the eigenvalue problem
\begin{align}
\widehat{H}I&=(\Lambda_0^{-1}+\cdots+\Lambda_n^{-1})I.
\end{align}
Here $\Lambda=(\Lambda_0,\cdots,\Lambda_n)$ are coordinates on the
maximal torus of diagonal matrices, i.e.,
$\Lambda_0\cdots\Lambda_n=1$.
\smallskip

The $q$\~Whittaker function due to Givental and Lee takes the form 
$I=p^{\ln Q/\ln q}J$, where $J$ is a formal power series
in the variables $Q_i$ with coefficients in
$K_G(\mathfrak{B})\otimes\mathbb{Q}(\Lambda,q)$.
Here $p=(p_1,\ldots,p_n)$, and the $p_i$ are elements of
$K_G(\mathfrak{B})$ represented by the pullbacks of the
Hopf line bundles over the projective factors involved
in the Pl\"{u}cker embedding of $\mathfrak{B}$.
The power series $J$ is the so-called equivariant $K$\~theoretic
$J$\~function of $\mathfrak{B}$; it is a generating function
for the $K$\~theoretic one-point genus zero Gromov-Witten invariants
with descendants.

In the case $n=1$, which is the main focus of this paper,
the $J$\~function from \cite{GiL} is as follows:
\begin{align}
J&=\sum_{d=0}^\infty\frac{e^{d(t_0-t_1)}}{\prod_{m=1}^d
(1-p\Lambda_0 q^m)(1-p\Lambda_0^{-1}q^m)}.
\label{A1Jfunction}
\end{align}
In $K_G(\mathfrak{B})$, the Hopf line bundle $p$ satisfies
\begin{align}
&(1-p\Lambda_0)(1-p\Lambda_0^{-1})=0;
\label{hopfrel}
\end{align}
the factor $\e=p^{\ln Q/\ln q}$ 
is understood formally
as a solution of the system
\begin{align}
&q^{\partial/\partial t_0}\e=p \e,\ \
q^{\partial/\partial t_1}\e=p^{-1}\e.
\label{ueqn}
\end{align}
We mention that (\ref{A1Jfunction}) and its $SL(3)$\~generalization
from \cite{GiL} are extended to certain explicit formulas to the
case of Grassmannians $SL(n+1)/P$ in \cite{T}.

\subsection{\sf The connection}
The $q$\~Whittaker function $\w(X,\La)$ is given in (\ref{Whitsym})
as a convergent series expressed in terms of the
$q$\~Hermite polynomials. This is the $B$\~model interpretation of
$\w$, as described in the Introduction. 
We now explain how the $J$\~function enters into the Harish-Chandra
asymptotic decomposition (\ref{Whitexpan}), thereby establishing
the $A$\~model interpretation of $\w$.

In view of (\ref{hopfrel}), we consider the
$J$\~function upon the substitution
$p=\Lambda_0$ or $p=\Lambda_0^{-1}$.
Let us write $\La=\La_0$ and consider only the first substitution,
as the latter is simply its conjugate under $s=s_1$ in 
the Weyl group $W$.
Set $\tilde{J}(X,\La)=J(qX^{-2},\La,q)|_{p\mapsto\La_0}$.

Thus the variable $Q_1=e^{t_0-t_1}$ is replaced by $qX^{-2}$,
where $X=q^x$.
The operators $q^{\partial/\partial t_1}$
and $q^{\partial/\partial t_0}$ are replaced by $\Gamma$ and
$\Gamma^{-1}$, respectively; recall that 
$\Gamma(f(X))=f(q^{1/2}X)$. The $q$\~Toda operator is then 
given by $\cT=(1-X^{-2})\Gamma+\Gamma^{-1}$.

Using $\Ga$, the relations (\ref{ueqn}) become
$\Gamma(\e)=\Lambda^{-1}\e$.
In terms of the Gaussian $\tga'$ from
Section \ref{sec:innprod}, we take
$
\e(X,\La)=\tga'(X\La)/\tga'(X)\tga'(\La).
$ 
It provides the leading term of the Harish-Chandra 
expansion of our global $q$\~Whittaker function
$\w(X,\La)$ from
(\ref{Whitexpan}). In terms of $\tilde{J}$,
the global function reads as follows.
For $|q|<1$ and $|X|>|q|^{1/2}$, 
\begin{align}\label{introwhitexpan}
&\w(X,\La)=\lan\overline{\mu}\ran
\sum_{w\in W}\overline{\si}(w(\La^{-1}))\,
\e(X,w(\La))\,\tilde{J}(X,w(\La)).
\end{align}

Note that our global function is the weighted sum of
$|W|$ asymptotic solutions (here $W=\S_2$), corresponding
to different choices of $p$ in the Givental-Lee
approach. However such sums do not appear in \cite{GiL} and
there is no presence there of the $q$\~Whittaker version of
the celebrated Harish-Chandra $c$\~function,
which is $\si$ in (\ref{introwhitexpan}). The latter is the
key in the Harish-Chandra theory and its (recent) 
$q$\~generalization.

Formula (\ref{introwhitexpan}) is a limiting case of the 
Harish-Chandra type asymptotic
decomposition of the global hypergeometric function 
(depending on $q,t$) for $A_1$; its existence for arbitrary root 
systems is due to \cite{Sto2}.
\medskip

\renewcommand\refname{\sc{References}}
\bibliographystyle{unsrt}

\begin{thebibliography} {ABCD}
\bibitem [BF] {BF}
{A.~Braverman}, and {M. Finkelberg},
{\em Semi-infinite Schubert varieties and quantum K-theory of
flag manifolds},
\ \ Preprint\ arXiv:1111.2266 (2011).

\bibitem [Ch1] {C101}
{I.~Cherednik},
{\em Double affine Hecke algebras},
London Mathematical Society Lecture
Note Series, { 319}, Cambridge University Press, Cambridge, 2006.

\bibitem [Ch2] {C1}
\bysame,
{\em Intertwining operators of double affine Hecke algebras},
Se\-lec\-ta Math. New ser. { 3} (1997), 459--495

\bibitem [Ch3] {C5}
\bysame,
{\em Difference Macdonald-Mehta conjecture},
IMRN { 10} (1997), 449--467.

\bibitem [Ch4] {ChL}
\bysame,
{\em Toward Harmonic Analysis on DAHA
(Integral formulas for canonical traces)},
Notes of the lecture delivered at University
of Amsterdam (May 30, 2008);
http://math.mit.edu/~etingof/hadaha.pdf

\bibitem [Ch5] {ChW}
\bysame,
{\em Whittaker limits of difference spherical functions}, 
IMRN, { 20} (2009), 3793--3842; arXiv:0807.2155 (2008).

\bibitem [Ch6] {ChA}
 {\em Affine extensions of Knizhnik-Zamolodchikov equations and
 Lusztig's isomorphisms}, in: Special functions.
Proceedings of the Hayashibara forum 1990, Okayama, 63-77,
Springer-Verlag (1991).

\bibitem [ChM] {C102}
\bysame, and {X.~Ma},
{\em Spherical and Whittaker functions via DAHA II
(Spinor Whittaker and Bessel functions in rank one)},
Selecta Mathematica (2012); Preprint arXiv: 0904.4324 (2009).

\bibitem [Et] {Et1}
{P.~Etingof},
{\em  Whittaker functions on quantum groups and $q$-deformed Toda
operators}, AMS Transl. Ser. 2, { 194}, 9--25, AMS, Providence,
Rhode Island, 1999.

\bibitem [FJM] {FJM}
{B.~Feigin}, {E.~Feigin}, {M.~Jimbo}, {T.~Miwa}, and {E.~Mukhin},
{\em Fermionic formulas for eigenfunctions of the difference
Toda Hamiltonian},
Lett. Math. Phys., {88} (2009), 39--77.

\bibitem [GLO] {GLO}
{A.~Gerasimov}, {D.~Lebedev}, and {S.~Oblezin},
{\em On $q$-deformed $\mathfrak{gl}_{\ell+1}$-Whittaker function III},
Lett. Math. Phys., {97} (2011), 1--24.

\bibitem [GiL] {GiL}
{A.~Givental}, and {Y.-P.~Lee},
{\em Quantum $K$-theory on flag
manifolds, finite-difference Toda lattices and quantum groups},
Invent. Math., { 151} (2003), 193--219.

\bibitem [GW] {GW}
{R.~Goodman}, and {N.~R.~Wallach},
{\em Conical vectors and Whittaker vectors},
J. Functional Analysis, { 39} (1980), 199--279.

\bibitem [HO] {HO}
{G.J.~Heckman}, and  {E.M.~Opdam},
{\em Root systems and hypergeometric functions I},
Comp. Math. { 64} (1987), 329--352.

\bibitem [Ion] {Ion}
{B.~Ion},
{\em Nonsymmetric Macdonald polynomials and Demazure characters}, 
Duke Mathematical Journal { 116}:2 (2003), 299--318.

\bibitem [KK] {KK}
{B.~Kostant}, and {S.~Kumar},
{\em T-Equivariant K-theory of generalized flag varieties},
J. Diff. Geometry {  32} (1990), 549--603.

\bibitem [Lub] {Lub}
{D.~Lubinsky},
{\em On q-exponential functions for $|q| = 1$},
Canad. Math. Bull. { 41}:1 (1998), 86--97.

\bibitem [Ma] {Ma4}
{I.~Macdonald},
{\em Symmetric Functions and Hall Polynomials}, Second Edition,
Oxford University Press (1999).

\bibitem [Me] {Me}
{M.~van Meer},
{\em
Bispectral quantum Knizhnik-Zamolodchikov equations for
arbitrary root systems},
Selecta Math., {17} (2011), 183-221.

\bibitem [MS] {MS}
\bysame, and {J.V.~Stokman},
{\em Double affine Hecke algebras and bispectral quantum Knizhnik-
Zamolodchikov equations}, Int. Math. Res. Not., { 6},
(2010) 969--1040.

\bibitem [O] {O5}
{E.~Opdam},
{\em On the spectral decomposition of affine
Hecke algebras}, J. Inst. Math. Jussieu (2004).

\bibitem [Rui] {Rui}
{S.~Ruijsenaars},
{\em Systems of Calogero-Moser type},
in: Proceedings of the 1994 Banff summer school
``Particles and fields", CRM series in mathematical physics,
(G. Semenoff, L. Vinet, Eds.), 251--352, Springer,
New York, 1999.

\bibitem [San] {San}
{Y.~Sanderson},
{\em On the Connection Between Macdonald Polynomials and
Demazure Characters},
J. of Algebraic Combinatorics, { 11} (2000), 269--275.

\bibitem [Sto1] {Sto}
{J.~Stokman},
{\em Difference Fourier transforms for nonreduced
root systems}, Sel. math., New ser. { 9} (2003) 409--494.

\bibitem [Sto2] {Sto2}
{J.~Stokman}
{\em The c-function expansion of a basic hypergeometric function
associated to root systems},
\ \ Preprint \ arXiv:1109.0613 (2011).

\bibitem [Sus] {Sus}
{S.~Suslov},
{\em Another addition theorem for the $q$-exponential function},
J. Phys. A: Math. Gen. { 33}: 41 (2000) L375-L380.

\bibitem [T] {T}
{K.~Taipale},
{\em K-theoretic J-functions of type A flag varieties},
\ \ Preprint \ arXiv:1110.3117 (2011).

\bibitem [Wa] {Wa}
{N.~R.~Wallach},
{\em Real Reductive Groups II},
Academic Press, Boston, 1992.

\end{thebibliography}

\end{document}